\documentclass[ms,nonblindrev]{informs3-VVM}

\OneAndAHalfSpacedXI % current default line spacing
\usepackage{ctable}
\usepackage{siunitx}

\usepackage{color}

\usepackage{algorithm}
\usepackage{algorithmic}

\usepackage{amsmath}
\usepackage{amssymb}

\usepackage{natbib}
\bibpunct[, ]{(}{)}{,}{a}{}{,}%
\def\bibfont{\small}%
\def\ab{\mathbf{a}}
\def\bb{\mathbf{b}}
\def\cb{\mathbf{c}}
\def\db{\mathbf{d}}
\def\eb{\mathbf{e}}
\def\pb{\mathbf{p}}
\def\qb{\mathbf{q}}
\def\sbee{\mathbf{s}}

\def\ub{\mathbf{u}}
\def\wb{\mathbf{w}}
\def\xb{\mathbf{x}}
\def\yb{\mathbf{y}}
\def\zb{\mathbf{z}}

\def\Ab{\mathbf{A}}
\def\Fb{\mathbf{F}}
\def\Mb{\mathbf{M}}

\def\alphab{\boldsymbol \alpha}
\def\betab{\boldsymbol \beta}
\def\gammab{\boldsymbol \gamma}
\def\lambdab{\boldsymbol \lambda}
\def\mub{\boldsymbol \mu}
\def\pib{\boldsymbol \pi}
\def\rhob{\boldsymbol \rho}

\def\xib{\boldsymbol \xi}
\def\etab{\boldsymbol \eta}

\def\zerob{\mathbf{0}}
\def\oneb{\mathbf{1}}

\def\Ebb{\mathbb{E}}
\def\Ibb{\mathbb{I}}

\def\Rbb{\mathbb{R}}

\def\Dcal{\mathcal{D}}
\def\Fcal{\mathcal{F}}
\def\Kcal{\mathcal{K}}
\def\Ncal{\mathcal{N}}
\def\Ocal{\mathcal{O}}
\def\Pcal{\mathcal{P}}
\def\Qcal{\mathcal{Q}}
\def\Rcal{\mathcal{R}}
\def\Scal{\mathcal{S}}
\def\Ucal{\mathcal{U}}
\def\Wcal{\mathcal{W}}
\def\Xcal{\mathcal{X}}

\def\DR{\mathrm{DR}}
\def\RR{\mathrm{RR}}
\def\Nom{\mathrm{N}}
\def\Hybrid{\mathrm{H}}
\def\UB{\mathrm{UB}}
\def\LB{\mathrm{LB}}

\def\WC{\mathrm{WC}}
\def\Uniform{\mathrm{Uniform}}

\def\Halmos{$\square$}
\def\CONV{\mathsf{CONV}}
\def\floor{\text{floor}}

\def\val{\mathsf{val}}
\def\train{\mathsf{train}}
\def\test{\mathsf{test}}

\def\PrimalCG{ \textsc{PrimalCG}\xspace}
\def\DualCG{ \textsc{DualCG}\xspace}

\def\RI{\mathrm{RI}}

\usepackage{xspace}

\def\orangeJuice{{\tt orangeJuice}\xspace}

\newcommand{\numvvm}[1]{ {\num[round-precision=3,round-mode=figures,scientific-notation=true]{#1}}}

\newenvironment{proofvvm}{\paragraph{Proof:}}{\\[1em]}

\usepackage{enumitem}

\RUNTITLE{Randomized Robust Price Optimization}

\TITLE{Randomized Robust Price Optimization} 
\ARTICLEAUTHORS{%
\AUTHOR{Xinyi Guan}
\AFF{Department of Logistics and Maritime Studies, Faculty of Business, The Hong Kong Polytechnic University, Hong Kong, \EMAIL{xinyi.guan@polyu.edu.hk}}
\AUTHOR{Velibor V. Mi\v{s}i\'{c}}
\AFF{UCLA Anderson School of Management, University of California, Los Angeles, California 90095, United States, \EMAIL{velibor.misic@anderson.ucla.edu}} %\textsc{Preliminary -- do not distribute}
}

\ABSTRACT{
The robust multi-product pricing problem is to determine the prices of a collection of products so as to maximize the worst-case revenue, where the worst case is taken over an uncertainty set of demand models that the firm expects could be realized in practice. A tacit assumption in this approach is that the pricing decision is a deterministic decision: the prices of the products are fixed and do not vary. In this paper, we consider a \emph{randomized} approach to robust pricing, where a decision maker specifies a distribution over potential price vectors so as to maximize its worst-case revenue over an uncertainty set of demand models. We formally define this problem -- the \emph{randomized robust price optimization} problem -- and analyze when a randomized price scheme performs as well as a deterministic scheme versus when it yields a benefit. We also propose solution methods for obtaining an optimal randomization scheme over a discrete set of candidate price vectors and show how these methods are applicable for common demand models, such as the linear, semi-log and log-log demand models. We numerically compare the randomized and deterministic approaches on a variety of synthetic and real problem instances; on instances derived from a real grocery retail scanner dataset, we show that the improvement in worst-case revenue can be as high as 92\%. Using the same grocery retail scanner dataset, we also show that the randomized approach can produce price prescriptions that achieve higher out-of-sample revenue than the nominal and deterministic robust approaches.
}

\TheoremsNumberedThrough     % Preferred (Theorem 1, Lemma 1, Theorem 2)
\ECRepeatTheorems

\EquationsNumberedThrough    % Default: (1), (2), ...
\MANUSCRIPTNO{\quad}


\begin{document}

\maketitle

\section{Introduction}

Price optimization is a key problem in modern business. The price optimization problem can be stated as follows: we are given a collection of products. We are given a demand model which tells us, for each product, what the expected demand for that product will be as a function of the price of that product as well as the price of the other products. Given this demand model, the price optimization problem is to decide a price vector -- i.e., what price to set for each product -- so as to maximize the total expected revenue arising from the collection of products.

%{\comebacktothis (Q: how much incremental revenue can be attributed to better pricing? is there a McKinsey report on this?)}

The primary input to a price optimization approach is a demand model, which maps the price vector to the vector of expected demands of a product. However, in practice, the demand model is never known exactly, and must be estimated from data. This poses a challenge because data is typically limited, and thus a firm often faces uncertainty as to what the demand model is. This is problematic because a mismatch between the demand model used for price optimization -- the nominal demand model -- and the demand model that materializes in reality can lead to suboptimal revenues. 

As a result, there has been much research in how to address demand model uncertainty in pricing. In the operations research community, a general framework for dealing with uncertainty is robust optimization. The idea of robust optimization is to select an uncertainty set, which is a set of values for the uncertain parameter that we believe could plausibly occur, and to optimize the worst-case value of the objective function, where the worst-case is taken over the uncertainty set. In the price optimization context, one would construct an uncertainty set of potential demand models and determine the prices that maximize the worst-case expected revenue, where the worst-case is the minimum revenue over all of the demand models in the uncertainty set. In applying such a procedure, one can ensure that the performance of the chosen price vector is good under a multitude of demand models, and that one does not experience the deterioration of a price vector optimized from a single nominal demand model.

Typically in robust optimization, one seeks the \emph{single} best decision that optimizes the worst-case value of the objective function; stated differently, one \emph{deterministically} implements a single decision. However, a recent line of research \citep{delage2019dice} has revealed that with regard to the worst-case objective, it is possible to obtain better performance than the traditional deterministic robust optimization approach by \emph{randomizing} over multiple solutions. Specifically, instead of selecting a single decision in some feasible set to optimize the worst-case objective, one selects a \emph{distribution} supported on the feasible set that informs the decision maker how to randomize. 

In this work, we propose a methodology for robust price optimization that is based on randomization. In particular, we propose solving a randomized robust price optimization (RRPO) problem, which outputs a probability distribution that specifies the frequency with which the firm should use different price vectors. From a practical perspective, such a randomization scheme has the potential to be implemented in modern retailing as a strategy for mitigating demand uncertainty. In particular, in an ecommerce setting, randomization is already used for A/B testing, which involves randomly assigning some customers to one experimental condition and other customers to a different experimental condition. Thus, it is plausible that the same form of randomization could be used to display different price vectors with certain frequencies. In the brick-and-mortar setting, one can potentially implement randomization by varying prices geographically or temporally. For example, if the RRPO solution is the three price vectors $\pb_1, \pb_2, \pb_3$ with probabilities 0.2, 0.3, 0.5, then for a set of 50 regions, we would choose 10 ($ = 50 \times 0.2$) to assign to price vector $\pb_1$, 15 ($= 50 \times 0.3$) to assign to $\pb_2$, and 25 ($ = 50 \times 0.5$) to assign to $\pb_3$. Similarly, if one were to implement the same randomization scheme temporally, then for a selling horizon of 20 weeks, one would use the price vector $\pb_1$ for 4 ($ = 20 \times 0.2$) weeks, $\pb_2$ for 6 ($ = 20 \times 0.3$) weeks and $\pb_3$ for 10 ($ = 20 \times 0.5$) weeks. 

We make the following specific contributions:

\begin{enumerate}
\item \textbf{Benefits of randomization}. We formally define the RRPO problem and analyze it under different conditions to determine when the underlying robust price optimization problem is \emph{randomization-receptive} -- there is a benefit from implementing a randomized decision over a deterministic decision -- versus when it is \emph{randomization-proof} -- the optimal randomized and deterministic decisions perform equally well. We show that the robust price optimization problem is randomization-proof in several interesting settings, which can be described roughly as follows: (1) when the set of feasible price vectors is convex and the set of uncertain revenue functions is concave; (2) when the set of feasible price vectors and the set of uncertain demand parameters are convex, and the revenue function obeys certain quasiconvexity and quasiconcavity properties with respect to the price vector and the uncertain parameter vector; and (3) when the set of feasible price vectors is finite, and a certain strong duality/saddle point property holds. We showcase a number of examples of special cases that satisfy the hypotheses of these results and consequently are randomization-proof. We also present several examples showing how these results can fail to hold when certain assumptions are relaxed and the problem thus becomes randomization-receptive.
\item \textbf{Tractable solution algorithms}. We propose algorithms for solving the RRPO problem in two different settings: 
\begin{enumerate}
\item In the first setting, we assume that the set of possible price vectors is a finite set and that the uncertainty set of demand function parameters is a convex set. In this setting, when the revenue function is convex in the uncertain parameters, we show that the RRPO problem can be solved via delayed constraint generation. The separation problem that is solved to determine which constraint to add is exactly the nominal pricing problem for a fixed uncertain parameter vector of the demand function. For the log-log and semi-log demand models, we show that this nominal pricing problem can be reformulated and solved to global optimality as a mixed-integer exponential cone program. We believe these reformulations are of independent interest as to the best of our knowledge, these are the first exact mixed-integer convex formulations for these problems in either the marketing or operations literatures, and they leverage recent advances in solution technology for mixed-integer conic programs (as exemplified in the conic solver Mosek). In this setting, we also prove that the complexity of randomization, which is the size of the support of the randomized price vector, is at most $M+1$, where $M$ is the dimension of the uncertain parameter vector. Lastly, we also show how our method in this setting can be extended to incorporate moment constraints on the price distribution.
\item In the second setting, we assume that both the price set and the uncertainty set are finite sets. In this setting, we show how the RRPO problem can be solved using a double column generation method, which involves iteratively generating new uncertainty realizations and price vectors by solving primal and dual separation problems, respectively. We show how the primal and dual separation problems can be solved exactly for the linear, semi-log and log-log demand models. 
\end{enumerate}
\item \textbf{Numerical evaluation with synthetic and real data}. We evaluate the effectiveness of randomized pricing on different problem instances generated synthetically and problem instances calibrated with real data. We show that randomized pricing can significantly improve worst-case revenues over deterministic pricing, with the benefit in our real data instances being as high as 92\%. Additionally, we show that for instances of realistic size (up to 20 products), our algorithm can solve the RRPO problem in a reasonable amount of time (no more than two minutes on average). We also consider a data-driven experiment to evaluate the out-of-sample prescriptive ability of randomized robust pricing using a retail scanner data set. We show that randomized robust pricing outperforms the nominal pricing approach in the log-log and semi-log demand models, is comparable to deterministic robust pricing in the log-log model, and outperforms deterministic robust pricing in the semi-log model.
\end{enumerate}

The rest of this paper is organized as follows. Section~\ref{sec:literature_review} reviews the related literature on pricing, robust optimization and randomized robust optimization. Section~\ref{sec:problem_definition} formally defines the nominal price optimization problem, the deterministic robust price optimization problem and the randomized robust price optimization problem. Section~\ref{sec:benefits} analyzes the robust price optimization problem and provides conditions under which the price optimization problem is randomization-proof. Section~\ref{sec:finitePcal_convexUcal} presents our constraint generation approach for solving the RRPO problem when the price set is a finite set and the uncertainty set is a convex set, discusses how this approach can be adapted for different families of demand models and presents our result on the complexity of randomization. Section~\ref{sec:finitePcal_finiteUcal_body} provides a brief overview of our methodology for solving the RRPO problem when the price set and uncertainty sets are finite sets, with the details provided in Section~\ref{sec:finitePcal_finiteUcal} of the companion. Section~\ref{sec:results} presents our numerical experiments. Lastly, in Section~\ref{sec:conclusions}, we conclude and highlight some directions for future research.

\section{Literature review}
\label{sec:literature_review}

Our paper is closely related to three streams of research: pricing, robust optimization and the use of randomized strategies in optimization. We discuss each of these three streams below. \\

\noindent \textbf{Pricing optimization and demand models}. Optimal pricing has been extensively studied in many fields such as revenue management and marketing research; for a general overview of this research area, we refer readers to \cite{soon2011review} and \cite{gallego2019revenue}. An important stream of pricing literature is on static pricing, which involves setting a fixed price for a product. The most commonly considered demand models in the static pricing literature are the linear and log-log models. For example, the papers of \cite{zenor1994profit} and \cite{bernstein2003pricing} assume linear demand functions in the study of pricing strategies. The papers of \cite{reibstein1984optimal}, and \cite{montgomery1999analyst} use log-log (multiplicative) demand functions to represent aggregate demand. The paper of \cite{kalyanam1996pricing} considers a semi-log demand model. Besides linear, semi-log and log-log demand functions, another type of demand form that is extensively discussed in pricing literature is based on an underlying discrete choice model. For example, \cite{hanson1996optimizing}, \cite{aydin2000product}, and \cite{hopp2005product} consider the product line pricing problem under the multinomial logit (MNL) model. \cite{keller2014efficient} consider attraction demand models which subsume MNL models. \cite{li2011pricing} study the pricing problem with the nested logit (NL) models, and show the concavity of the profit function with respect to market share holds. \cite{gallego2014multiproduct} characterize the optimal pricing structure under the general nested logit model with product-differentiated price sensitivities and arbitrary nest coefficients. The papers of \cite{keller2013tractable} and \cite{zhang2018multiproduct} study the multiproduct pricing problem under the family of generalized extreme value (GEV) models which includes MNL and NL models as special cases. Our work differs from this prior work on multiproduct pricing in that the demand model is not assumed to be known, and that there is an uncertainty set of plausible demand models that could actually materialize. Correspondingly, the firm is concerned not with expected revenue under a single, nominal demand model, but with the worst-case revenue with respect to this uncertainty set of demand models. \\

%Another significant stream of pricing literature is to consider multiple-period pricing decisions where the prices of products change over time and there is a fixed inventory of each product; we refer readers to \cite{mcgill1999revenue}, \cite{elmaghraby2003dynamic}, \cite{bitran2003overview}, and \cite{talluri2004theory} for a comprehensive review on dynamic pricing strategies with inventory considerations. As in the static pricing literature, studies in dynamic pricing also vary in the types of demand models. \cite{besbes2015surprising} assume linear demand in a multiperiod single product pricing problem, and show that the corresponding pricing policy can perform well even under model misspecification. \cite{caro2012clearance} consider multiplicative models where the demand rate and price discount have the logarithmic relationship, and consider a multiproduct clearance pricing optimization problem for the fast-fashion retailer Zara. \cite{akccay2010joint} consider dynamic pricing under MNL models for horizontally differentiated products and show that the profit function is unimodal in prices, while \cite{song2021demand} and \cite{dong2009dynamic} reformulate the MNL profit as a concave function of its market share rather than prices. Our work focuses on the static, single-period setting where there is no inventory consideration, and is thus not directly related to this stream of the pricing literature. \\

\noindent \textbf{Robust Optimization}. Within the literature mentioned above, the demand models are assumed to be known exactly or are estimated from historical data. In practice, the decision maker often does not have complete knowledge of the demand model. Also, in many real applications, the lack of sales data makes it hard to obtain a good estimation of demand models, which leads to model misspecification and thus suboptimal pricing decisions. In operations research, this type of challenge is most commonly addressed using the framework of robust optimization, where uncertain parameters are assumed to belong to some uncertainty set and one optimizes for the worst-case objective of parameters within the set. We refer readers to \cite{ben2009robust}, \cite{bertsimas2011theory}, \cite{gabrel2014recent} and \cite{bertsimas2022robust} for a detailed overview of this approach. Robust optimization has been widely applied in various problem settings such as assortment optimization \citep{rusmevichientong2012robust,bertsimas2017robust,sturt2021value}, inventory management \citep{bertsimas2006robust,govindarajan2021distribution} and financial option pricing \citep{bandi2014robust,sturt2021nonparametric}.

Within this literature, our paper contributes to the substream that considers robust optimization for pricing. \cite{thiele2009multi} consider tractable robust counterparts to the deterministic multiproduct pricing problem with the budget-of-resource-consumption constraint in the case of additive demand uncertainty, and investigate the impact of uncertainty on the optimal prices of multiple products sharing capacitated resources. \cite{mai2019robust} consider robust multiproduct pricing optimization under the generalized extreme value (GEV) choice model and characterize the robust optimal solutions for unconstrained and constrained pricing problems. \cite{hamzeei2021robust} study the robust pricing problem with interval uncertainty of the price sensitivity parameters under the multi-product linear demand model. For robust dynamic pricing problems, \cite{lim2007relative} and \cite{lim2008robust} use relative entropy to represent uncertainty in the demand rate and thus the demand uncertainty can be expressed through a constraint on relative entropy. \cite{perakis2006competitive}, \cite{adida2010dynamic} and \cite{chen2018robust} all study robust dynamic pricing problems with demand uncertainty modeled by intervals. 
%\cite{perakis2006competitive} specifically study dynamic pricing of perishable products in the presence of competition where the uncertainty of demand is constrained within a closed and convex uncertainty set. \cite{adida2010dynamic} assume a linear demand function where uncertain parameters are within intervals constrained by a kind of budget-of-uncertainty constraints. \cite{chen2018robust} assume that the demand functional form is unknown and the uncertain demand are constrained within three types of interval bounds. 
\cite{harsha2019dynamic} study robust dynamic price optimization on an omnichannel network with cross-channel interactions in demand and supply where demand uncertainty is modeled through budget constraints.
From a different perspective, \cite{cohen2018dynamic} develop a data-driven framework for solving the robust dynamic pricing problem by directly using samples in the optimization. Specifically, the paper considers three types of robust objective (max-min, min-max regret and max-min ratio), and uses the given sampled scenarios to approximate the uncertainty set by a finite number of constraints. Another paper that considers data-driven pricing is the paper of \cite{yan2022representative}, which models customer demand using a representative agent model. In a representative agent model, choice probabilities arise as solutions to a constrained concave optimization problem over the unit simplex, where one finds a discrete probability distribution that maximizes the expected utility minus a convex perturbation function. By appropriately setting the perturbation function, one can recover various types of discrete choice models. The paper of \cite{yan2022representative} considers a method for identifying the perturbation function nonparametrically using price-demand data, and then subsequently solving a mixed-integer program to obtain the optimal prices using the estimated representative agent model. While our paper shares a high-level similarity with this work in that one seeks to address model misspecification, this prior work does not incorporate uncertainty around the estimated demand model, and does not consider randomness in the pricing decision. 

Our work differs from the majority of this body of work that takes a robust optimization approach to pricing in that the decision we seek to make is no longer a deterministic decision, but a randomized one. Within this body of work, the papers closest to our work are \cite{allouah2021optimal} and \cite{allouah2022pricing}. The paper of \cite{allouah2021optimal} considers a pricing problem under a valuation-based model of demand, where each customer has a valuation drawn from an unknown cumulative distribution function on the positive real line, and the firm only has one historical data point. The paper considers the pricing problem from a max-min ratio standpoint, where the firm seeks to find a pricing mechanism that maximizes the percentage attained of the true maximum revenue under the worst-case (minimum) valuation distribution consistent with the data point. The paper shows the approximation rate that is attainable given knowledge of different quantiles of the valuation distribution when the valuation distribution is a regular distribution or a monotone non-decreasing hazard rate distribution. The mechanisms that are proposed in the paper, which are mappings from the data point to a price, include deterministic ones that offer a fixed price, as well as randomized ones that offer different prices probabilistically. The paper of \cite{allouah2022pricing} considers a similar setting, where instead of knowing a point on the valuation CDF exactly or to within an interval, one has access to an IID sample of valuations drawn from the unknown valuation CDF, and similarly proposes deterministic and randomized mechanisms for this setting. 

With regard to \cite{allouah2021optimal} and \cite{allouah2022pricing}, our setup differs in a number of ways. First, our methodology focuses on a max-min revenue objective, as opposed to a max-min ratio objective that considers performance relative to oracle-optimal revenue. Our methodology also does not start from a valuation model, but instead starts from an aggregate demand model, and additionally considers the multi-product case in the general setting. Additionally, we do not take data as a starting point, but instead assume that the aggregate demand model is uncertain. Lastly, the overarching goals are different: while \cite{allouah2021optimal,allouah2022pricing} seek to understand the value of data, and how well one can do with limited data, our goal is to demonstrate that from the perspective of worst-case revenue, a randomized pricing strategy can be preferable over a deterministic fixed price, and to develop tractable computational methods for computing such strategies under commonly used demand models in the multi-product setting. \\

\noindent \textbf{Randomized strategies in optimization under uncertainty.} The conventional robust optimization problems mentioned above only consider deterministic solutions. In recent years, the benefit of using randomized strategies has received increasing attention in the literature on decision making under uncertainty and robust optimization. \cite{mastin2015randomized} study randomized strategies for min-max regret combinatorial optimization problems in the cases of interval uncertainty and uncertainty representable by discrete scenarios, and provide bounds on the gains from randomization for these two cases. \cite{bertsimas2016power} consider randomness in a network interdiction min-max problem where the interdictor can benefit from using a randomized strategy to select arcs to be removed. 

The paper of \cite{delage2019dice} considers the problem of making a decision whose payoff is uncertain and minimizing a risk measure of this payoff, and studies under what circumstances a randomized decision leads to lower risk than a deterministic decision. The paper characterizes the classes of randomization-receptive and randomization-proof risk measures in the absence of distributional ambiguity (i.e., classical stochastic programs), and discusses conditions under which problems with distributional ambiguity can benefit from randomized decisions. 

Subsequently, the paper of \cite{delage2022value} studies the value of randomized solutions for mixed-integer distributionally robust optimization problems. The paper develops bounds on the magnitude of improvement given by randomized solutions over deterministic solutions, and proposes a two-layer column generation method for solving single-stage and two-stage linear DRO problems with randomization. Our paper relates to \cite{delage2022value} in that we apply a similar two-layer column generation approach for solving the randomized robust pricing problem when the price set and the uncertainty set are both finite; we discuss this connection in more detail in our discussion of the paper of \cite{wang2020randomized} below. 
The paper of \cite{sadana2023value} develops a randomization approach for solving a distributionally robust maximum flow network interdiction problem with a conditional-value-at-risk objective, which is also solved using a column generation approach.

Our work is most closely related to the paper of \cite{wang2020randomized}. The paper of \cite{wang2020randomized} introduces randomization into the robust assortment optimization and characterizes the conditions under which a randomized strategy strictly improves worst-case expected revenues over a deterministic strategy. The paper proposes several different solution methods for finding an optimal distribution over assortments for the MNL, Markov chain and ranking-based models. % For the MNL model in particular, the paper adapts the two-layer column generation method of \cite{delage2022value} to solve the randomized robust assortment optimization problem when the uncertainty set is discrete.
Our paper shares a high-level viewpoint with the paper of \cite{wang2020randomized} in that revenue management decisions, such as assortment decisions and pricing decisions, are subject to uncertainty and from an operational point of view, have the potential to be randomized and to benefit from randomization. From a technical standpoint, several of our results on the benefit of randomization when the price set is discrete that are stated in Section~\ref{subsec:benefits_finite} are generalizations of results in \cite{wang2020randomized} to the pricing setting that we study. In terms of methodology, the solution approach we apply when the price set and uncertainty set are discrete in Section~\ref{sec:finitePcal_finiteUcal_body} (described more fully in Section~\ref{sec:finitePcal_finiteUcal}) is related to the approach in \cite{wang2020randomized} for the MNL model, as we also use the two-layer column generation scheme of \cite{delage2022value}. The main difference between the method in \cite{wang2020randomized} and our method lies in the nature of the subproblems. In the paper of \cite{wang2020randomized}, the primal subproblem is a binary sum of linear fractional functions problem, and the dual subproblem is essentially a mixture of multinomial logits assortment problem that can be reformulated as a mixed-integer linear program. In our paper, the primal and dual subproblems that are used to generate new price vectors and uncertainty realizations come from the underlying pricing problem and the structure of different demand models (linear, semi-log and log-log), which lead to different subproblems than in the assortment setting. In particular, in the semi-log and log-log cases, both the primal and dual subproblems can be formulated as mixed-integer exponential cone programs. In the semi-log and log-log cases specifically, the formulations of the dual subproblems, which are used to identify new price vectors to add, require a logarithmic transformation together with a biconjugate representation of the log-sum-exp function. This technique is also used to develop a constraint generation scheme for the randomized robust pricing problem when the price set is discrete and the uncertainty set is convex (Section~\ref{sec:finitePcal_convexUcal}); the subproblem in this case involves solving a nominal pricing problem under the semi-log or log-log model, which we are also able to reformulate exactly as a mixed-integer exponential cone program. As noted in the introduction, these are the first exact formulations of these problems using mixed-integer conic programming. \\

\noindent \textbf{Other work on randomization}. Lastly, we comment on several streams of work that use randomization but are unrelated to our paper. %We note that there is a substantial literature on solution methods that use randomness for deterministic optimization problems. For example, the optimization community has considered solving large-scale optimization problems by randomly sampling constraints/columns \citep{de2004constraint,calafiore2005uncertain,chen2020column}. Our work differs from this lterature in that randomization is used in an operational way to manage risk, and not as an algorithmic tool.
Within the revenue management community, there are instances where randomization is an operational aspect of the algorithm. For example, in network revenue management, the heuristic of probabilistic allocation control involves using the primal variable values in the deterministic linear program (DLP) to decide how frequently requests should be accepted or rejected \citep{jasin2012re}. Here, randomization is used to ensure that the long run frequency with which different requests are accepted or rejected is as close as possible to the DLP solution, which corresponds to an idealized upper bound on expected revenue. As another example, randomization is often also a part of methods for problems that involve learning. For example, \cite{ferreira2018online} propose a method for network revenue management where there is uncertainty in demand rates based on Thompson sampling, which is a method from the bandit literature that involves taking a random sample from the posterior distribution of an uncertain parameter and taking the action that is optimal with respect to that sample. Here, randomization is a way of ensuring that the decision maker explores possibly suboptimal actions. In our work, the focus is not on using randomization to achieve better expected performance or using randomization to achieve a balance between exploration and exploitation, but rather to operationalize randomization to achieve better worst-case performance.

\section{Problem definition}
\label{sec:problem_definition}
In this section, we begin by defining the nominal price optimization problem in Section~\ref{subsec:problem_definition_NPO}. We subsequently define the deterministic robust price optimization problem in Section~\ref{subsec:problem_definition_DRPO}. Lastly, we define the randomized robust price optimization problem in Section~\ref{subsec:problem_definition_RRPO}.

\subsection{Nominal price optimization problem}
\label{subsec:problem_definition_NPO}

We assume that the firm offers $I$ products, indexed from 1 to $I$. We let $p_i$ denote the price of product $i \in [I]$, where we use the notation $[n] = \{1,\dots,n\}$ for any positive integer $n$. We use $\pb = (p_1,\dots, p_I)$ to denote the vector of prices. We assume that the price vector $\pb$ is constrained to lie in the set $\Pcal \subseteq \Rbb_+^I$, where $\Rbb_+$ is the set of nonnegative real numbers. 

We let $d_i$ denote the demand function of product $i$, so that $d_i(\pb)$ denotes the demand of product $i$ when the price vector $\pb$ is chosen. The revenue function $R(\cdot)$ can then be written as $R(\pb) = \sum_{i =1}^I p_i \cdot d_i(\pb)$.

The nominal price optimization (NPO) problem can be written simply as
\begin{equation*}
\textsc{NPO}: \quad \max_{\pb \in \Pcal} R(\pb).
\end{equation*}

%{\comebacktothis -- give examples of different demand models. (Linear, semi-log, log-log, MNL.) }

There are numerous demand models that can be used in practice, which lead to different price optimization problems; we briefly review some of the more popular ones here.

\begin{enumerate}
	\item \emph{Linear demand model}: A linear demand model is defined by parameters $\alphab \in \Rbb^I$, $\betab \in \Rbb^{I}$, $\gammab = (\gamma_{i,j})_{i, j \in [I], i \neq j} \in \Rbb^{I \cdot(I-1)}$. The demand function $d_i(\cdot)$ of each product $i \in [I]$ has the form
	\begin{equation}
	d_i(\pb) = \alpha_i - \beta_i p_i + \sum_{j \neq i} \gamma_{i,j} p_j,
	\end{equation}
	where $\beta_i \geq 0$ is the own-price elasticity parameter of product $i$, which indicates how much demand for product $i$ is affected by the price of product $i$, whereas $\gamma_{i,j}$ is a cross-price elasticity parameter that describes how much demand for product $i$ is affected by the price of a different product $j$. Note that $\gamma_{i,j}$ can be positive, which generally corresponds to products $i$ and $j$ being substitutes (i.e., when the price of product $j$ increases, customers tend to switch to product $i$), or negative, which corresponds to products $i$ and $j$ being complements (i.e., products $i$ and $j$ tend to be purchased together, so when the price of product $j$ increases, this causes a decrease in demand for product $i$).
	The corresponding revenue function $R(\cdot)$ is then
	\begin{equation}
	R(\pb) = \sum_{i=1}^I p_i \cdot (\alpha_i - \beta_i p_i + \sum_{j \neq i} \gamma_{i,j} p_j).
	\end{equation}
	\item \emph{Semi-log demand model}: A semi-log demand model is defined by parameters $\alphab \in \Rbb^I$, $\betab \in \Rbb^{I}$, $\gammab = (\gamma_{i,j})_{i, j \in [I], i \neq j} \in \Rbb^{I \cdot(I-1)}$. In a semi-log demand model, the logarithm of the demand function $d_i(\cdot)$ of each product $i \in [I]$ has a linear form in the prices $p_1,\dots, p_I$: 
	\begin{equation}
	\log d_i(\pb) = \alpha_i - \beta_i p_i + \sum_{j \neq i} \gamma_{i,j} p_j.
	\end{equation}
	This implies that the demand function is 
	\begin{equation}
	d_i(\pb) = e^{\alpha_i - \beta_i p_i + \sum_{j \neq i} \gamma_{i,j} p_j}.
	\end{equation}
	The corresponding revenue function $R(\cdot)$ is then 
	\begin{equation}
	R(\pb) = \sum_{i=1}^I p_i \cdot e^{ \alpha_i - \beta_i p_i + \sum_{j \neq i} \gamma_{i,j} p_j}.
	\end{equation}
	\item \emph{Log-log demand model}: A log-log demand model is defined by parameters $\alphab \in \Rbb^I$, $\betab \in \Rbb^{I}$, $\gammab = (\gamma_{i,j})_{i, j \in [I], i \neq j} \in \Rbb^{I \cdot(I-1)}$. In a log-log demand model, the logarithm of the demand function $d_i(\cdot)$ of each product $i \in [I]$ has a linear form in the log-transformed prices $\log p_1, \dots, \log p_I$: 
	\begin{equation}
	\log d_i(\pb) = \alpha_i - \beta_i \log p_i + \sum_{j \neq i} \gamma_{i,j} \log p_j.
	\end{equation}
	This implies that the demand function for product $i$ is 
	\begin{align}
	d_i(\pb) & = e^{ \alpha_i - \beta_i \log p_i + \sum_{j \neq i} \gamma_{i,j} \log p_j } \\
			&  = e^{\alpha_i} p_i^{-\beta_i} \cdot \prod_{j \neq i} p_j^{\gamma_{i,j}},
	\end{align}
	and that the revenue function is therefore
	\begin{align}
	R(\pb) & = \sum_{i=1}^I p_i \cdot e^{\alpha_i} p_i^{-\beta_i} \cdot \prod_{j \neq i} p_j^{\gamma_{i,j}} \\
	& = \sum_{i=1}^I e^{\alpha_i} \cdot p_i^{1 - \beta_i} \cdot \prod_{j \neq i} p_j^{\gamma_{i,j}}.
	\end{align}
%	\item \emph{Multinomial logit demand model}: A multinomial logit (MNL) demand model is defined by parameters  $\alphab \in \Rbb^I$, $\betab \in \Rbb^{I}$. Without loss of generality, we assume that the market size is normalized to 1. The demand function $d_i(\cdot)$ of each product $i \in [I]$ is given by 
%	\begin{equation}
%	d_i(\pb) = \frac{ e^{\alpha_i - \beta_i p_i}}{1 + \sum_{i' = 1}^I e^{\alpha_{i'} - \beta_{i'} p_{i'}} }.
%	\end{equation}
%	The revenue function is therefore 
%	\begin{equation}
%	R(\pb) = \frac{ \sum_{i=1}^I p_i e^{\alpha_i - \beta_i p_i}}{ 1 + \sum_{i = 1}^I e^{\alpha_{i} - \beta_{i} p_{i}} }.
%	\end{equation}
\end{enumerate}
%We make one important remark about the setup of our price optimization problem. In our price optimization problem, we are solely focused on maximizing revenue; our formulation is not focused on maximizing profit and does not incorporate marginal costs. The reason for this is 

\subsection{Deterministic robust price optimization problem}
\label{subsec:problem_definition_DRPO}

We now define the deterministic robust price optimization (DRPO) problem. To define this problem abstractly, we let $\Rcal$ denote an uncertainty set of possible revenue functions. The DRPO problem is to then maximize the worst-case revenue, where the worst-case is the minimum revenue of a given price vector taken over all revenue functions in $\Rcal$. Mathematically, this problem can be written as
\begin{equation*}
\textsc{DRPO}: \quad \max_{\pb \in \Pcal} \min_{R \in \Rcal} R(\pb).
\end{equation*}
We use $Z^*_{\DR}$ to denote the optimal objective value of the DRPO problem. 

%The revenue function uncertainty set $\Rcal$ can arise in a couple of different ways, which we outline below. 

%arises out of a demand function uncertainty set $\Dcal$. In particular, $\Rcal = \{ R \equiv \sum_{i=1}^I p_i \cdot d_i(\pb) \mid \db(\cdot) \in \Dcal \}$. The demand function uncertainty set $\Dcal$ itself can be defined in several ways, depending on the type of uncertainty we wish to protect against. We provide several examples below. 

Although $\Rcal$ can be defined in many different ways, we will now focus on one general case that we will assume for most of our subsequent results in Sections~\ref{sec:benefits}, \ref{sec:finitePcal_convexUcal} and \ref{sec:finitePcal_finiteUcal_body}. Suppose that we fix the demand model to a specific parametric family, such as a log-log demand model. Let $\ub$ denote the vector of demand model parameters. For example, for log-log, $\ub$ would be the tuple $\ub = (\alphab, \betab, \gammab)$. Let $\Ucal \subseteq \Rbb^M$ denote a set of possible values of $\ub$, where the integer $M$ is the dimension of $\ub$; $\Ucal$ is then an uncertainty set of model parameters. With a slight abuse of notation, let $d_i(\pb, \ub)$ denote the demand for product $i$ when the demand model parameters are specified by $\ub$. Then $\Rcal$ can be defined as:
\begin{equation}
\Rcal = \{ R(\pb) \equiv \sum_{i = 1}^I p_i \cdot d_i(\pb, \ub) \mid \ub \in \Ucal\},
\end{equation}
i.e., it is all of the possible revenue functions spanned by the uncertain parameter vector $\ub$ in $\Ucal$. 

To express the DRPO problem in this setting more conveniently, we will abuse our notation slightly and use $R(\pb, \ub)$ to denote the revenue function evaluated at a price vector $\pb$ with a particular parameter vector $\ub$ specified. With this abuse of notation, the DRPO problem can be written as
\begin{equation*}
\textsc{DRPO}: \quad \max_{\pb \in \Pcal} \min_{\ub \in \Ucal} R(\pb, \ub).
\end{equation*}

%{\comebacktothis -- let's think more carefully about notation for uncertainty sets. (Use $\Ucal$ or $\Dcal$?) Revenue functions are uncertain; demand functions are uncertain; demand function parameters are uncertain. We should be careful not to mix these up.}

\subsection{Randomized robust price optimization problem}
\label{subsec:problem_definition_RRPO}

In the DRPO problem, we assume that the decision maker will deterministically implement a single price vector $\pb$ in the face of uncertainty in the revenue function. In the RRPO problem, we instead assume that the decision maker will randomly select a price vector $\pb$ according to some distribution $F$ over the feasible price set $\Pcal$. Under this assumption, we can write the RRPO problem as
\begin{equation*}
\textsc{RRPO}: \quad \max_{F \in \Fcal} \min_{R \in \Rcal} \int_{\Pcal} R(\pb) \, dF(\pb),
\end{equation*}
where $\Fcal$ is the set of all distributions supported on $\Pcal$. We use $Z^*_{\RR}$ to denote the optimal objective value of the RRPO problem. Note that $Z^*_{\RR} \geq Z^*_{\DR}$. This is because for every $\pb' \in \Pcal$, the distribution $F(\cdot) = \delta_{\pb'}(\cdot)$, where $\delta_{\pb'}(\cdot)$ is the Dirac delta function at $\pb'$, is contained in $\Fcal$. For this distribution, $\min_{R \in \Rcal} \int_{\Pcal} R(\pb) \, dF(\pb) = \min_{R \in \Rcal} R(\pb')$, which is exactly the worst-case revenue of deterministically selecting $\pb'$.

A special instance of this problem arises when $\Pcal$ is a discrete, finite set. In this case, $F$ is a discrete probability distribution, and one can re-write the outer problem as an optimization problem over a discrete probability distribution $\pib = (\pi_{\pb})_{\pb \in \Pcal}$ over $\Pcal$:
\begin{equation*}
\textsc{RRPO-D}: \quad \max_{\pib \in \Delta_{\Pcal}} \min_{R \in \Rcal} \sum_{\pb \in \Pcal} R(\pb) \pi_{\pb},
\end{equation*}
where we use $\Delta_S$ to denote the $(|S|-1)$-dimensional unit simplex, i.e., $\Delta_S = \{ \pib \in \Rbb^S \mid \sum_{i \in S} \pi_i = 1, \ \pi_i \geq 0 \ \forall \ i \in S \}$.

Lastly, under the assumption that $\Rcal$ is the set of revenue functions of a fixed demand model family whose parameter vector $\ub$ belongs to a parameter uncertainty set $\Ucal$, we can restate the RRPO problem when $\Pcal$ is a generic set and when $\Pcal$ is finite as 
\begin{align}
\textsc{RRPO}: \quad & \max_{F \in \Fcal} \min_{\ub \in \Ucal} \int_{\Pcal} R(\pb, \ub) \, dF(\pb), \\
\textsc{RRPO-D}: \quad & \max_{\pib \in \Delta_{\Pcal}} \min_{\ub \in \Ucal} \sum_{\pb \in \Pcal} R(\pb,\ub) \pi_{\pb}.
\end{align}

\section{Benefits of randomization}
\label{sec:benefits}

In this section, we analyze when randomization can be beneficial. To aid us, we introduce some additional nomenclature in this section, which follows the terminology established in the prior literature on randomized robust optimization \citep{delage2019dice,delage2022value,wang2020randomized}. We say that a robust price optimization (RPO) problem is \emph{randomization-receptive} if $Z^*_{\RR} > Z^*_{\DR}$, that is, randomizing over price vectors leads to a higher worst-case revenue than deterministically selecting a single price vector. Otherwise, we say that a RPO problem is \emph{randomization-proof} if $Z^*_{\RR} = Z^*_{\DR}$, that is, there is no benefit from randomizing over price vectors.

In the following three sections, we derive three classes of results that establish when the RPO problem is randomization-proof. The first condition (Section~\ref{subsec:benefits_concave}) for randomization-proofness applies in the case when $\Pcal$ is a convex set and $\Rcal$ is an arbitrary set of revenue functions each of which is concave in $\pb$. The second and third conditions apply to the case where $\Rcal$ arises out of a single demand model, where the parameter vector $\ub$ belongs to an uncertainty set. The second condition (Section~\ref{subsec:benefits_quasiconvexity_quasiconcavity}) applies in the case when $\Pcal$ and $\Ucal$ are compact, convex sets and $R(\pb, \ub)$ obeys certain quasiconvexity and quasiconcavity properties. The third condition (Section~\ref{subsec:benefits_finite}) is for the case when $\Pcal$ is finite and involves a certain minimax condition being met; as corollaries, we show that randomization-proofness occurs if the DRPO problem satisfies a strong duality property, and that randomization-proofness occurs if the DRPO solution solves the nominal price optimization problem solution at the worst-case $\ub^*$. Along the way, we also give examples where our results establish that certain families of RPO problems are randomization-proof, and also highlight how the results fail to hold when certain hypotheses are relaxed. 

The main takeaway from these results is that the set of RPO problems that are randomization-proof is small. As we will see, the conditions under which a RPO problem will be randomization-proof are delicate and quite restrictive, and are satisfied only for certain very special cases; in most other realistic cases, the RPO problem will be randomization-receptive. Consequently, in Sections~\ref{sec:finitePcal_convexUcal} and \ref{sec:finitePcal_finiteUcal_body}, we will develop algorithms for solving the RRPO problem when the candidate price vector $\Pcal$ is finite, and in Sections~\ref{sec:results} we will show a wide range of both synthetic and real data instances in which the RPO problem is randomization-receptive.

Before we move on to our results, we comment on the connection between randomization-proofness and results from game theory. In particular, the robust pricing problem that we study can be viewed as a two-player zero-sum game: the first player is the firm, who chooses an action $\pb$ from the strategy space $\Pcal$, and the second player is nature, who chooses an action $\ub$ from the strategy space $\Ucal$. The payoff function that the first player maximizes is $R(\pb, \ub)$, while the payoff function that the second player maximizes is $-R(\pb, \ub)$, which is equivalent to minimizing $R(\pb, \ub)$. In game theory, it is common to analyze when a game has a \emph{pure strategy Nash equilibrium}. A pure strategy Nash equilibrium in our setting is a pair of actions $(\pb^*, \ub^*)$ such that each player cannot improve their payoff by unilaterally deviating from their action; mathematically, it satisfies
\begin{align}
R(\pb, \ub^*) \leq R(\pb^*, \ub^*), \quad \forall \ \pb \in \Pcal, \\
R(\pb^*, \ub^*) \leq R(\pb^*, \ub), \quad \forall \ \ub \in \Ucal.
\end{align}
Note that a pure strategy Nash equilibrium $(\pb^*, \ub^*)$ is exactly a saddle point of the function $R$. By standard optimization arguments, the existence of a saddle point is equivalent to strong duality holding, that is,
\begin{equation}
\max_{\pb \in \Pcal} \min_{\ub \in \Ucal} R(\pb, \ub) = \min_{\ub \in \Ucal} \max_{\pb \in \Pcal} R(\pb, \ub).
\end{equation}
When strong duality holds, the RPO problem is randomization-proof, because
%\begin{align}
%Z^*_{\RR} & = \max_{F \in \Fcal} \min_{\ub \in \Ucal} \int R(\pb, \ub) \, dF(\pb) \\
%& \leq \min_{\ub \in \Ucal} \max_{F \in \Fcal}  \int R(\pb, \ub) \, dF(\pb) \\ 
%& = \min_{\ub \in \Ucal} \max_{\pb \in \Pcal}  R(\pb, \ub) \\ 
%& = \max_{\pb \in \Pcal} \min_{\ub \in \Ucal} R(\pb, \ub) \\
%& = Z^*_{\DR},
%\end{align}
\begin{align}
Z^*_{\RR} = \max_{F \in \Fcal} \min_{\ub \in \Ucal} \int R(\pb, \ub) \, dF(\pb) \leq \min_{\ub \in \Ucal} \max_{F \in \Fcal}  \int R(\pb, \ub) \, dF(\pb)  = \min_{\ub \in \Ucal} \max_{\pb \in \Pcal}  R(\pb, \ub)  = \max_{\pb \in \Pcal} \min_{\ub \in \Ucal} R(\pb, \ub)  = Z^*_{\DR},
\end{align}
i.e., $Z^*_{\RR} \leq Z^*_{\DR}$, and together with $Z^*_{\RR} \geq Z^*_{\DR}$, this immediately implies that $Z^*_{\RR} = Z^*_{\DR}$. Thus, when one views the RPO problem as a two-player zero-sum game, if a pure strategy Nash equilibrium exists for this game, then the RPO problem is randomization-proof. We leverage this connection to establish one of our results, Theorem~\ref{theorem:quasiconcave_quasiconvex}, where we invoke the well-known Debreu-Glicksberg-Fan theorem (see Theorem~1.2 of \citealt{fudenberg1991game}) to assert the existence of a pure strategy Nash equilibrium. However, the two conditions are not equivalent: it can be the case that a pure strategy Nash equilibrium does not exist, and yet the RPO problem is randomization-proof. In Section~\ref{subsec:benefits_concave}, we note an instance where the RPO problem is randomization-proof but a pure strategy Nash equilibrium does not exist, and in this instance the price set $\Pcal$ is a continuous interval. In Section~\ref{subsec:benefits_finite}, we note a similar instance where the price set $\Pcal$ is finite.  %Our first result (Theorem~\ref{theorem:concave_revenue_function_set}) gives one setting where the RPO problem can be randomization-proof, but a pure strategy Nash equilibrium does not exist.

\subsection{Concave revenue function uncertainty sets}
\label{subsec:benefits_concave}

Our first major result is for the case where $\Rcal$ consists of concave revenue functions. 
\begin{theorem}
Suppose that $\Pcal$ is a convex set and that $\Rcal$ is such that every $R \in \Rcal$ is a concave function of $\pb$. Then the RPO problem is randomization-proof, that is, $Z^*_{\RR} = Z^*_{\DR}$. \label{theorem:concave_revenue_function_set}
\end{theorem}
The proof of this result (see Section~\ref{proof_theorem:concave_revenue_function_set} of the ecompanion) follows from a simple application of Jensen's inequality. We pause to make a few important comments about this result. First, one aspect of this result that is special is that $\Rcal$ can be a very general set: it could be countable or uncountable, and it could consist of revenue functions corresponding to different families of demand models. This will not be the case for our later results in Section~\ref{subsec:benefits_quasiconvexity_quasiconcavity} and \ref{subsec:benefits_finite}, which require that $\Rcal$ is defined based on a single demand model family.

Second, we remark that the conditions that all functions in $\Rcal$ be concave and that $\Pcal$ be a convex set cannot be relaxed in general. In Section~\ref{subsec:examples_benefits_concave} of the ecompanion, we provide examples to illustrate how Theorem~\ref{theorem:concave_revenue_function_set} can fail to hold if $\Rcal$ contains a non-concave function, or if $\Pcal$ is a non-convex set.

Third, we remark on the connection between Theorem~\ref{theorem:concave_revenue_function_set} and results in game theory. As alluded to at the end of Section~\ref{sec:benefits}, Theorem~\ref{theorem:concave_revenue_function_set} is \emph{not} a consequence of any known result on the existence of pure strategy Nash equilibria in game theory. The reason for this is that in general, strong duality does not hold in this setup, i.e., it is \emph{not} the case that
\begin{equation}
\max_{\pb \in \Pcal} \min_{R \in \Rcal} R(\pb) = \min_{R \in \Rcal} \max_{\pb \in \Pcal} R(\pb).
\end{equation}
In Section~\ref{subsec:examples_benefits_concave} of the ecompanion we provide an example (Example~\ref{example:concave_strong_duality_fails}) of an RPO problem that is randomization-proof, but for which strong duality does not hold (and hence, the corresponding game does not have a pure strategy Nash equilibrium). 

To provide a different perspective, as we previously noted, the RPO problem can be viewed as a game between the firm and nature. In Theorem~\ref{theorem:concave_revenue_function_set}, the set of revenue functions is not indexed by any uncertain parameter $\ub$, but one can artificially identify each function $R \in \Rcal$ with some $\ub$, i.e., with some abuse of notation $R(\pb) = R(\pb, \ub)$ for some $\ub$. If one takes this viewpoint, then the set $\Ucal$ of indexing $\ub$'s can be arbitrary, and need not be convex or compact. Additionally, for a fixed $\pb$, the function $R(\pb, \ub)$ could be arbitrary in $\ub$ and need not be convex/quasi-convex or even continuous. We are not aware of any general result asserting that a pure strategy Nash equilibrium exists in this setting, and in light of the example presented in Section~\ref{subsec:examples_benefits_concave}, such a result cannot exist. 

Lastly, Theorem~\ref{theorem:concave_revenue_function_set} has a number of implications for different classes of demand models. In Section~\ref{subsec:examples_benefits_concave_applications}, we document several applications of Theorem~\ref{theorem:concave_revenue_function_set}. In particular, we present conditions under which single-product pricing under linear demand, multi-product pricing under linear demand, single-product pricing under semi-log demand and single-product pricing under log-log demand are randomization-proof.

\subsection{Quasiconcavity in $\pb$ and quasiconvexity in $\ub$}
\label{subsec:benefits_quasiconvexity_quasiconcavity}

The second result we establish concerns the RRPO problem when there is a demand parameter uncertainty set $\Ucal$. In this case, the RRPO and DRPO problems are 
\begin{align*}
\textsc{RRPO}: & \ \max_{F \in \Fcal} \min_{\ub \in \Ucal} \int_{\Pcal} R(\pb, \ub) \ dF(\pb), \\
\textsc{DRPO}: & \ \max_{\pb \in \Pcal} \min_{\ub \in \Ucal}  R(\pb, \ub).
\end{align*}
%We make the following assumption about $R$.
%\begin{assumption}
%$R$ is a continuous function of $(\pb,\ub)$. \label{assumption:R_pb_ub_continuous}
%%For any fixed $\pb \in \Pcal$, $R(\pb, \ub)$ is continuous in $\ub$.  \label{assumption:R_ub_continuous}
%%\end{assumption}
%%\begin{assumption}
%%For any fixed $\ub \in \Ucal$, $R(\pb, \ub)$ is continuous in $\pb$. \label{assumption:R_pb_continuous}
%\end{assumption}

We then have the following result. 
\begin{theorem} \label{theorem:quasiconcave_quasiconvex} 
Suppose that $\Pcal \subseteq \Rbb^I$ and $\Ucal \subseteq \Rbb^M$ are compact convex sets. %Suppose that $\int_{\Pcal} R(\pb, \ub) \, dF(\pb)$ is a quasiconvex function of $\ub$ on $\Ucal$ for any $F \in \Fcal$. 
Suppose that $R(\pb, \ub)$ is quasiconcave in $\pb$ on $\Pcal$ for any $\ub \in \Ucal$, quasiconvex in $\ub$ on $\Ucal$ for any $\pb \in \Pcal$ and continuous in $(\pb, \ub)$. Then, the robust price optimization problem is randomization-proof, that is, $Z^*_{\DR} = Z^*_{\RR}$. 
\end{theorem}
The proof of Theorem~\ref{theorem:quasiconcave_quasiconvex} (see Section~\ref{proof_theorem:quasiconcave_quasiconvex} of the ecompanion) follows from viewing the robust price optimization problem as a two-player zero-sum game, and applying the Debreu-Glicksberg-Fan theorem (see \citealt{fudenberg1991game}) to assert the existence of a pure-strategy Nash equilibrium. As we noted in  Section~\ref{sec:benefits}, such a Nash equilibrium is a saddle point of the function $R$ and automatically implies that strong duality holds, which establishes that the problem is randomization-proof.

This result allows us to show that a larger number of RPO problems are randomization-proof. In Section~\ref{subsec:examples_benefits_quasiconcave_quasiconvex_applications} of the ecompanion, we show that single-product pricing under semi-log and log-log demand is randomization-proof even in settings where the revenue function is not concave in $\pb$ for some values of $\ub$ and one cannot invoke Theorem~\ref{theorem:concave_revenue_function_set}. With regard to these two examples, we note that they critically rely on the revenue function being log-concave and therefore quasiconcave, which is only the case for single product price optimization problems. Log-concavity and quasiconcavity are in general not preserved under addition (i.e., the sum of quasiconcave functions is not always quasiconcave, and the sum of log-concave functions is not always log-concave), and so Theorem~\ref{theorem:quasiconcave_quasiconvex} will in general not be applicable for multiproduct pricing problems involving the semi-log or log-log demand model.

\subsection{Finite price set $\Pcal$}
\label{subsec:benefits_finite}

In this section, we analyze randomization-proofness when $\Pcal$ is a finite set. To study this setting, let us define the set $\Qcal$ as the set of all probability distributions supported on $\Ucal$. We note that these results are adaptations of several results from \cite{wang2020randomized} to the pricing setting that we study, which develop analogous conditions for randomization-proofness for the robust assortment optimization problem. 

Our first result establishes that randomization-proofness is equivalent to the existence of a distribution $Q$ over $\Ucal$ under which any price vector's expected performance is no better than the deterministic robust optimal value. 

\begin{theorem} \label{theorem:finite_Pcal}
%Suppose that $R(\pb, \ub)$ is continuous in $\ub$ for any fixed $\pb \in \Pcal$. 
A robust price optimization problem with finite $\Pcal$ is randomization-proof if and only if there exists a distribution $Q \in \Qcal$ such that for all $\pb \in \Pcal$, 
\begin{align}
\int_{\Ucal} R(\pb, \ub) \, dQ(\ub) \leq Z^*_{\DR}. \label{eq:finite_Pcal_intQ_leq_Z_DR}
\end{align}
\end{theorem}

To prove this result, we use Sion's minimax theorem to establish that 
\begin{equation}
Z^*_{\RR} = \inf_{Q \in \Qcal} \max_{\pb \in \Pcal} \int_{\Ucal} R(\pb, \ub) \, dQ(\ub);
\end{equation}
with that result in hand, the condition in Theorem~\ref{theorem:finite_Pcal} is equivalent to establishing that 
\begin{align*}
Z^*_{\DR} \geq \inf_{Q \in \Qcal} \max_{\pb \in \Pcal} \int_{\Ucal} R(\pb, \ub) \, dQ(\ub)  = Z^*_{\RR},
\end{align*}
which, together with the inequality $Z^*_{\DR} \leq Z^*_{\RR}$ immediately yields randomization-proofness. We note that this result is analogous to Theorem 1 in \cite{wang2020randomized}, which provides a similar necessary and sufficient condition for randomization-proofness in the context of robust assortment optimization. Our proof, which relies on Sion's minimax theorem, is perhaps slightly more direct than the proof of Theorem 1 in \cite{wang2020randomized}, although this is a matter of taste.

Our next two results are consequences of this theorem. The first states that a RPO problem will be randomization-proof if it obeys strong duality. The second states that, a robust price optimization problem is randomization-proof if the deterministic robust price vector $\pb^*_{\DR}$ is an optimal solution of the nominal price optimization problem under a worst-case $\ub^*$ that attains the worst-case objective under $\pb^*_{\DR}$, and also presents a partial converse if one additionally assumes that $\ub^*$ is unique. We note that these results are both analogous to Corollaries 1 and 2 in \cite{wang2020randomized}.

\begin{corollary} \label{corollary:finite_Pcal_strong_duality}
Suppose that $\Pcal$ is finite. 
A robust price optimization problem with finite $\Pcal$ is randomization-proof if it satisfies strong duality:
\begin{align}
\max_{\pb \in \Pcal} \min_{\ub \in \Ucal} R(\pb, \ub) = \min_{\ub \in \Ucal} \max_{\pb \in \Pcal}  R(\pb, \ub). \label{eq:strong_duality}
\end{align}
\end{corollary}

\begin{corollary}
Suppose that $\Pcal$ is finite. Suppose that $\pb^*_{\DR} \in \arg \max_{\pb \in \Pcal} \min_{\ub \in \Ucal} R(\pb, \ub)$ is an optimal solution of the deterministic robust price optimization problem, and suppose that $\ub^* \in \arg \min_{\ub \in \Ucal} R(\pb^*_{\DR}, \ub)$ is an optimal solution of the worst-case problem for $\pb^*_{\DR}$. 
\begin{enumerate}
\item[a)] If $\pb^*_{\DR} \in \arg \max R(\pb, \ub^*)$, then the robust price optimization problem is randomization-proof. 
\item[b)] If the robust price optimization problem is randomization-proof and $\ub^*$ is the unique solution of $\min_{\ub \in \Ucal} R(\pb^*_{\DR}, \ub)$, then $\pb^*_{\DR} \in \arg \max R(\pb, \ub^*)$. 
\end{enumerate} 
\label{corollary:finite_Pcal_pDR}
\end{corollary}

Corollary~\ref{corollary:finite_Pcal_strong_duality} follows because the left hand side in \eqref{eq:strong_duality} is exactly $Z^*_{\DR}$, while the right hand side is an upper bound on $Z^*_{\RR}$. Part a) of Corollary~\ref{corollary:finite_Pcal_pDR} follows because $(\pb^*_{\DR}, \ub^*)$ is a saddle point, which implies that strong duality holds. Part b) of Corollary~\ref{corollary:finite_Pcal_pDR} follows by using Theorem~\ref{theorem:finite_Pcal} to assert the existence of a distribution $Q$ over $\Ucal$ satisfying~\eqref{eq:finite_Pcal_intQ_leq_Z_DR}, using the uniqueness of $\ub^*$ to argue that this distribution must place unit mass on $\ub^*$, and simplifying \eqref{eq:finite_Pcal_intQ_leq_Z_DR} using this fact to show that $\pb^*_{\DR}$ is the required optimal solution.

With regard to part (b) of Corollary~\ref{corollary:finite_Pcal_pDR}, we note that the uniqueness requirement for $\ub^*$ cannot in general be relaxed. In Section~\ref{counterexample_corollary:finite_Pcal_pDR} of the ecompanion, we show an instance where $\min_{\ub \in \Ucal} R(\pb^*_{\DR}, \ub)$ has multiple optimal solutions and the problem is randomization-proof, but $\pb^*_{\DR}$ does not solve $\max_{\pb \in \Pcal} R(\pb, \ub')$ for every $\ub'$ that solves $\min_{\ub \in \Ucal} R(\pb^*_{\DR}, \ub)$. Additionally, this instance also illustrates that even when $\Pcal$ is finite, it is possible that the robust price optimization problem may be randomization-proof without strong duality holding.

An interesting follow-on question to these results is under what conditions are strong duality and randomization-proofness equivalent when $\Pcal$ is finite. It turns out that this is the case when $R(\pb,\ub)$ is convex in $\ub$ and $\Ucal$ is convex. We state and prove this result in Section~\ref{proof_proposition:equivalence_SD_randproof_finitePcal} of the ecompanion. Additionally, we demonstrate through an example in Section~\ref{subsec:examples_benefits_DR_RR_dual_all_different} that the deterministic robust objective $Z^*_{\DR}$, the randomized robust objective $Z^*_{\RR}$ and the dual objective $\min_{\ub \in \Ucal} \max_{\pb \in \Pcal} R(\pb, \ub)$ can in general be distinct for the finite $\Pcal$ setting.

We conclude this section with the remark that the sufficient conditions for randomization-proofness in Corollaries~\ref{corollary:finite_Pcal_strong_duality} and \ref{corollary:finite_Pcal_pDR} are rather stringent and demanding. With regard to Corollary~\ref{corollary:finite_Pcal_strong_duality}, strong duality is in general unlikely to hold given that $\Pcal$ is a finite set. With regard to Corollary~\ref{corollary:finite_Pcal_pDR}, we note that in general, the solution of the deterministic robust price optimization problem is unlikely to also be an optimal solution of an appropriately defined non-trivial nominal price optimization problem; this is frequently not the case in many applications of robust optimization outside of pricing.\footnote{ For example, the nominal problem $\max\{ \cb^T \xb  \mid \oneb^T \xb = 1, \ \xb \geq \zerob\}$ usually has solutions of the form $\xb^* = \eb_i$, where $\eb_i$ is the standard unit vector in $\Rbb^n$, but the robust problem $\max\{ \min_{\cb \in \Ucal} \cb^T \xb  \mid \oneb^T \xb = 1, \ \xb \geq \zerob\}$ could in general have solutions that are in the interior of the unit simplex, which would never be optimal for the nominal problem except for the trivial case where all $c_i$'s are the same. As another example, Section~4.5 of \cite{bertsimas2017robust} describes an experiment involving a robust assortment optimization problem with a finite uncertainty set of three choice models, where the robust problem and the nominal problems under each of the choice models are solved heuristically. The robust solution is not an optimal solution for the nominal problem under any of the three choice models (in fact, the robust solution in this case turns out to be a ``blend'' of the three nominal product lines).}  Given this, these conditions are suggestive of the fact that most robust price optimization problems will be randomization-receptive. This motivates our study of solution algorithms for numerically solving the RRPO problem in the next two sections.

%* Solution algorithm for discrete price set \Pcal and convex uncertainty set \Ucal 

\section{Solution algorithm for finite price set $\Pcal$, convex uncertainty set $\Ucal$}
\label{sec:finitePcal_convexUcal}

In this section, we describe a general solution algorithm for solving the RRPO problem when the price set $\Pcal$ is a finite set, and the uncertainty set $\Ucal$ is a general convex uncertainty set. Section~\ref{subsec:finitePcal_convexUcal_algorithm} describes the general solution algorithm, which is a constraint generation algorithm that involves solving a nominal pricing problem over $\Pcal$ as a subroutine. Sections~\ref{subsec:finitePcal_convexUcal_linear}, \ref{subsec:finitePcal_convexUcal_semilog} and \ref{subsec:finitePcal_convexUcal_loglog} describe how the solution algorithm specializes to the cases of the linear, semi-log and log-log demand models, respectively, and in particular, how the nominal pricing problem can be solved for each of these three cases; the formulations we present for the semi-log and log-log models here may be of independent interest as they are, to the best of our knowledge, the first exact mixed-integer convex formulations for the multi-product pricing problem under a finite price set for these demand models. %Lastly, Section~\ref{subsec:finitePcal_convexUcal_MNL} describes why the proposed solution algorithm does not extend to the case of MNL demand. 
In Section~\ref{subsec:finitePcal_convexUcal_complexity}, we develop a theoretical result that characterizes the support of the optimal price distribution in the finite $\Pcal$, convex $\Ucal$ setting. In Section~\ref{subsec:finitePcal_convexUcal_constrained}, we briefly describe how to extend the methodology in Section~\ref{subsec:finitePcal_convexUcal_algorithm} to the case where the price distribution satisfies additional constraints, with further details provided in the ecompanion (Section~\ref{sec:constrained}). Lastly, in Section~\ref{subsec:finitePcal_convexUcal_continuous_vs_discrete}, we discuss our focus on a finite set $\Pcal$ versus a convex price set $\Pcal$. 

\subsection{General solution approach}
\label{subsec:finitePcal_convexUcal_algorithm}

The first general solution scheme that we consider is when $\Pcal$ is a discrete set and the uncertainty set $\Ucal$ is a convex uncertainty set. In this case, if the revenue function $R(\pb, \ub)$ is convex in $\ub \in \Ucal$, then the RRPO problem can be reformulated as follows: 
\begin{align*}
& \max_{\pib \in \Delta_{\Pcal}} \min_{\ub \in \Ucal} \sum_{\pb \in \Pcal} \pi_{\pb} R(\pb, \ub) \\
& = \min_{\ub \in \Ucal} \max_{\pib \in \Delta_{\Pcal}}  \sum_{\pb \in \Pcal} \pi_{\pb} R(\pb, \ub) \\ 
& = \min_{\ub \in \Ucal} \max_{\pb \in \Pcal}  R(\pb, \ub) ,
\end{align*}
where the first equality follows by Sion's minimax theorem, and the second equality follows by the fact that the inner maximum is attained by setting $\pi_{\pb} = 1$ for some $\pb$ and setting $\pi_{\pb'} = 0$ for all $\pb' \neq \pb$. This last problem can be written in epigraph form as 
\begin{subequations}
\begin{alignat}{2}
& \underset{t,\ub}{\text{minimize}} & \quad & t \\
& \text{subject to} & & t \geq R(\pb, \ub), \quad \forall \ \pb \in \Pcal, \label{prob:finiteP_convexU_robust_constraint} \\
& & & \ub \in \Ucal.
\end{alignat}
\label{prob:finiteP_convexU_epigraph}%
\end{subequations}
Problem~\eqref{prob:finiteP_convexU_epigraph} can be solved using constraint generation. In such a scheme, we start with constraint~\eqref{prob:finiteP_convexU_robust_constraint} enforced only at a subset $\hat{\Pcal} \subset \Pcal$, and solve problem~\eqref{prob:finiteP_convexU_epigraph} to obtain a solution $(t,\ub)$. At this solution, we solve the problem $\max_{\pb \in \Pcal} R(\pb, \ub)$, and compare this objective value to the current value of $t$. If it is less than or equal to $t$, we conclude that $(t,\ub)$ satisfies constraint~\eqref{prob:finiteP_convexU_robust_constraint} and terminate with $(t,\ub)$ as the optimal solution. Otherwise, if it is greater than $t$, we have identified a $\pb$ for which constraint~\eqref{prob:finiteP_convexU_robust_constraint} is violated and we add the new constraint to $\hat{\Pcal}$. We then re-solve problem~\eqref{prob:finiteP_convexU_epigraph} with the new $\hat{\Pcal}$ to obtain a new solution $(t,\ub)$ and repeat the process until we can no longer identify any violated constraints. To recover the optimal randomization scheme from the solution of this problem (i.e., the distribution $\pib$), we simply consider the optimal dual variable of each constraint $t \geq R(\pb, \ub)$. 

The viability of this solution approach critically depends on our ability to solve the separation problem $\max_{\pb \in \Pcal} R(\pb, \ub)$ efficiently, and to solve the problem~\eqref{prob:finiteP_convexU_epigraph} efficiently for a fixed subset $\hat{\Pcal} \subset \Pcal$. In what follows, we shall demonstrate that this problem can actually be solved practically for the linear, semi-log and log-log problems. %Unfortunately, as we will discuss in Section~\ref{subsec:finitePcal_convexUcal_MNL}, the reformulation of the RRPO problem as a min-max problem does not apply to the MNL demand model. 

To develop our approaches for the linear, semi-log and log-log models, we will make the following assumption about the price set $\Pcal$, which simply states that $\Pcal$ is a Cartesian product of finite sets of prices for each of the products. 
\begin{assumption}
$\Pcal = \Pcal_1 \times \dots \times \Pcal_I$, where $\Pcal_i$ is a finite subset of $\Rbb_+$ for each $i$.  \label{assumption:Pcal_Cartesian_product}
\end{assumption}

%Lastly, we discuss how one can obtain the optimal randomization scheme from the solution of problem~\eqref{prob:finiteP_convexU_epigraph}. \comebacktothis -- mention that optimal dual variable of constraint $t \leq R(\pb, \ub)$ gives the probability of selecting that $\pb$. Is there a way to do this generically, without any further assumptions on $\Ucal$? 

% START OF LINEAR DEMAND CONVEX UCAL
\subsection{Linear demand model}
\label{subsec:finitePcal_convexUcal_linear}
We begin by showing how our solution approach for convex $\Ucal$ applies to the linear demand model case. Recall that the linear model revenue function is 
\begin{equation}
R(\pb, \ub) = \sum_{i=1}^I p_i (\alpha_i - \beta_i p_i + \sum_{j \neq i} \gamma_{i,j} p_j ).
\end{equation}
For a fixed $\pb$, the function $R(\pb, \ub)$ is linear and therefore convex in $\ub = (\alphab, \betab, \gammab)$. Thus, given a subset $\hat{\Pcal} \subset \Pcal$, the problem~\eqref{prob:finiteP_convexU_epigraph} should be easy to solve, assuming that $\Ucal$ is also a sufficiently tractable convex set. For example, if $\Ucal$ is a polyhedron, then since each constraint~\eqref{prob:finiteP_convexU_robust_constraint} is linear in $\ub$, problem~\eqref{prob:finiteP_convexU_epigraph} would be a linear program.

The separation problem for the linear demand model case is 
\begin{align*}
& \max_{\pb \in \Pcal} \sum_{i=1}^I p_i (\alpha_i - \beta_i p_i + \sum_{j \neq i} \gamma_{i,j} p_j ) \\
& = \max_{\pb \in \Pcal} \sum_{i=1}^I p_i \alpha_i - \sum_{i=1}^I \beta_i p^2_i + \sum_{i=1}^I \sum_{j \neq i} \gamma_{i,j} p_i  p_j 
\end{align*}
Since $\Pcal = \Pcal_1 \times \dots \times \Pcal_I$, we can formulate this as a mixed-integer program. Let $x_{i,t}$ be a binary variable that is 1 if product $i$ has price $t \in \Pcal_i$, and 0 otherwise. Similarly, let $y_{i,j,t_1,t_2}$ be a binary decision variable that is 1 if product $i$ is given price $t_1$ and product $j$ is given price $t_2$ for $i \neq j$, and 0 otherwise. Then the separation problem can be straightforwardly written as 
\begin{subequations}
\begin{alignat}{2}
& \underset{\xb, \yb}{\text{maximize}} & & \sum_{i=1}^I \sum_{t \in \Pcal_i} \alpha_i \cdot t \cdot x_{i,t} - \sum_{i=1}^I \sum_{t \in \Pcal_i} t^2 \cdot \beta_i \cdot x_{i,t} + \sum_{i=1}^I \sum_{ j \neq i} \sum_{t_1 \in \Pcal_i} \sum_{t_2 \in \Pcal_j} \gamma_{i,j} \cdot t_1 \cdot t_2 \cdot y_{i,j,t_1,t_2} \\
& \text{subject to} & \quad & \sum_{t \in \Pcal_i} x_{i,t} = 1, \quad \forall \ i \in [I], \\
& & & \sum_{t_2 \in \Pcal_j} y_{i,j,t_1, t_2} = x_{i,t_1}, \quad \forall \ i, j \in [I], j \neq i, t_1 \in \Pcal_i, \\
& & & \sum_{t_1 \in \Pcal_i} y_{i,j,t_1, t_2} = x_{j,t_2}, \quad \forall \ i, j \in [I], j \neq i, t_2 \in \Pcal_j, \\
& & & x_{i,t} \in \{0,1\}, \quad \forall \ i \in [I], \ t \in \Pcal_i, \\
& & & y_{i,j,t_1,t_2} \in \{0,1\}, \quad \forall \ i, j\in [I], i \neq j, t_1 \in \Pcal_i, t_2 \in \Pcal_j,
\end{alignat}%
\label{prob:linear_MILP}%
\end{subequations}
where the first constraint simply enforces that exactly one price is chosen for each product, while the second and third constraints require that the $y_{i,j,t_1,t_2}$ variables are equal to $x_{i,t_1} \cdot x_{j, t_2}$. In terms of the number of variables and constraints, if we assume that $|\Pcal_1| = \dots = |\Pcal_I| = P$ (all price sets are of the same size), then the number of constraints is $I + 2 I(I-1) \cdot P$ and the number of variables is $I \cdot P + I(I-1) \cdot P^2$.

% END OF LINEAR DEMAND CONVEX UCAL

\subsection{Semi-log demand model}
\label{subsec:finitePcal_convexUcal_semilog}

We will now show how the solution approach we have defined earlier applies to the semi-log demand model. 
%\comebacktothis -- is this a good place for this assumption?? 
Recall that the semi-log revenue function is 
\begin{equation}
R(\pb, \ub) = \sum_{i=1}^I p_i \cdot e^{\alpha_i - \beta_i p_i + \sum_{j\neq i} \gamma_{i,j} p_j}.
\end{equation}
Observe that for a fixed $\pb$, the function $R(\pb, \ub)$ is convex in $\ub$, since it is the nonnegative weighted combination of exponentials of linear functions of $\ub = (\alphab, \betab, \gammab)$. Thus, given a subset $\hat{\Pcal} \subset \Pcal$, solving problem~\eqref{prob:finiteP_convexU_epigraph} should again be ``easy'', assuming also that $\Ucal$ is a sufficiently tractable convex set. (In particular, the function $R(\pb,\ub)$ can be represented using $I$ exponential cones; assuming that $\Ucal$ is also representable using conic constraints, problem~\eqref{prob:finiteP_convexU_epigraph} will thus be some type of continuous conic program.)

We now turn our attention to the separation problem, $\max_{\pb \in \Pcal} R(\pb, \ub)$. Specifically, this problem is
\begin{equation*}
\max_{\pb \in \Pcal} \sum_{i=1}^I p_i \cdot e^{\alpha_i - \beta_i p_i + \sum_{j\neq i} \gamma_{i,j} p_j}.
\end{equation*}
Observe that since the function $f(t) = \log(t)$ is monotonic, the set of optimal solutions remains unchanged if we consider the same problem with a log-transformed objective 
\begin{align}
& \max_{\pb \in \Pcal} \log \left( \sum_{i=1}^I p_i \cdot e^{\alpha_i - \beta_i p_i + \sum_{j\neq i} \gamma_{i,j} p_j} \right) \nonumber \\
& = \max_{\pb \in \Pcal} \log \left( \sum_{i=1}^I  e^{\alpha_i + \log p_i - \beta_i p_i + \sum_{j\neq i} \gamma_{i,j} p_j} \right). \label{eq:semilog_LSE_precursor}
\end{align}
To now further re-formulate this problem, we observe that the objective function can be re-written using the function $g(\yb) = \log( \sum_{i=1}^I e^{y_i})$. The function $g$ is what is known as the \emph{log-sum-exp} function, which is a convex function \citep{boyd2004convex}. More importantly, a standard result in convex analysis is that any proper, lower semi-continuous, convex function is equivalent to its \emph{biconjugate} function, which is the convex conjugate of its convex conjugate \citep{rockafellar1970convex}. For the log-sum-exp function, this in particular means that $g(\yb)$ can be written as 
\begin{align*}
g(\yb) = \max_{ \mub \in \Delta_{[I]} } \{ \mub^T \yb - \sum_{i=1}^I \mu_i \log \mu_i \}.
\end{align*}
The function $h(x) = x \log x$ is the negative entropy function \citep{boyd2004convex}, and is a convex function; thus, the function inside the $\max\{ \cdot \}$ is a linear function minus a sum of convex functions, and is a concave function. 

For our problem, this means that \eqref{eq:semilog_LSE_precursor} can be re-written as
\begin{align}
& \max_{\pb \in \Pcal} \log R(\pb, \ub) \nonumber \\
& = \max_{\pb \in \Pcal} \log \left( \sum_{i=1}^I  e^{\alpha_i + \log p_i - \beta_i p_i + \sum_{j\neq i} \gamma_{i,j} p_j} \right) \nonumber \\
& = \max_{\pb \in \Pcal} \max_{ \mub \in \Delta_{[I]} } \left\{ \sum_{i=1}^I \mu_i (\alpha_i + \log p_i - \beta_i p_i + \sum_{j\neq i} \gamma_{i,j} p_j) - \sum_{i=1}^I \mu_i \log \mu_i \right\} \nonumber \\
& = \max_{\pb \in \Pcal,\  \mub \in \Delta_{[I]} } \left\{ \sum_{i=1}^I \mu_i (\alpha_i + \log p_i - \beta_i p_i + \sum_{j\neq i} \gamma_{i,j} p_j) - \sum_{i=1}^I \mu_i \log \mu_i \right\}.\label{eq:semilog_biconjugate_abstract}
\end{align}
To further reformulate this problem, we now make use of Assumption~\ref{assumption:Pcal_Cartesian_product}, which states that $\Pcal$ is the Cartesian product of finite sets. Let us introduce a new binary decision variable $x_{i,t}$ which is 1 if product $i$'s price is set to price $t \in \Pcal_i$, and 0 otherwise. Using this new decision variable, observe that we can replace $p_i$ wherever it occurs with $\sum_{t \in \Pcal_i} t \cdot x_{i,t}$. We can also similarly replace $\log p_i$ with $\sum_{t \in \Pcal_i} \log t \cdot x_{i,t}$. Therefore, problem~\eqref{eq:semilog_biconjugate_abstract} can be further reformulated as 
\begin{subequations}
\begin{alignat}{2}
& \underset{\xb, \mub}{\text{maximize}} & \quad & \sum_{i=1}^I \mu_i \left(\alpha_i + \sum_{t \in \Pcal_i} \log t \cdot x_{i,t} - \beta_i \cdot \sum_{t \in \Pcal_i} t  \cdot x_{i,t} + \sum_{j\neq i} \gamma_{i,j} \sum_{t \in \Pcal_j} t \cdot x_{j,t} \right) - \sum_{i=1}^I \mu_i \log \mu_i \\ 
& \text{subject to} & & \sum_{i=1}^I \mu_i = 1, \\
& & &  \sum_{t \in \Pcal_i} x_{i,t} = 1, \quad \forall \ i \in [I], \\
& & & x_{i,t} \in \{0,1\}, \quad \forall \ i \in [I], \ t \in \Pcal_i, \\
& & & \mu_i \geq 0, \quad \forall \ i \in [I].
\end{alignat}
\end{subequations}
This last problem is \emph{almost} a mixed-integer convex program: as noted earlier, the expression $- \sum_{i=1}^I \mu_i \log \mu_i$ is concave in $\mu$. The main wrinkle is the presence of the bilinear terms in the objective function, specifically terms of the form $\mu_i \cdot x_{j,t}$. Fortunately, we can circumvent this difficulty by introducing a new decision variable, $w_{i,j,t}$, which is the linearization of $\mu_i \cdot x_{j,t}$, for each $i, j \in [I]$, $t \in \Pcal_j$. By adding this new decision variable and additional constraints, we arrive at our final formulation, which is a mixed-integer convex program. 
\begin{subequations}
\begin{alignat}{2}
& \underset{\mub, \wb, \xb}{\text{maximize}} & \quad & \sum_{i=1}^I \mu_i \alpha_i + \sum_{i=1}^I \sum_{t \in \Pcal_i} \log t \cdot w_{i,i,t} - \sum_{i=1}^I \beta_i \cdot \sum_{t \in \Pcal_i} t \cdot w_{i,i,t} + \sum_{i=1}^I \sum_{j\neq i} \gamma_{i,j} \cdot (\sum_{t \in \Pcal_j} t \cdot w_{i,j,t}) - \sum_{i=1}^I \mu_i \log \mu_i \\ 
& \text{subject to} & &  \sum_{t \in \Pcal_j} w_{i,j,t} = \mu_i, \quad \forall \ i \in [I], \ j \in [I], \label{prob:semilog_MIECP_w_sums_to_mu} \\
& & & \sum_{i=1}^I w_{i,j,t} = x_{j,t}, \quad \forall \ j \in [I], \ t \in \Pcal_j, \label{prob:semilog_MIECP_w_sums_to_x} \\
& & & \sum_{i=1}^I \mu_i = 1, \\
& & & \sum_{t \in \Pcal_i} x_{i,t} = 1, \quad \forall \ i \in [I], \\
& & & w_{i,j,t} \geq 0, \quad \forall \ i \in [I],\ j \in [I], t \in \Pcal_j, \label{prob:semilog_MIECP_w_nonnegative} \\
& & & x_{i,t} \in \{0,1\}, \quad \forall \ i \in [I], \ t \in \Pcal_i, \\
& & & \mu_i \geq 0, \quad \forall \ i \in [I]. \label{prob:semilog_MIECP_mu_nonnegative}
\end{alignat}
\label{prob:semilog_MIECP}%
\end{subequations}
There are a few important points to observe about this formulation. First, note that because the $\mu_i$'s sum to 1 over $i$, and the $x_{j,t}$'s are binary and sum to 1 over $t \in \Pcal_j$ for any $j$, then ensuring that $w_{i,j,t} = \mu_i \cdot x_{j,t}$ can be done simply through constraints~\eqref{prob:semilog_MIECP_w_sums_to_mu} and \eqref{prob:semilog_MIECP_w_sums_to_x}. This is different from the usual McCormick envelope-style linearization technique, which in this case would involve the four inequalities:
\begin{align}
w_{i,j,t} & \leq x_{j,t}, \\
w_{i,j,t} & \leq \mu_i, \\
w_{i,j,t} & \geq x_{j,t} + \mu_i - 1, \\
w_{i,j,t} & \geq 0,
\end{align}
for every $i \in [I]$, $j \in [I]$, $t \in \Pcal_j$. It is not difficult to show that these constraints are implied by constraints~\eqref{prob:semilog_MIECP_w_sums_to_mu}, \eqref{prob:semilog_MIECP_w_sums_to_x} and \eqref{prob:semilog_MIECP_w_nonnegative}. 

Second, at the risk of belaboring the obvious, the optimal objective value of problem~\eqref{prob:semilog_MIECP} is the value of $\max_{\pb \in \Pcal} \log R(\pb, \ub)$, where $R$ is the semi-log revenue function. Upon solving problem~\eqref{prob:semilog_MIECP} to obtain the objective value $Z'$, we can obtain the optimal objective value of the untransformed problem $\max_{\pb \in \Pcal} R(\pb, \ub)$ as $e^{Z'}$. 

Third, this formulation is notable because, to our knowledge, this is the first exact mixed-integer convex formulation of the nominal multi-product pricing problem under semi-log demand and a price set defined as the Cartesian product of finite sets. To date, virtually all research that has considered solving this type of problem in the marketing and operations management literatures has involved heuristics (see, for example, Section EC.3 of \citealt{misic2020optimization}, which solves log-log and semi-log multi-product pricing problems for a collection of stores using local search). From this perspective, although we developed this formulation as part of the overall solution approach for the RRPO problem, we believe it is of more general interest.

Building on the previous point, problem~\eqref{prob:semilog_MIECP} can be formulated as a mixed-integer exponential cone program. We state the exact formulation in Section~\ref{sec:miecp_explicit_models} of the ecompanion; assuming that $|\Pcal_1| = \dots = |\Pcal_I| = P$, this formulation has $I^2 \cdot P + I \cdot P + 2I$ variables, $I^2 + I \cdot P + I + 1$ linear constraints and $I$ exponential cone constraints. Such problems are garnering increasing attention from the academic and industry sides. In particular, since 2019, the MOSEK solver \citep{mosek} supports the exponential cone and can solve mixed-integer conic programs that involve the exponential cone to global optimality. Although the solution technology for mixed-integer conic programs is not as developed as that of mixed-integer linear programs (as exemplified by state-of-the-art solvers such as Gurobi and CPLEX), it is reasonable to expect that these solvers will continue to improve and allow larger and larger problem instances to be solved to optimality in the future.

Lastly, we comment that the same reformulation technique used above -- taking the logarithm, replacing the log-sum-exp function with its biconjugate, and then linearizing the products of the binary decision variables and the probability mass function values (the $\mu_i$ variables) that arise from the biconjugate -- can also be used to derive an exact formulation of the deterministic robust price optimization problem. By taking the same approach, one obtains a max-min-max problem, and one can use Sion's minimax theorem again to swap the inner maximization over $\mub$ with the minimization over $\ub$ to obtain a robust counterpart that can then be further reformulated using duality or otherwise solved using delayed constraint generation. We provide the details of this derivation in Section~\ref{subsec:convexUcal_DRPO_semilog} of the ecompanion.

\subsection{Log-log demand model}
\label{subsec:finitePcal_convexUcal_loglog}

To now show how the solution scheme in Section~\ref{subsec:finitePcal_convexUcal_algorithm} applies to the log-log approach, we again recall the form of the log-log revenue function:
\begin{align}
R(\pb, \ub) & = \sum_{i=1}^I p_i \cdot e^{\alpha_i - \beta_i \log p_i + \sum_{j\neq i} \gamma_{i,j} \log p_j} \\
& = \sum_{i=1}^I e^{\alpha_i + \log p_i - \beta_i \log p_i + \sum_{j\neq i} \gamma_{i,j} \log p_j}
\end{align}
Using the same biconjugate trick as with the semi-log approach, we can show that
\begin{align}
& \max_{\pb \in \Pcal} \log R(\pb, \ub) \nonumber \\
& = \max_{\pb \in \Pcal} \log \left( \sum_{i=1}^I  e^{\alpha_i + \log p_i - \beta_i \log p_i + \sum_{j\neq i} \gamma_{i,j} \log p_j} \right) \nonumber \\
& = \max_{\pb \in \Pcal} \max_{\mub \in \Delta_{[I]}} \left\{ \sum_{i=1}^I \mu_i (\alpha_i + \log p_i - \beta_i \log p_i + \sum_{j\neq i} \gamma_{i,j} \log p_j) - \sum_{i=1}^I \mu_i \log \mu_i \right\} \label{prob:loglog_maxmax} 
%& = \max_{\pb \in \Pcal, \mub \in \Delta_{[I]}} \left\{ \sum_{i=1}^I \mu_i (\alpha_i + \log p_i - \beta_i \log p_i + \sum_{j\neq i} \gamma_{i,j} \log p_j) - \sum_{i=1}^I \mu_i \log \mu_i \right\} 
\end{align}
If we now invoke Assumption~\ref{assumption:Pcal_Cartesian_product}, then we can introduce the same decision variables $x_{i,t}$ and $w_{i,j,t}$ as in problem~\eqref{prob:semilog_MIECP} to obtain a mixed-integer convex formulation of the log-log price optimization problem, which has the same feasible region as the semi-log formulation~\eqref{prob:semilog_MIECP}:
\begin{subequations}
\begin{alignat}{2}
& \underset{\wb, \xb, \mub}{\text{maximize}} & \quad & \sum_{i=1}^I \mu_i \alpha_i + \sum_{i=1}^I \sum_{t \in \Pcal_i} \log t \cdot w_{i,i,t} - \sum_{i=1}^I \beta_i \cdot \sum_{t \in \Pcal_i} \log t \cdot w_{i,i,t} + \sum_{i=1}^I \sum_{j\neq i} \gamma_{i,j} \sum_{t \in \Pcal_j} \log t \cdot w_{i,j,t} - \sum_{i=1}^I \mu_i \log \mu_i \\ 
& \text{subject to} & &  \text{constraints~\eqref{prob:semilog_MIECP_w_sums_to_mu} -- \eqref{prob:semilog_MIECP_mu_nonnegative}}. 
%\sum_{t \in \Pcal_j} w_{i,j,t} = \mu_i, \quad \forall \ i \in [I], \ j \in [I], \label{prob:loglog_MIECP_w_sums_to_mu} \\
%& & & \sum_{i=1}^I w_{i,j,t} = x_{j,t}, \quad \forall \ j \in [I], \ t \in \Pcal_j, \label{prob:loglog_MIECP_w_sums_to_x} \\
%& & & \sum_{i=1}^I \mu_i = 1, \\
%& & & \sum_{t \in \Pcal_i} x_{i,t} = 1, \quad \forall \ i \in [I], \\
%& & & w_{i,j,t} \geq 0, \quad \forall \ i \in [I],\ j \in [I], t \in \Pcal_j, \label{prob:loglog_MIECP_w_nonnegative} \\
%& & & x_{i,t} \in \{0,1\}, \quad \forall \ i \in [I], \ t \in \Pcal_i, \\
%& & & \mu_i \geq 0, \quad \forall \ i \in [I]. \label{prob:loglog_MIECP_mu_nonnegative}
\end{alignat}
\label{prob:loglog_MIECP}%
\end{subequations}
While the feasible region of problem~\eqref{prob:loglog_MIECP} is the same as that of  \eqref{prob:semilog_MIECP}, the objective function of \eqref{prob:loglog_MIECP} is different. Just like problem~\eqref{prob:semilog_MIECP}, problem~\eqref{prob:loglog_MIECP} can be written as a mixed-integer exponential cone program (see Section~\ref{sec:miecp_explicit_models} of the ecompanion for the exact formulation), and the number of variables, linear constraints and exponential cone constraints is the same as that of \eqref{prob:semilog_MIECP}. Similarly, to the best of our knowledge, this is the first exact mixed-integer convex formulation of the log-log multi-product price optimization problem under a Cartesian product price set. Lastly, just like the semi-log problem~\eqref{prob:semilog_MIECP}, one can easily modify the formulation to obtain an exact formulation of the deterministic robust price optimization problem under log-log demand (see Section~\ref{subsec:convexUcal_DRPO_loglog} of the ecompanion). 

We note that the log-log separation problem has an interesting property, which is that there exist optimal solutions that are extreme, in the sense that each product's price is set to either its lowest or highest allowable price. This property is formalized in the following proposition (see Section~\ref{proof_proposition:loglog_extremal} for the proof).
\begin{proposition}
Suppose that Assumption~\ref{assumption:Pcal_Cartesian_product} holds. There exists an optimal solution $(\mub, \pb)$ of problem~\eqref{prob:loglog_maxmax} such that for each $i \in [I]$, either $p_i = \min \Pcal_i$ or $p_i = \max \Pcal_i $. \label{proposition:loglog_extremal}
\end{proposition}
The proof follows because when $\mub$ is fixed, the log-transformed separation problem \eqref{prob:loglog_maxmax} is linear in $(\log p_1,\dots, \log p_I)$; thus, given an optimal but possibly non-extremal solution $(\mub^*, \pb^*)$, we can optimize over $\pb$ with $\mub^*$ fixed to obtain a new solution $(\mub^*, \pb')$ that must also be optimal, and for which $\pb'$ is extremal.

%\subsection{MNL demand model}
%\label{subsec:finitePcal_convexUcal_MNL}
%
%While the general solution approach described in Section~\ref{subsec:finitePcal_convexUcal_algorithm} nicely applies to both the semi-log and log-log demand models, it is sadly not applicable to the MNL demand model. In particular, in the MNL demand model, the revenue function is 
%\begin{equation}
%R(\pb, \ub) = \frac{ \sum_{i=1}^I p_i e^{\alpha_i - \beta_i p_i }}{1 + \sum_{i=1}^I e^{\alpha_i - \beta_i p_i  }}
%\end{equation}
%The roadblock that one encounters in considering the MNL revenue function is that in general, it is not quasiconvex in $\ub = (\alphab, \betab)$, so the transformation of the RRPO problem to the problem $\min_{\ub \in \Ucal} \max_{\pb \in \Pcal} R(\pb, \ub)$ does not go through using the same argument described in Section~\ref{subsec:finitePcal_convexUcal_algorithm}. Notwithstanding this, it is also not clear whether problem~\eqref{prob:finiteP_convexU_epigraph} is tractable, as the constraint $t \geq R(\pb,\ub)$ will in general not be a convex constraint in $\ub$. Nevertheless, we will see that when $\Ucal$ is a discrete set, then it does becomes possible to solve the MNL-based RRPO problem, thanks to a previously developed exponential cone formulation of MNL choice probabilities \citep{akcakus2021exact}. 

\subsection{Complexity of randomized robust pricing with finite $\Pcal$}
\label{subsec:finitePcal_convexUcal_complexity}

In this section, we discuss the \emph{complexity} of randomized robust pricing strategies. By complexity, we refer the size of the support of the resulting distribution $\pib$ that solves the RRPO problem. This is an important quantity to analyze for a number of reasons. First, a randomized pricing strategy that randomizes over a small number of price vectors will be easier to implement compared to one that randomizes over a large number of price vectors. Second, from the perspective of the customer, randomized pricing can potentially raise concerns around fairness and price discrimination, as some customers may be shown more favorable prices for some products than other customers. In this regard, a randomized pricing strategy that randomizes over a smaller set of price vectors could potentially be regarded as more fair than one that randomizes over a larger set of price vectors. 

Our key result in this section is an upper bound on the complexity of randomization in the context of RRPO with a finite $\Pcal$ and convex $\Ucal$.

\begin{theorem}
Suppose $R(\pb, \ub)$ is convex in $\ub$ for any $\pb \in \Pcal$, that $\Ucal$ is a non-empty convex compact subset of $\Rbb^M$, and that $\Pcal$ is a finite set. Then there exists an optimal solution $\pib^*$ to the RRPO problem $\max_{\pib \in \Delta_{\Pcal}} \min_{\ub \in \Ucal} \sum_{\pb \in \Pcal} \pi_{\pb} R(\pb, \ub)$, such that $\pib^*$ has a support of at most $M + 1$, i.e., the number of price vectors $\pb$ for which $\pi_{\pb} > 0$ is at most $M + 1$. \label{theorem:complexity_convexUcal}
\end{theorem}

The proof of this result is provided in Section~\ref{proof_theorem:complexity_convexUcal} of the electronic companion. The proof of this result is rather involved, and requires the careful application of results from the paper of  \cite{calafiore2005uncertain} on support constraints. In particular, in a generic convex program with a linear objective function, a \emph{support constraint} is a constraint that, once removed, results in a change in the optimal objective value of the problem. A key result in \cite{calafiore2005uncertain} is that any convex program with an $n$-dimensional decision variable has at most $n$ support constraints. In the context of the dual, any active constraint in the dual corresponds to a potentially non-zero primal variable (i.e., a $\pi_{\pb}$ for some $\pb \in \Pcal$) by complementary slackness. However, while every support constraint is an active constraint at any optimal solution, an active constraint at a given optimal solution need not be a support constraint. The proof of our result thus requires (1) iteratively reducing the dual problem~\eqref{prob:finiteP_convexU_epigraph} by removing constraints one at a time to obtain an equivalent dual problem with a smaller number of constraints, (2) asserting that all active constraints in this equivalent dual problem are support constraints, and (3) asserting that the number of such constraints is no more than $M+1$. 

Besides the significance of this result in our randomized pricing context, the proof technique is general and could potentially be applied in the context of other convex programs with a linear objective function. On this point, we note that the analogous result for the finite $\Pcal$ and finite $\Ucal$ case (Section~\ref{subsec:finitePcal_finiteUcal_complexity} of the ecompanion) is considerably simpler to establish, because the RRPO problem can be written as a linear program and one can invoke standard properties of basic feasible solutions to bound the support of $\pib$. In contrast, problem~\eqref{prob:finiteP_convexU_epigraph} is not a linear program, and to the best of our knowledge, we are not aware of any  characterization of the number of active constraints in general convex programs that one could trivially invoke to obtain the above result.

This result immediately allows us to understand the complexity of randomization for the three demand models that we study. In particular, the dimension $M$ for the linear, semi-log and log-log demand models is at most equal to $I + I + I(I-1) = I^2 + I$, where $I$ is the number of products. Thus, Theorem~\ref{theorem:complexity_convexUcal} implies that any solution to the RRPO problem when $\Pcal$ is finite and $\Ucal$ is convex will not have a support of size greater than $I^2 + I + 1$, which scales polynomially in the number of products. In practice, the complexity of a randomized robust pricing strategy could be much smaller than this. For example, in Section~\ref{subsec:results_orangejuice}, we show an example of an optimal solution to the RRPO problem for the log-log demand model with $I = 11$ products, where the complexity of randomization is only 6.  

\subsection{Constrained price distributions}
\label{subsec:finitePcal_convexUcal_constrained}

Besides the complexity of randomization, it is also interesting to consider solving the RRPO problem with constraints on the set of allowable distributions. By incorporating constraints on the set of distributions that one randomizes over, one can more directly control how prices are randomized and control the variability of the prices. It turns out that moment constraints on the distribution $\pib$ can be easily incorporated through a suitable modification of the dual problem~\eqref{prob:finiteP_convexU_epigraph}. Due to space constraints, we relegate our discussion of this to Section~\ref{sec:constrained} of the ecompanion, which also includes a small set of numerical experiments to show that the price of incorporating such constraints is not too large. 

\subsection{Convex versus finite price sets}
\label{subsec:finitePcal_convexUcal_continuous_vs_discrete}

We close Section~\ref{sec:finitePcal_convexUcal} by discussing our assumption of $\Pcal$ being finite. Thus far, we have assumed that $\Pcal$ is a finite set, rather than a convex set; a natural question is whether it is possible to generalize or modify the approach to handle the case where $\Pcal$ is a convex set.

From a technical standpoint, it is not straightforward to handle the case where $\Pcal$ is a convex set. To see why, observe that the RRPO problem can be reformulated in the following way using Sion’s minimax theorem:
\begin{align*}
& \max_{F \in \Fcal} \min_{\ub \in \Ucal} \int R(\pb, \ub) \, dF(\pb) \\ 
& = \min_{\ub \in \Ucal} \max_{F \in \Fcal} \int R(\pb, \ub) \, dF(\pb) \\ 
& = \min_{\ub \in \Ucal} \max_{\pb \in \Pcal} R(\pb, \ub)
\end{align*}
where we make the additional assumption that $R(\pb,\ub)$ is convex in $\ub$ and that $\Ucal$ is a convex compact set. This last problem can be rewritten as 
\begin{subequations}
\begin{alignat}{2}
& \underset{t, \ub}{\text{minimize}} & \quad &  t \\
& \text{subject to} & &  t \geq R(\pb,\ub), \quad \forall \ \pb \in \Pcal, \\
& & & \ub \in \Ucal,
\end{alignat}
\end{subequations}
which is similar to problem~\eqref{prob:finiteP_convexU_epigraph}. The difficulty now arises from $\Pcal$ being a convex set. For the log-log and semi-log demand models, when $\Pcal$ is discrete, the separation problem $\max_{\pb \in \Pcal} R(\pb,\ub)$ can be reformulated using the biconjugate technique, because one can linearize products of binary variables and the biconjugate variables $(\mu_1,\dots, \mu_I)$ using standard integer programming ideas. For the present case where $\Pcal$ is convex, applying the biconjugate reformulation immediately leads to a roadblock because now one has to solve a problem with bilinear terms. It is not clear whether an exact reformulation exists in this case. Although one could solve the separation problem heuristically, this would not guarantee that one would obtain a provably optimal solution. 

Beside this technical issue, we believe that the restriction to a finite $\Pcal$ is not overly restrictive, and is compatible with practice. First, in most retail scenarios, prices are naturally discrete due to how currencies are denominated. For example, the smallest denomination for US dollars is \$0.01 (the cent), and one does not see prices that are quoted to a finer precision than a hundredth of a dollar in typical retail settings (e.g., at grocery stores). Second, even beside this natural discretization, some retailers will usually consider coarser discretizations due to marketing and behavioral reasons. For example, in marketing it is well-recognized that consumers often respond more positively to prices that end in the number 9 (for example, a price such as \$149 versus \$150; see \citealt{anderson2003effects}). Relatedly, there is research in the marketing literature that suggests that customers anchor on the first few significant digits of the price and de-emphasize the remaining digits \citep{bizer2005direct,shlain2022more}; for a product priced in the thousands, this suggests that one should use prices such as \$1999 or \$1499, and not prices such as \$1937 or \$1415. Similarly, other research suggests that customers regard round prices as indicative of higher quality \citep{stiving2000price}. Finally, it is unlikely that a retailer will consider candidate prices that span multiple orders of magnitude for a product. A computer retailer like Apple might price laptops in the thousands of dollars, but likely will not consider prices in the tens of thousands or hundreds of thousands of dollars. This also limits how large a discrete set of price vectors $\Pcal$ can be. 

Aside from the practical considerations we have described, we note that it is also possible to show that one can approximate a continuous RRPO problem with a discrete RRPO problem, and to understand theoretically how fine of a discretization one needs to consider. In Section~\ref{sec:approximation_convex_discrete} of the ecompanion, we show theoretically that for any RRPO problem over a continuous price set, there exists a finite price set such that the objective value of the resulting discrete RRPO problem is arbitrarily close to the true continuous RRPO problem. In particular, the proof of this result implies that restricting the price of each product to a uniform grid, where the coarseness of the grid is informed by the Lipschitz constant of the revenue function, is sufficient to obtain a good approximation of the true continuous problem.

\section{Solution method for finite $\Pcal$, finite $\Ucal$}
\label{sec:finitePcal_finiteUcal_body}

In addition to the case where $\Ucal$ is convex, we also consider the case where $\Ucal$ is a finite discrete set. Due to page limitations, our presentation of our solution method for this case is relegated to Section~\ref{sec:finitePcal_finiteUcal} of the ecompanion. At a high level, the foundation of our approach is double column generation, which alternates between solving the primal version of the RRPO problem, which is $\max_{\pib \in \Delta_{\Pcal}} \min_{\ub \in \Ucal} \sum_{\pb \in \Pcal} \pi_{\pb} R(\pb, \ub)$, and the dual version of the RRPO problem, which is $\min_{\lambdab \in \Delta_{\Ucal}} \max_{\pb \in \Pcal} \sum_{\ub \in \Ucal} \lambda_{\ub} R(\pb, \ub)$. In each iteration, we solve the primal problem with $\Pcal$ replaced by a subset $\hat{\Pcal} \subseteq \Pcal$, where we use constraint generation to handle the inner minimization over $\ub \in \Ucal$; this gives rise to a finite set of uncertainty realizations $\hat{\Ucal} \subseteq \Ucal$. We then solve the dual problem with $\Ucal$ replaced by $\hat{\Ucal}$, where we use constraint generation to handle the inner maximization over $\pb \in \Pcal$, which gives rise to a finite set of price vectors $\hat{\Pcal} \subseteq \Pcal$. At each step of the algorithm, the objective value of the primal problem restricted to $\hat{\Pcal}$ is a lower bound on the true optimal objective, while the objective value of the dual problem restricted to $\hat{\Ucal}$ is an upper bound on the optimal objective; the algorithm terminates when these two bounds are equal or are otherwise within a pre-specified tolerance. 

To implement this approach for the demand models that we consider, one needs to be able to solve the primal separation problem (solve $\min_{\ub \in \Ucal} \sum_{\pb \in \hat{\Pcal}} \pi_{\pb} R(\pb,\ub)$) and the dual separation problem (solve $\max_{\pb \in \Pcal} \sum_{\ub \in \hat{\Ucal}} \lambda_{\ub} \cdot R(\pb, \ub)$). We show how both of these problems can be reformulated as mixed-integer exponential cone programs for the semi-log and log-log demand models, and as mixed-integer linear programs for the linear demand model.

\section{Numerical experiments}
\label{sec:results}

In this section, we conduct several sets of experiments involving synthetic problem instances to understand the tractability of the RRPO approach and the improvement in worst-case revenue of the randomized robust pricing strategy over the deterministic robust pricing strategy. In Section~\ref{subsec:results_convexUcal}, we consider instances involving the linear, semi-log and log-log models where the uncertainty set $\Ucal$ is a convex set. We also consider instances involving the linear, semi-log and log-log models where the uncertainty set $\Ucal$ is a discrete set; due to space considerations, these results are presented in Section~\ref{subsec:results_discreteUcal_loglog_semilog} of the ecompanion. %In Section~\ref{subsec:results_discreteUcal_MNL}, we consider instances involving the MNL model where the uncertainty set $\Ucal$ is discrete. 
In Section~\ref{subsec:results_orangejuice}, we consider log-log and semi-log robust price optimization instances derived from a real data set on sales of orange juice products from a grocery store chain. 
Lastly, in Section~\ref{subsec:R1_results_orangejuice}, we consider a data-driven experiment using the same orange juice data set that illustrates the benefit of randomized robust price strategies over deterministic robust and non-robust price prescriptions in out-of-sample performance.

All of our code is implemented in the Julia programming language \citep{bezanson2017julia}. All optimization models are implemented using the JuMP package \citep{lubin2015computing}. All linear and mixed-integer linear programs are solved using Gurobi \citep{gurobi} and all mixed-integer exponential cone programs are solved using Mosek \citep{mosek}, with a maximum of 8 threads per program. All of our experiments are conducted on Amazon Elastic Compute Cloud (EC2), on a single instance of type {\tt m6a.48xlarge} (AMD EPYC 7R13 processor, with 192 virtual CPUs and 768 GB of memory).

\subsection{Experiments with convex $\Ucal$ and linear, log-log and semi-log demand models} 
\label{subsec:results_convexUcal}

In our first set of experiments, we consider the log-log and semi-log demand models, and specifically consider a $L1$-norm uncertainty set $\Ucal$:
\begin{equation}
\Ucal=\left\{\ub=(\alphab,\betab,\gammab) \ \vline \  \| \tilde{\ub} \|_1\leq\theta; \ \tilde{\ub}_k=\frac{\ub_k-\ub_{0k}}{\ub_{0k}},\  \forall k\in\{1,..,I+I^2\}  \right\}, \label{eq:L1_norm_Ucal}
\end{equation}
where $\ub_0$ is the vector of nominal values of the uncertain parameters $\ub = (\alphab, \betab,\gammab)$, and $\theta$ is the budget of the uncertainty set. Informally, the rationale for this type of uncertainty set is that the nominal demand model parameter estimates are almost certainly wrong when they are estimated from a small sample of data. The perspective that one takes with this type of uncertainty set is that while there could be error in any individual demand model parameter, we expect that the aggregate relative error should not be too large (for example, we do not expect a large number of demand model parameters to exhibit large relative errors). We note that this type of budget uncertainty set is common in the robust optimization literature (see for example \citealt{bertsimas2004price}, \citealt{bertsimas2006robust}).

For each of the three demand models (linear, semi-log and log-log), the number of products $I$ varies in $\{5, 10, 15, 20\}$. For each value of $I$, we generate 24 random instances, where the values of $\alphab, \betab, \gammab$ are independently randomly generated as follows:
\begin{enumerate}
%\item \emph{Linear demand}. Each $\alpha_i \sim \Uniform(200, 300)$, $\beta_i \sim \Uniform(5,15)$, $\gamma_{i,j} \sim \Uniform(-0.1,+0.1)$.
%\item \emph{Semi-log demand}. Each $\alpha_i \sim \Uniform(4,7)$, $\beta_i \sim \Uniform(1,1.5)$, $\gamma_{i,j} \sim \Uniform(-0.4, +0.4)$.
%\item \emph{Log-log demand}. Each $\alpha_i \sim \Uniform(10,14)$, $\beta_i \sim \Uniform(1,2)$, $\gamma_{i,j} \sim \Uniform(-0.6, +0.6)$. 
\item \emph{Linear demand}. Each $\alpha_i \sim \Uniform(100, 200)$, $\beta_i \sim \Uniform(5,15)$, $\gamma_{i,j} \sim \Uniform(-0.1,+0.1)$.
\item \emph{Semi-log demand}. Each $\alpha_i \sim \Uniform(8,10)$, $\beta_i \sim \Uniform(1.5,2)$, $\gamma_{i,j} \sim \Uniform(-0.5, +0.5)$.
\item \emph{Log-log demand}. Each $\alpha_i \sim \Uniform(10,14)$, $\beta_i \sim \Uniform(1.5,2)$, $\gamma_{i,j} \sim \Uniform(-0.8, +0.8)$. 
\end{enumerate}
For each product $i \in [I]$, we set $\Pcal_i = \{1,2,3,4,5\}$. 

%For the uncertainty set $\Ucal$, the budget parameter $\theta$ varies in $\{0.1, 0.5, 1, 1.5, 2\}$ for each instance. %The nominal value of each uncertain parameter is randomly generated. For the instances instances with semi-log demand, $\alphab_0$ is generated from $\Uniform[5, 8]$, $\betab_0$ is generated from $\Uniform[0.5,1.6]$, and $\gammab_0$ is generated from $\Uniform[0.005,0.5]$. In the instances with log-log demand, $\alphab_0$ is generated from $\Uniform[10,15]$, $\betab_0$ is generated from $\Uniform[1.5,2.5]$, and $\gammab_0$ is generated from $\Uniform[-0.8,0.8]$. 

For each instance, we solve the nominal problem, the DRPO problem and the RRPO problem. For DRPO and RRPO, we vary the budget parameter $\theta$ that defines the uncertainty within the set $\{0.1, 0.5, 1, 1.5, 2\}$. To solve the RRPO problem for each instance, we execute the constraint generation solution algorithm described in Section~\ref{sec:finitePcal_convexUcal}. For instances with log-log demand, we take advantage of Proposition~\ref{proposition:loglog_extremal} and thus simplify the price set $\Pcal$ to contain the highest and lowest price levels for each product only. To solve the DRPO problem, we formulate it as either a mixed-integer linear program (for linear demand) or a mixed-integer exponential cone program (for semi-log and log-log demand) via the log-sum-exp biconjugate-based technique described in Section~\ref{sec:convexUcal_DRPO} of the ecompanion, and use standard LP duality techniques to reformulate the objective function of the resulting problem (formulation~\eqref{prob:DRPO_semilog_MIECP} and formulation~\eqref{prob:DRPO_loglog_MIECP} in Section~\ref{sec:convexUcal_DRPO}). Due to the prohibitive computation times that we encountered for the DRPO problem with log-log and semi-log demand, we impose a computation time limit of 20 minutes. From our experimentation with the DRPO problem for log-log and semi-log, it is often the case that an optimal or nearly optimal solution is found early on, and the bulk of the remaining computation time, which can be in the hours, is required by Mosek to prove optimality and close the gap.   %To reduce the computation time of DRPO with log-log demand, we again use the extreme price property to simplify $\Pcal$ to contain extreme price levels only. 
Finally, to solve the nominal problem for each instance, we either solve it as a mixed-integer linear program (for linear demand) or otherwise use the biconjugate technique to formulate it as a mixed-integer exponential cone program (for log-log and semi-log demand).

We present the objective value as well as the computation time of each RRPO, DRPO and nominal problem. We additionally compute several other metrics. We compute $\Ebb[R(\pb^*_{\RR},\ub_0)]$, which is the expected revenue of the randomized RPO solution assuming that the nominal parameter values are realized. We also compute $R(\pb^*_{\DR},\ub_0)$, the nominal revenue of DRPO solution, and $Z_{\Nom,\WC} = \min_{\ub\in\Ucal} R(\pb^*_{\Nom},\ub)$, the worst-case revenue of the nominal solution. We use the following metric, RI (relative improvement), to show the benefit of randomized strategy in robust price optimization:
\begin{equation}
	\text{RI} = (Z_{\RR}^* - Z_{\DR}^* )/Z_{\DR}^* \times 100\%
\end{equation}
We additionally record the solution time for the randomized robust, deterministic robust and nominal approaches, which we denote by $t_{\RR}$, $t_{\DR}$ and $t_{\Nom}$, respectively. For each metric, we compute its average over the 24 instances for each value of $I$ and $\theta$.

Tables~ \ref{table:result_convexU_linear}, \ref{table:result_convexU_semi_log} and \ref{table:result_convexU_log_log} shows the results for the linear, semi-log and log-log demand models, respectively. For linear demand, we find that the improvement by randomized robust pricing over deterministic robust pricing is modest; the largest average improvement is 8.49\% for $I = 5$, $\theta = 2$. We note that we experimented with other forms of uncertainty sets and choices of the nominal parameter values for the linear demand model, but we generally did not encounter large improvements of the same size as we did for the other two demand models.

\begin{table}[ht]
	\centering \small
	\begin{tabular}{llllllllllllll}
		\toprule
		$I$ & $\theta$ & $t_{\RR}$ & $Z_{\RR}^*$ & $\Ebb[R(\pb^*_{\RR},\ub_0)]$ & $t_{\DR}$ & $Z^*_{\DR}$ & RI(\%) & $R(\pb^*_{\DR},\ub_0)$ & $t_{\Nom}$ & $Z_{\Nom}^*$ & $Z_{\Nom,\WC}$ \\ 
		\midrule
		5 & 0.1 & 0.07 & 2382.23 & 2473.30 & 0.03 & 2382.23 & 0.00 & 2473.30 & 0.03 & 2473.30 & 2382.23 \\   
		5 & 0.5 & 0.10 & 2022.13 & 2464.54 & 0.03 & 2019.55 & 0.16 & 2468.20 & -- & -- & 2017.94 \\   
		5 & 1.0 & 0.24 & 1602.35 & 2397.08 & 0.04 & 1583.01 & 1.34 & 2424.11 & -- & -- & 1562.58 \\   
		5 & 1.5 & 0.39 & 1231.02 & 2271.26 & 0.04 & 1186.52 & 3.97 & 2321.85 & -- & -- & 1135.17 \\  
		5 & 2.0 & 0.61 & 911.15 & 2084.73 & 0.04 & 844.02 & 8.49 & 2114.73 & -- & -- & 707.77 \\  
		\midrule
		10 & 0.1 & 0.20 & 4912.69 & 5009.03 & 0.11 & 4912.69 & 0.00 & 5009.03 & 0.10 & 5009.03 & 4912.69 \\   
		10 & 0.5 & 0.25 & 4530.01 & 5001.22 & 0.11 & 4527.34 & 0.06 & 5009.03 & -- & -- & 4527.34 \\   
		10 & 1.0 & 0.53 & 4072.51 & 4969.01 & 0.12 & 4050.19 & 0.57 & 4980.27 & -- & -- & 4045.64 \\   
		10 & 1.5 & 0.91 & 3639.15 & 4901.17 & 0.12 & 3594.08 & 1.28 & 4941.42 & -- & --& 3592.30 \\   
		10 & 2.0 & 1.44 & 3229.58 & 4836.17 & 0.14 & 3155.97 & 2.38 & 4864.44 & -- & -- & 3138.96 \\   
		\midrule
		15 & 0.1 & 0.36 & 7416.52 & 7513.95 & 0.22 & 7416.52 & 0.00 & 7513.95 & 0.17 & 7513.95 & 7416.52 \\   
		15 & 0.5 & 0.47 & 7028.22 & 7508.13 & 0.23 & 7026.80 & 0.02 & 7513.95 & -- & -- & 7026.80 \\   
		15 & 1.0 & 1.03 & 6557.32 & 7483.57 & 0.25 & 6539.66 & 0.28 & 7513.95 & -- & -- & 6539.66 \\   
		15 & 1.5 & 1.70 & 6105.29 & 7445.90 & 0.26 & 6061.74 & 0.74 & 7472.86 & -- & -- & 6071.54 \\   
		15 & 2.0 & 2.56 & 5664.92 & 7365.34 & 0.29 & 5593.44 & 1.31 & 7459.25 & -- & -- & 5603.43 \\ 
		\midrule
		20 & 0.1 & 0.62 & 9840.40 & 9936.95 & 0.37 & 9840.40 & 0.00 & 9936.95 & 0.29 & 9936.95 & 9840.40 \\   
		20 & 0.5 & 0.83 & 9455.48 & 9932.58 & 0.40 & 9454.20 & 0.01 & 9936.95 &-- & -- & 9454.20 \\   
		20 & 1.0 & 1.96 & 8983.11 & 9921.01 & 0.43 & 8971.45 & 0.13 & 9936.95 & -- & -- & 8971.45 \\   
		20 & 1.5 & 3.25 & 8518.85 & 9893.43 & 0.45 & 8489.54 & 0.35 & 9931.17 & -- & -- & 8498.79 \\   
		20 & 2.0 & 4.77 & 8067.97 & 9856.06 & 0.47 & 8008.99 & 0.75 & 9931.17 & -- & -- & 8026.12 \\
		\bottomrule
	\end{tabular}
	\caption{Results for linear instances with convex $\Ucal$. \label{table:result_convexU_linear}}
\end{table}

\begin{table}[ht]
	\centering \small
	\begin{tabular}{llllllllllll}
		\toprule
		$I$&$\theta$ & $t_{\RR}$ & $Z_{\RR}^*$ & $\Ebb[R(\pb^*_{\RR},\ub_0)]$ & $t_{\DR}$ & $Z_{\DR}^*$ & RI(\%) & $R(\pb^*_{\DR},\ub_0)$ & $t_{\Nom}$ & $Z_{\Nom}^*$ & $Z_{\Nom,\WC}$ \\ 
		\midrule
		5 & 0.1 & 0.15 & \numvvm{93744.23} & \numvvm{221222.45} & 0.10 & \numvvm{91972.69} & 4.89 & \numvvm{224749.25} & 0.05 & \numvvm{227228.32} & \numvvm{91511.98} \\   
		5 & 0.5 & 0.31 & \numvvm{10610.01} & \numvvm{113888.98} & 0.16 & \numvvm{6578.80} & 64.55 & \numvvm{77412.76} & -- & -- & \numvvm{3614.58} \\   
		5 & 1.0 & 0.37 & \numvvm{2878.63} & \numvvm{78905.26} & 0.17 & \numvvm{1499.12} & 88.40 & \numvvm{47232.48} & -- & -- & \numvvm{364.41} \\   
		5 & 1.5 & 0.40 & \numvvm{1069.28} & \numvvm{69933.11} & 0.18 & \numvvm{490.14} & 106.08 & \numvvm{22333.65} & -- & -- & \numvvm{84.42} \\   
		5 & 2.0 & 0.41 & \numvvm{424.44} & \numvvm{66460.80} & 0.18 & \numvvm{184.07} & 116.57 & \numvvm{14405.04} & --& -- & \numvvm{29.43} \\ 
		\midrule
		10 & 0.1 & 0.70 & \numvvm{30101269.92} & \numvvm{63310339.83} & 0.44 & \numvvm{28687848.21} & 12.02 & \numvvm{71575092.33} & 0.17 & \numvvm{71616384.73} & \numvvm{28664093.65} \\   
		10 & 0.5 & 1.31 & \numvvm{4374988.35} & \numvvm{39000059.01} & 1.47 & \numvvm{1302758.42} & 202.19 & \numvvm{45130358.67} & -- & -- & \numvvm{947933.08} \\   
		10 & 1.0& 1.91 & \numvvm{933086.04} & \numvvm{32764388.33} & 4.45 & \numvvm{139419.46} & 442.91 & \numvvm{16213017.94} & -- & -- & \numvvm{56786.95} \\   
		10 & 1.5 & 3.05 & \numvvm{315490.74} & \numvvm{19926269.74} & 10.29 & \numvvm{33067.32} & 652.95 & \numvvm{7478100.45} & -- & -- & \numvvm{9555.91} \\   
		10 & 2.0 & 3.94 & \numvvm{135738.48} & \numvvm{18371400.02} & 18.12 & \numvvm{11318.56} & 844.92 & \numvvm{10971494.48} & --& -- & \numvvm{4046.86} \\ 
		\midrule
		15 & 0.1 & 1.80 & \numvvm{2036467052.53} & \numvvm{4987407476.52} & 1.77 & \numvvm{1988832731.72} & 9.24 & \numvvm{5053720805.52} & 0.51 & \numvvm{5053760084.53} & \numvvm{1988827186.56} \\   
		15 & 0.5 & 3.66 & \numvvm{201766127.46} & \numvvm{3844202585.82} & 13.40 & \numvvm{62746668.45} & 295.04 & \numvvm{3516364162.02} & -- & -- & \numvvm{50465973.03} \\   
		15 & 1.0 & 7.39 & \numvvm{30283875.43} & \numvvm{2233011100.70} & 133.51 & \numvvm{4651888.91} & 879.80 & \numvvm{1400004784.66} & -- & -- & \numvvm{1548998.64} \\   
		15 & 1.5 & 14.07 & \numvvm{9139774.94} & \numvvm{1607384007.53} & 595.51 & \numvvm{939718.73} & 1451.79 & \numvvm{2664400108.83} & -- & -- & \numvvm{380439.00} \\   
		15 & 2.0 & 18.21 & \numvvm{3985648.35} & \numvvm{1246260207.71} & 800.14 & \numvvm{280858.38} & 1957.47 & \numvvm{2172193845.69} & -- & -- & \numvvm{151955.12} \\ 
		\midrule
		20 & 0.1 & 6.25 & \numvvm{157307103440.96} & \numvvm{409858123415.04} & 5.47 & \numvvm{156726301951.09} & 13.12 & \numvvm{410911493712.79} & 1.97 & \numvvm{410937627349.81} & \numvvm{156726332506.78} \\
		20 & 0.5 & 18.29 & \numvvm{10536473168.75} & \numvvm{205653309577.32} & 217.95 & \numvvm{3358904126.29} & 577.40 & \numvvm{409024478409.74} & -- & -- & \numvvm{3343464710.58} \\   
		20 & 1.0 & 24.63 & \numvvm{1534890339.83} & \numvvm{192432938377.34} & 972.80 & \numvvm{102516889.68} & 2848.23 & \numvvm{25473860743.89} & -- & --& \numvvm{37798037.01} \\   
		20 & 1.5 & 41.60 & \numvvm{473276643.33} & \numvvm{132421422000.32} & 1126.14 & \numvvm{19035402.83} & 6113.91 & \numvvm{20538109091.69} & -- & -- & \numvvm{13776536.50} \\   
		20 & 2.0 & 46.11 & \numvvm{200652809.23} & \numvvm{100394840422.63} & 1140.78 & \numvvm{6223231.30} & 8874.96 & \numvvm{365898733369.71} & -- & --& \numvvm{5332544.91} \\ 
		\bottomrule
	\end{tabular}
	\caption{Results for semi-log instances with convex $\Ucal$. \label{table:result_convexU_semi_log}}
\end{table}

\begin{table}[ht]
	\centering \small
	\begin{tabular}{llllllllllll}
		\toprule
		$I$&$\theta$ & $t_{\RR}$ & $Z_{\RR}^*$ & $\Ebb[R(\pb^*_{\RR},\ub_0)]$ & $t_{\DR}$ & $Z_{\DR}^*$ & RI(\%) & $R(\pb^*_{\DR},\ub_0)$ & $t_{\Nom}$ & $Z_{\Nom}^*$ & $Z_{\Nom,\WC}$ \\ 
		\midrule
		5 & 0.1 & 0.12 & \numvvm{1772151.82} & \numvvm{3879101.48} & 0.11 & \numvvm{1655591.33} & 8.35 & \numvvm{4021199.83} & 0.07 & \numvvm{4306037.70} & \numvvm{1553745.85} \\   
		5 & 0.5 & 0.18 & \numvvm{395089.22} & \numvvm{2638106.01} & 0.15 & \numvvm{344027.61} & 18.10 & \numvvm{2445531.14} & -- & -- & \numvvm{220735.33} \\   
		5 & 1.0 & 0.19 & \numvvm{108282.06} & \numvvm{2452594.37} & 0.15 & \numvvm{94214.30} & 18.67 & \numvvm{2249083.97} & -- & -- & \numvvm{54212.14} \\   
		5 & 1.5 & 0.18 & \numvvm{32007.04} & \numvvm{2405452.04} & 0.15 & \numvvm{27368.12} & 19.99 & \numvvm{2249083.97} & -- & -- & \numvvm{15683.51} \\   
	    5 & 2.0 & 0.19 & \numvvm{9502.35} & \numvvm{2405522.18} & 0.16 & \numvvm{8085.17} & 20.17 & \numvvm{2153057.90} & -- & --& \numvvm{4625.74} \\
		\midrule
		10 & 0.1 & 0.51 & \numvvm{28675429.28} & \numvvm{60304775.39} & 1.41 & \numvvm{23186403.66} & 27.15 & \numvvm{74350900.63} & 0.37 & \numvvm{76582442.47} & \numvvm{23059278.94} \\   
		10 & 0.5 & 1.15 & \numvvm{6160902.13} & \numvvm{36656659.46} & 15.88 & \numvvm{3135534.90} & 90.19 & \numvvm{26122270.54} & -- & -- & \numvvm{1871543.35} \\   
		10 & 1.0 & 1.81 & \numvvm{2195211.74} & \numvvm{25815944.68} & 46.27 & \numvvm{1019616.82} & 111.15 & \numvvm{17704670.75} & -- & -- & \numvvm{551838.24} \\   
		10 & 1.5 & 2.19 & \numvvm{997459.29} & \numvvm{22968673.10} & 59.82 & \numvvm{466412.79} & 113.12 & \numvvm{13479232.93} & -- & --& \numvvm{234807.75} \\   
		10 & 2.0 & 2.54 & \numvvm{503920.87} & \numvvm{21204878.12} & 63.96 & \numvvm{238256.91} & 113.39 & \numvvm{14226191.86} & -- & -- & \numvvm{114156.02} \\ 
		\midrule
		15 & 0.1 & 1.59 & \numvvm{250597781.96} & \numvvm{746268214.27} & 8.39 & \numvvm{234557915.09} & 19.38 & \numvvm{828146116.26} & 1.71 & \numvvm{838470429.14} & \numvvm{231526331.95} \\   
		15 & 0.5 & 5.71 & \numvvm{43277093.11} & \numvvm{286443608.89} & 728.17 & \numvvm{16669947.97} & 202.98 & \numvvm{245775886.78} & -- & -- & \numvvm{8211144.28} \\   
		15 & 1.0 & 9.62 & \numvvm{16657463.44} & \numvvm{201501570.35} & 1155.56 & \numvvm{5072757.24} & 257.70 & \numvvm{181962649.00} & -- & -- & \numvvm{2292462.48} \\   
		15 & 1.5 & 13.02 & \numvvm{8594507.97} & \numvvm{155864720.76} & 1200.60 & \numvvm{2305581.63} & 289.39 & \numvvm{120020635.42} & -- & -- & \numvvm{1032472.81} \\   
		15 & 2.0 & 16.49 & \numvvm{4990160.88} & \numvvm{135978319.02} & 1200.66 & \numvvm{1249090.44} & 309.56 & \numvvm{51478973.94} & -- & -- & \numvvm{543608.49} \\ 
		\midrule
		20 & 0.1 & 4.82 & \numvvm{2520583780.43} & \numvvm{7863789479.40} & 79.04 & \numvvm{2273994801.77} & 41.87 & \numvvm{8629448674.87} & 5.17 & \numvvm{8633354302.79} & \numvvm{2270476162.92} \\   
		20 & 0.5 & 14.16 & \numvvm{371946898.60} & \numvvm{3037516223.67} & 1200.35 & \numvvm{66754433.94} & 512.22 & \numvvm{1385762909.75} & --& -- & \numvvm{27686996.61} \\   
		20 & 1.0 & 24.35 & \numvvm{141671404.39} & \numvvm{1749282913.12} & 1200.74 & \numvvm{18521018.11} & 745.22 & \numvvm{742282408.96} & -- & -- & \numvvm{5841538.52} \\   
		20 & 1.5 & 37.00 & \numvvm{74367613.77} & \numvvm{1343743086.38} & 1200.87 & \numvvm{8692959.77} & 870.06 & \numvvm{342449021.68} & -- & -- & \numvvm{2753130.91} \\   
		20 & 2.0 & 44.84 & \numvvm{44502090.74} & \numvvm{1128394734.62} & 1200.96 & \numvvm{4741264.36} & 925.14 & \numvvm{323582530.32} & -- & -- & \numvvm{1572476.25} \\
		\bottomrule
	\end{tabular}
	\caption{Results for log-log instances with convex $\Ucal$. \label{table:result_convexU_log_log}}
\end{table}

Besides linear demand, these results also show that for semi-log and log-log demand, there can be a very large difference between the randomized and deterministic robust pricing schemes. The benefit of randomization, quantified by the metric RI, ranges from about $5\%$ to as much as $8875\%$ for semi-log instances, and from about $8\%$ to $925\%$ for log-log instances. Note that the magnitude of RI for the semi-log instances is larger than that for log-log, because the logarithm of demand in the semi-log model has a linear dependence on price which results in an exponential dependence of demand on price, but in log-log, the logarithm of demand is linear in the logarithm of price, resulting in a milder polynomial dependence of demand on price. For semi-log and log-log demand, both the worst-case revenue of RRPO solution and the worst-case revenue of DRPO solution decrease as the uncertainty set becomes larger, and the rate of reduction becomes less as the uncertainty budget $\theta$ is larger. In addition, for linear, semi-log and log-log demand, the RI generally increases as the uncertainty budget $\theta$ increases. Also, as we expect, $Z_{\RR}^* \geq Z_{\DR}^* \geq Z_{\Nom,\WC}$. Interestingly, the randomized robust pricing scheme can sometimes achieve better performance than the deterministic robust scheme under the nominal demand model (for example, compare $\Ebb[R(\pb^*_{\RR},\ub_0)]$ and $R(\pb^*_{\DR},\ub_0)$ for log-log demand with $I = 10$, $\theta \geq 0.5$); this appears to be the case for almost all $(I,\theta)$ combinations for log-log, and for a smaller set of $(I, \theta)$ combinations for semi-log. 

With regard to the computation time, we observe that the computation time generally grows with the number of products for both RRPO and DRPO. For linear demand, both RRPO and DRPO can be solved extremely quickly (no more than 5 seconds on average, even with $I = 20$ products). For log-log and semi-log, when the number of products is held constant, the amount of time required to solve either RRPO generally becomes larger as the uncertainty set becomes larger, but in all cases the average time is under a minute.  What we also find is that for both log-log and semi-log demand, RRPO generally requires much less time to solve to complete optimality than DRPO; this is likely because the nominal problem (which is a key piece of the constraint generation method for RRPO when $\Ucal$ is convex) can be solved rapidly, whereas the robust version of this mixed-integer exponential cone program is more challenging for Mosek.

\subsection{Results using real data instances}

\label{subsec:results_orangejuice}

In our next set of experiments, we evaluate the effectiveness of solution algorithms on problem instances calibrated with real data. For these experiments, we consider the \orangeJuice data set from \cite{montgomery1997creating}, which was accessed via the {\tt bayesm} package in R \citep{rossi2022bayesm}. This data set contains price and sales data for $I=11$ different orange juice brands at the Dominick’s Finer Foods chain of grocery stores in the Chicago area. Each observation in the data set consists of: the store $s$; the week $t$; the log of the number of units sold $\log(q_{t,s,i})$ for brand $i$;  the prices $p_{t,s,1}$...$p_{{t,s,11}}$ of the eleven orange juice brands; the dummy variable $d_{t,s,i}$ indicating whether brand $i$ had any in-store displays at store $s$ in week $t$; and the variable $f_{t,s,i}$ indicating if brand $i$ was featured/advertised at store $s$ in week $t$. 
We fit log-log and semi-log regression models for each brand $i$ according to the following specifications:
\begin{align} 
\text{(semi-log)}\quad \log(q_{t,s,i}) & = \alpha_i - \beta_i p_{t,s,i} + \sum_{j\neq i} \gamma_{ij} p_{t,s,j} + \psi_i d_{t,s,i} + \phi_i f_{t,s,i} + \epsilon_{t,s,i}, \\
\text{(log-log)} \quad \log(q_{t,s,i}) & = \alpha_i - \beta_i \log(p_{t,s,i}) + \sum_{j\neq i} \gamma_{ij} \log(p_{t,s,j}) + \psi_i d_{t,s,i} + \phi_i f_{t,s,i}+ \epsilon_{t,s,i},
\end{align}
where $\{ \epsilon_{t,s,i} \}_{t,s,i}$ is a collection of IID normally distributed error terms. The point estimates of the model parameters are provided in Section~\ref{subsec:extra_results_orangejuice_estimates} of the ecompanion.  We note that prior work has considered the estimation of both of these types of models (see the examples in \citealt{rossi2022bayesm}; see also \citealt{montgomery1997creating} and \citealt{misic2020optimization}). 

We consider the problem of obtaining a price vector $\pb = (p_1, p_2,\dots, p_{11})$ for this collection of 11 products. To formulate the price vector set $\Pcal$, we assume that each product $i$ has five allowable prices, which are shown in Table~\ref{table:price_set}. These prices correspond to the 0th (i.e., minimum), 25th, 50th, 75th and 100th (i.e., maximum) percentiles of the observed prices in the dataset. 
\begin{table}[ht]
	\centering
	\begin{tabular}{rrrrrrrrrrrr}
		\toprule 
		Product & (1)& (2) & (3) & (4) & (5) & (6) & (7) & (8) & (9) & (10) & (11) \\ \midrule
%		& 1.29 & 1.91 & 1.25 & 0.99 & 0.88 & 1.84 & 0.91 & 0.91 & 0.69 & 0.52 & 1.00 \\ 
%		& 2.49 & 2.79 & 2.69 & 1.99 & 1.99 & 2.45 & 1.99 & 1.99 & 1.79 & 1.58 & 1.50 \\ 
%		& 2.99 & 3.17 & 2.89 & 2.35 & 2.17 & 2.64 & 2.39 & 2.19 & 1.99 & 1.59 & 1.79 \\ 
%		& 3.19 & 3.33 & 3.12 & 2.49 & 2.49 & 2.99 & 2.56 & 2.39 & 2.36 & 1.99 & 2.00 \\ 
%		& 3.87 & 3.88 & 3.35 & 3.06 & 3.17 & 3.39 & 3.07 & 2.69 & 3.08 & 2.69 & 2.50 \\ \bottomrule
		& 1.29 & 2.86& 1.25 & 0.99 & 0.88 & 2.76 & 0.91 & 0.91 & 0.69 & 0.52 & 1.99 \\ 
		& 2.49 & 4.19 & 2.69 & 1.99 & 1.99 & 3.67 & 1.99 & 1.99 & 1.79 & 1.58 & 2.99\\ 
		& 2.99 & 4.75 & 2.89 & 2.35 & 2.17 & 3.96& 2.39 & 2.19 & 1.99 & 1.59 & 3.59 \\ 
		& 3.19 & 4.99 & 3.12 & 2.49 & 2.49 & 4.49 & 2.56 & 2.39 & 2.36 & 1.99 & 3.99 \\ 
		& 3.87 & 5.82 & 3.35 & 3.06 & 3.17 & 5.09 & 3.07 & 2.69 & 3.08 & 2.69 & 4.99 \\ \bottomrule
	\end{tabular}
	\caption{Possible price levels for products in \orangeJuice experiment instances. \label{table:price_set}}
\end{table}

For each type of demand model, we consider two forms of uncertainty set: a convex $L1$-norm uncertainty set (as in equation~\eqref{eq:L1_norm_Ucal}) and a discrete budget uncertainty set (as in equation~\eqref{eq:discrete_budget_Ucal}). Specifically, for the discrete budget uncertainty set, we assume that $\bar{\alphab} = 1.2\alphab$, $\underline{\alphab} = 0.8\alphab$, $\bar{\betab} = 1.3\betab$,  $\underline{\betab} = 0.7\betab$, $\bar{\gammab} = 1.4\gammab$, and $\underline{\gammab} = 0.6\gammab$. Tables~\ref{table:orangejuice_convexU_semi_log} and \ref{table:orangejuice_convexU_log_log} below present the results under the convex $L1$-norm uncertainty set for the semi-log and log-log demand models, respectively. Due to page considerations, the results for the discrete $\Ucal$ case are provided in Section~\ref{subsec:extra_results_orangejuice_discreteUcal} of the ecompanion.

We can see from Tables \ref{table:orangejuice_convexU_semi_log} and \ref{table:orangejuice_convexU_log_log} that the randomized pricing strategy performs significantly better than the deterministic pricing solution under the worst-case demand model, with the RI ranging from $17.86\%$ to $47.81\%$ for semi-log demand and from $27.71\%$ to $92.31\%$ for log-log demand. In addition, for the same demand type and uncertainty set, the computation time of RRPO is comparable to that of DRPO. With regard to the discrete uncertainty set case, the results shown in Section~\ref{subsec:extra_results_orangejuice_discreteUcal} are qualitatively similar, with the randomized robust pricing strategy similarly outperforming the deterministic robust solution. We do also observe that under both demand models, the time to solve the RRPO problem is generally smaller than for the DRPO problem (RRPO requires no more than about 17 seconds, whereas DRPO can take up to 300 seconds on average in the worst case). 

%solving RRPO with the discrete uncertainty set requires more time than solving it with convex uncertainty set, although the overall time is still reasonable (in the most extreme case, RRPO for the discrete uncertainty set can take up to approximately 300 seconds, and DRPO requires up to 600 second, compared to 60 seconds for both RRPO and DRPO for the $L1$-norm uncertainty set). 

%We can see from Tables \ref{table:orangejuice_convexU_semi_log}, \ref{table:orangejuice_convexU_log_log}, \ref{table:orangejuice_discreteU_semi_log} and \ref{table:orangejuice_discreteU_log_log} that randomized pricing performs much better than deterministic pricing under the worst cases, with the RI ranging from $17.86\%$ to $60.81\%$ for semi-log demand and from $27.71\%$ to $92.31\%$ for log-log demand. In addition, for the same demand type and uncertainty set, the computation time of RRPO is comparable to that of DRPO. We also see that, for both semi-log and log-log demand models, solving RRPO with discrete uncertainty set requires much more time than solving it with convex uncertainty set. 

\begin{table}[ht]
	\centering
	\small
		\begin{tabular}{lllllllllll}
			\toprule
			$\theta$ & $t_{\RR}$ & $Z_{\RR}^*$ & $\Ebb[R(\pb^*_{\RR},\ub_0)]$ & $t_{\DR}$ & $Z_{\DR}^*$ & RI(\%) & $R(\pb^*_{\DR},\ub_0)$ & $t_{\Nom}$ & $Z_{\Nom}^*$ & $Z_{\Nom,\WC}$ \\ 
			\midrule
			0.10 & 16.86 & 342357.06 & 481125.25 & 15.98 & 290474.67 & 17.86 & 590546.51 & 0.83 & 590547.01 & 290474.76 \\   0.50 & 41.42 & 197517.06 & 373973.40 & 47.65 & 147748.35 & 33.68 & 294483.28 & -- & -- & 96016.90 \\   
			0.80 & 42.61 & 149709.04 & 352742.22 & 67.83 & 105734.14 & 41.59 & 294483.28 & -- & -- & 67924.78 \\   
			1.00 & 67.47 & 125987.02 & 349644.33 & 74.03 & 86977.24 & 44.85 & 265173.68 & -- & --& 55394.70 \\   
			1.50 & 61.70 & 82880.96 & 348467.74 & 44.65 & 56474.64 & 46.76 & 265173.68 & -- & -- & 34864.43 \\   
			2.00 & 65.21 & 54665.15 & 348466.26 & 52.10 & 37164.75 & 47.09 & 265173.68 & -- & -- & 22615.70 \\ 
			\bottomrule
		\end{tabular}
	\caption{Results for \orangeJuice pricing problem with semi-log demand and convex $\Ucal$. \label{table:orangejuice_convexU_semi_log}}
\end{table}

\begin{table}[ht]
	\centering
	\small
		\begin{tabular}{lllllllllll}
			\toprule
			$\theta$ & $t_{\RR}$ & $Z_{\RR}^*$ & $\Ebb[R(\pb^*_{\RR},\ub_0)]$ & $t_{\DR}$ & $Z_{\DR}^*$ & RI(\%) & $R(\pb^*_{\DR},\ub_0)$ & $t_{\Nom}$ & $Z_{\Nom}^*$ & $Z_{\Nom,\WC}$ \\ 
			\midrule
			0.10 & 0.71 & 722647.22 & 1051269.80 & 1.97 & 565866.71 & 27.71 & 922172.97 & 0.88 & 1112050.59 & 560812.30 \\   0.50 & 1.42 & 342614.34 & 687632.84 & 4.99 & 233387.10 & 46.80 & 782893.68 & -- & -- & 152881.89 \\   
			0.80 & 1.91 & 260049.66 & 672481.74 & 6.77 & 162276.97 & 60.25 & 782893.68 & -- & --& 102893.20 \\   
			1.00 & 2.08 & 217580.86 & 683576.99 & 9.41 & 128220.45 & 69.69 & 782893.68 & -- & -- & 81427.57 \\   
			1.50 & 2.26 & 142307.66 & 670759.12 & 13.02 & 75897.66 & 87.50 & 599315.12 & -- & -- & 48983.56 \\   
			2.00 & 2.39 & 94847.37 & 670758.38 & 13.49 & 49319.21 & 92.31 & 377932.71 & -- & -- & 31055.19 \\ 
			\bottomrule
		\end{tabular}
	\caption{Results for \orangeJuice pricing problem with log-log demand and convex $\Ucal$. \label{table:orangejuice_convexU_log_log}}
\end{table}

Lastly, it is also interesting to compare the randomized robust pricing strategy to the deterministic robust price vector. Taking the log-log demand model and the convex $L1$ uncertainty set with $\theta = 0.8$ as an example, the solution of the RRPO problem is the following randomized pricing strategy:
\begin{equation}
	\pb = \left\{ \begin{array}{llll} 
		(3.87, 5.82, 1.25, 0.99, 3.17, 5.09, 3.07, 0.91, 0.69, 2.69, 1.99) &&& \text{w.p.}\ 0.1628, \\
		(1.29, 5.82, 3.35, 3.06, 0.88, 2.76, 3.07, 2.69, 3.08, 2.69, 4.99)&&& \text{w.p.}\ 0.1752, \\
		(3.87, 2.86, 1.25, 3.06, 3.17, 5.09, 0.91, 2.69, 3.08, 0.52, 4.99) &&& \text{w.p.}\ 0.2658,  \\
		(3.87, 5.82, 3.35, 3.06, 0.88, 2.76, 3.07, 0.91, 3.08, 2.69, 4.99) &&& \text{w.p.}\ 0.0381, \\
		(3.87, 2.86, 1.25, 3.06, 3.17, 2.76, 3.07, 2.69, 0.69, 2.69, 4.99)&&& \text{w.p.}\ 0.3258, \\
		(3.87, 2.86, 1.25, 3.06, 3.17, 5.09, 0.91, 0.91, 3.08, 0.52, 4.99) &&& \text{w.p.}\ 0.0323. \\
	\end{array} \right.
\end{equation}
% semi-log, theta = 0.8 below:
%\begin{equation}
%	\pb = \left\{ \begin{array}{llll} 
%		(3.87,  3.88,  1.25,  0.99,  3.17,  3.39,  3.07,  0.91 , 0.69,  2.69 , 1.00) &&& \text{w.p.}\ 0.1431, \\
%		(3.87 , 1.91 , 3.35  ,3.06  ,3.17,  3.39  ,0.91  ,0.91  ,3.08  ,0.52 , 2.50)&&& \text{w.p.}\ 0.1366, \\
%		(3.87 , 3.88,  1.25,  3.06,  3.17 , 1.84,  3.07 , 2.69  ,0.69 , 2.69 , 2.50) &&& \text{w.p.}\ 0.1741, \\
%		(3.87 , 3.88  ,3.35  ,3.06,  0.88,  1.84 , 3.07,  0.91 , 3.08  ,2.69,  2.50) &&& \text{w.p.}\ 0.0467, \\
%		(3.87 , 1.91  ,1.25 , 3.06 , 3.17,  1.84 , 3.07,  2.69 , 0.69 , 2.69 , 2.50)&&& \text{w.p.}\ 0.1546, \\
%		(3.87 , 1.91 , 1.25 , 3.06  ,3.17  ,1.84,  0.91 , 0.91 , 3.08 , 0.52  ,2.50) &&& \text{w.p.}\ 0.0403, \\
%		(1.29 , 3.88 , 3.35 , 3.06,  0.88 , 1.84,  3.07 , 0.91 , 3.08,  2.69 , 2.50) &&& \text{w.p.}\ 0.1817, \\
%		(3.87 , 1.91  ,1.25 , 3.06 , 3.17  ,3.39  ,0.91 , 0.91 , 3.08  ,0.52 , 2.50) &&& \text{w.p.}\ 0.1229.\\
%	\end{array} \right.
%\end{equation}
Observe that in this randomized pricing strategy, each price vector is such that the product is set to either its lowest or highest allowable price. This is congruent with Proposition~\ref{proposition:loglog_extremal}, which suggests that the nominal problem under the log-log demand model will always have a solution that involves setting each product to its highest or lowest price; since our solution algorithm is based on constraint generation using this nominal problem as a separation procedure, it makes sense that the randomized price vector will be supported on such extremal price vectors. On the other hand, the solution of the DRPO problem is the price vector $\pb_{\DR} = (3.87, 2.86, 1.25, 3.06, 3.17, 2.76, 0.91, 2.69, 0.69, 0.52, 4.99)$, for which we observe that the chosen prices are also either the lowest or highest for each product.

% ORANGE JUICE DDE 

\subsection{Data-driven evaluation of out-of-sample performance}
\label{subsec:R1_results_orangejuice}

In this section, we present the results of a data-driven experiment to demonstrate the benefit of randomized robust decisions in terms of out-of-sample revenue. The previous numerical experiments show that from a risk standpoint, randomized robust pricing decisions can lead to better worst-case performance than deterministic robust pricing decisions and nominal/non-robust deterministic pricing decisions. In contrast, this section focuses on a different question, which we state as follows. Suppose we are given a data set of historical price vectors and demands, and we use this data set to craft an appropriate uncertainty set for the randomized robust price optimization problem; how does the resulting randomized robust pricing decision perform with respect to the true, unknown demand model? 

To understand this question, we will again consider the \orangeJuice data set used previously in Section~\ref{subsec:results_orangejuice}, and compare the randomized robust, deterministic robust and nominal approaches. The high-level idea of this experiment is as follows: for each store, divide the available price-demand data into training, validation and test sets. We use the training set to estimate a nominal demand model. We use the validation set to tune the uncertainty set for the randomized robust and deterministic robust approaches. We then use the test set to evaluate the randomized robust, deterministic robust and nominal approaches. 

A key challenge that we face in using the \orangeJuice data set (or any real data set, for that matter) is that the historical data only gives demand (and hence revenue) information on a small, finite number of price vectors. For any price vector that was not offered in a given store, we do not know what the revenue would have been. Thus, to perform the validation and test steps in our experiment, we proceed by considering a \emph{restricted pricing problem}. Namely, in the validation step, we solve restricted versions of the RRPO and DRPO problems where we can only select from price vectors that appear in the validation set. Then, given the pricing strategy that is returned by RRPO/DRPO, we evaluate its true revenue using the \emph{real} revenues of the corresponding price vectors in the validation set. Similarly, in the test step, we solve restricted versions of the RRPO, DRPO and NPO problems where we can only use price vectors that appear in the test set, where RRPO and DRPO are solved using the uncertainty sets designed in the validation step. Given the pricing strategies produced by the restricted RRPO, DRPO and NPO problems, we calculate their revenues using the true revenues for their price vectors in the test set. 

To formalize this experiment, we let $\Dcal = \{ (\pb_i, \db_i) \}_{i=1}^{n'}$ be the historical data set for a single store in the data set, where $\pb_i$ is the vector of prices in week $i$ and $\db_i$ is the vector of demands in week $i$; consequently, $\pb_i^T \db_i$ is the revenue garnered in week $i$. Let $\Pcal = \{ \pb_1,\dots, \pb_{n} \}$ be the set of unique price vectors corresponding to $\Dcal$. 

Let $\Pcal_{\train}$, $\Pcal_{\val}$, $\Pcal_{\train}$ be a partition of the set of indices $[n] = \{1,\dots,n\}$, such that $\Pcal = \Pcal_{\train} \cup \Pcal_{\val} \cup \Pcal_{\test}$, where each $\Pcal_{\cdot}$ set indicates which price vectors are assigned to the training, validation and test sets. Let the functions $\bar{R}: \Pcal \to \Rbb$ denote the average revenue of each price vector in the data set, which is defined as
\begin{equation}
\bar{R}(\pb) = \frac{1}{ | \{i \in [n'] \mid \pb_i = \pb \} | } \cdot \sum_{ \substack{i = 1:\\ \pb_i = \pb}}^{n'} \pb^T \db_i, \\
\end{equation}
% \begin{align}
% \bar{R}_{\val}(\pb) & = \frac{1}{ | \{i \in [n'] \mid \pb_i = \pb \} | } \cdot \sum_{i = 1: \pb_i = \pb}^{n'} \pb^T \db_i, \\
% \bar{R}_{\test}(\pb) = \frac{1}{ |  \{i \in [n'] \mid \pb_i = \pb \} | } \cdot \sum_{i = 1: \pb_i = \pb}^{n'} \pb^T \db_i.
% \end{align}
Note that in cases where $\pb$ only occurs in a single week in the validation or test sets, then $\bar{R}(\pb)$ will be the revenue of $\pb$ in that single historical week. We also point out here that $\bar{R}(\pb)$ is distinct from $R(\pb,\ub)$: the former is our data-driven estimate of the average revenue of the price vector $\pb$, whereas $R(\pb, \ub)$ is the predicted revenue from a demand model with parameters $\ub$. We additionally let $\Dcal_{\train} = \{ (\pb, \db) \in \Dcal \mid \pb \in \Pcal_{\train} \}$ denote the price-demand observations for the price vectors corresponding to the training set.

Let us assume that the uncertainty set for the DRPO and RRPO approaches is parametrized by a vector $\wb$, so that $\Ucal(\wb)$ corresponds to the uncertainty set with the control parameters $\wb$. We let $\Wcal$ denote the set of control parameters that we will choose from in our validation step. 

%Suppose that each price vector in $\Dcal$ is distinct (i.e., no price vector is repeated), and given a price vector $\pb_i \in \Pcal_{\train} \cup \Pcal_{\val} \cup \Pcal_{\test}$, let $\bar{R}(\pb_i) = \pb_i^T \db_i$ be the historical revenue in the data set. Note that $\bar{R}(\pb)$ is the historical revenue for price vector $\pb$, whereas we will use $R(\pb, \ub)$ to denote the revenue prediction of the underlying demand model structure for a price vector $\pb$ and demand model parameter vector $\ub$. 

We now describe our overall experimental procedure as Algorithm~\ref{alg:DDE}.  
\begin{algorithm}
\begin{algorithmic}[1]
\STATE Estimate nominal demand model parameter vector $\ub^{\Nom}$ using $\Dcal_{\train}$.
\FOR{$\wb \in \Wcal$}
	\STATE Solve $\pb^{\DR, \wb} \in \arg \max_{\pb \in \Pcal_{\val}} \min_{\ub \in \Ucal(\wb)} R(\pb, \ub)$. 
	\STATE Solve $\pib^{\RR, \wb} \in \arg \max_{\pib \in \Delta_{\Pcal_{\val}}} \min_{\ub \in \Ucal(\wb)} \sum_{\pb \in \Pcal_{\val}} \pi_{\pb} \cdot R(\pb, \ub)$. 
\ENDFOR
\STATE Set $\wb_{\DR} = \arg \max_{\wb \in \Wcal} \bar{R}(\pb^{\DR, \wb})$.
\STATE Set $\wb_{\RR} = \arg \max_{\wb \in \Wcal} \sum_{\pb \in \Pcal_{\val}} \pi^{\RR,\wb}_{\pb} \bar{R}(\pb)$. \\[0.5em]
\STATE Solve $\pb^{\Nom,\test} \in \arg \max_{\pb \in \Pcal_{\test}} R(\pb, \ub^{\Nom})$.
\STATE Calculate $R^{\Nom, \test} = \bar{R}(\pb^{\Nom, \test})$. \\[0.5em]
\STATE Solve $\pb^{\DR,\test} \in \arg \max_{\pb \in \Pcal_{\test}} \min_{\ub \in \Ucal(\wb_{\DR})} R(\pb, \ub)$. 
\STATE Calculate $R^{\DR, \test} = \bar{R}(\pb^{\DR,\test})$. \\[0.5em]
\STATE Solve $\pi^{\RR, \test} \in \arg \max_{\pib \in \Delta_{\Pcal_{\test}}} \min_{\ub \in \Ucal(\wb_{\RR})} \sum_{\pb \in \Pcal_{\test}} \pi_{\pb} R(\pb, \ub)$.
\STATE Calculate $R^{\RR, \test} = \sum_{\pb \in \Pcal_{\test}} \pi^{\RR, \test}_{\pb} \bar{R}(\pb)$.
\RETURN $R^{\Nom, \test}$, $R^{\DR, \test}$, $R^{\RR, \test}$. 
\end{algorithmic}
\caption{Data-driven experiment procedure for \orangeJuice data set. \label{alg:DDE} }
\end{algorithm}

We apply Algorithm~\ref{alg:DDE} to the price-demand data of each store in the orange juice data set. For each store, we partition the set of $n$ historical price vectors so that approximately the first 50\% of price vectors constitute the training set $\Pcal_{\train}$, the next 25\% of price vectors constitute the validation set $\Pcal_{\val}$, and the final 25\% of price vectors constitute the test set $\Pcal_{\test}$. Note that the price vectors in $\Pcal$ are sorted in order of temporal appearance; thus, the first occurrence in time of each price vector in $\Pcal_{\train}$ is before that of any price vector in $\Pcal_{\val}$, and the first occurrence in time of each price vector in $\Pcal_{\val}$ is before that of any price vector in $\Pcal_{\test}$. 

For the uncertainty set $\Ucal$, we consider the following budget uncertainty set:
\begin{equation}
\Ucal(\wb) = \left\{ \ub = (\ab, \betab, \gammab) \ \vline \ \| \tilde{\ub} \|_1 \leq w_1;\ \| \tilde{\ub} \|_{\infty} \leq w_2;\ \tilde{u}_k = \frac{u_k - u^{\Nom}_k}{ u^{\Nom}_k }, \forall k \in \{1,\dots,I+I^2\}  \right\},
\end{equation}
where the term $\tilde{u}_k$ measures the relative deviation of the $k$th demand model parameter from its nominal value, the parameter $w_1$ is an upper bound on the aggregate deviation of all of the demand model parameters (i.e., the budget of uncertainty), and the parameter $w_2$ is an upper limit on the allowable relative deviation of each demand model parameter from its nominal value. We set $\Wcal = \{0,0.2,0.4,\dots,10.0\} \times \{0,0.1,\dots,1.0\}$. We test both the semi-log and log-log demand models.  

In addition to the deterministic robust and randomized robust approaches, we also consider an approach that combines the two, which we refer to as the hybrid approach. The hybrid approach works as follows. In the validation phase, we calculate $R^{\RR, \val}$ and $R^{\DR, \val}$ as the best values of the randomized robust and deterministic robust approaches, respectively, on the validation set:
\begin{align}
R^{\RR, \val} & = \max_{\wb \in \Wcal} \bar{R}(\pb^{\RR, \wb}), \\
R^{\DR, \val} & = \max_{\wb \in \Wcal} \bar{R}(\pb^{\DR, \wb}).
\end{align}
Then in the test phase, we implement the randomized robust approach if $R^{\RR, \val} \geq R^{\DR, \val}$ and the deterministic robust approach otherwise (i.e., if $R^{\RR, \val} < R^{\DR, \val}$). The performance of this hybrid method on the test set can thus be written as 
\begin{equation}
R^{\Hybrid, \test} = R^{\RR, \test} \cdot \Ibb\{ R^{\RR, \val} \geq R^{\DR, \val} \} + R^{\DR, \test} \cdot \Ibb\{R^{\RR, \val} < R^{\DR, \val} \}.
\end{equation}
The rationale of this approach is that when the randomized robust pricing outperforms deterministic robust pricing on the validation set, it is reasonable to expect that randomized robust pricing would outperform deterministic robust pricing on the test set; and conversely, when the deterministic approach outperforms the randomized approach on the validation set, one would expect this to also be true for the test set. We denote the hybrid method by $\Hybrid$.

%In addition to Algorithm~\ref{alg:DDE}, we also consider the following variant of Algorithm~\ref{alg:DDE}. In lines XXX and YYY of Algorithm~\ref{alg:DDE}, the tuned control parameters $\wb_{\DR}$ and $\wb_{\RR}$ are those which give the best performance of the DR and RR strategies on the validation set. Instead of picking the values of $\wb$ that give the best performance, we can instead consider the $\wb$ with the smallest value of $w_1$ that is within a certain fraction $\rho$ of the optimal validation set performance. In particular, we can set 

To compare the performance of all of the approaches, we compute a number of result metrics. For each store in the data set, the metric $\RI_{\Nom \to m}$, where $m \in \{\DR, \RR, \Hybrid\}$,  is defined as
\begin{align}
\RI_{\Nom \to m} & = \frac{ R^{m, \test} - R^{\Nom, \test}}{R^{\Nom, \test}} \times 100\%, 
%\RI_{\Nom \to \RR} & = \frac{ R^{\RR, \test} - R^{\Nom, \test}}{R^{\Nom, \test}} \times 100\%,
\end{align}
which represents the relative improvement of approach $m$ over the nominal approach.  %We calculate the average value of these metrics over all of the stores, which we denote as $\RI^{\avg}_{\Nom \to \DR}$ and $\RI^{\avg}_{\Nom \to \RR}$. 

For each store, we also calculate the gap metric $G_m$ for a given method $m \in \{\Nom, \DR, \RR, \Hybrid \}$ as 
\begin{align}
G_m = \frac{ R^{\max, \test} - R^{m, \test}}{ R^{\max, \test}} \times 100\%,
\end{align}
where $R^{\max, \test} = \max_{m' \in \{\Nom, \DR, \RR, \Hybrid \}} R^{m', \test}$ is the best test set performance of the four methods. 

Table~\ref{table:R1_DDE} below displays the mean values of the $\RI_{\Nom \to m}$ and $G_m$ metrics. In addition to the ordinary mean, we also report the trimmed mean with $\alpha = 0.05$ and $\alpha = 0.10$, where $\alpha$ denotes the fraction of observations to remove from the lower and upper tails (i.e., $\alpha = 0.05$ would mean removing the bottom 2.5\% of values and the top 2.5\% values, and taking the mean of the remaining middle 95\% of values). The reason for considering the trimmed mean is that for some stores, both the DR and RR price prescriptions can exhibit large under or overperformance relative to the nominal pricing strategy. By removing the extreme values at the tails, we obtain a more representative measure of the performance. 

\begin{table}
\centering
\begin{tabular}{llccccccc} \toprule
Demand model & Type of mean & $\RI_{\Nom \to \RR}$ & $\RI_{\Nom \to \DR}$ & $\RI_{\Nom \to \Hybrid}$ & $G_{\Nom}$ & $G_{\DR}$ & $G_{\RR}$ & $G_{\Hybrid}$ \\ 
& & (\%) & (\%) & (\%) & (\%) & (\%) & (\%) & (\%) \\ \midrule
  % Semi-log & Mean                & 4.65 & 1.51 & 2.53 & 5.67 & 5.78 & 3.08 & 4.91 \\ 
  % Semi-log & Trimmed Mean (0.05) & 1.76 & 0.17 & 0.89 & 3.50 & 3.79 & 1.31 & 2.68 \\ 
  % Semi-log & Trimmed Mean (0.10) & 1.40 & 0.21 & 0.81 & 2.43 & 2.48 & 0.85 & 1.47 \\[0.25em]
  % Log-log & Mean                 & 11.88 & 14.31 & 16.34 & 13.29 & 8.06 & 8.25 & 6.78 \\ 
  % Log-log & Trimmed Mean (0.05)  & 8.52 & 8.53 & 10.78 & 11.36 & 6.43 & 6.18 & 5.09 \\ 
  % Log-log & Trimmed Mean (0.10)  & 5.37 & 3.80 & 6.32 & 9.44 & 5.18 & 4.56 & 4.06 \\ \bottomrule
Semi-log  & Mean                  & 4.65 & 1.51 & 3.09 & 5.67 & 5.78 & 3.08 & 4.35 \\ 
  Semi-log  & Trimmed Mean (0.05) & 1.76 & 0.17 & 1.20 & 3.50 & 3.79 & 1.31 & 2.06 \\ 
  Semi-log  & Trimmed Mean (0.10) & 1.40 & 0.21 & 1.08 & 2.43 & 2.48 & 0.85 & 1.12 \\[0.25em]
  Log-log & Mean                  & 11.88 & 14.31 & 16.61 & 13.29 & 8.06 & 8.25 & 6.52 \\ 
  Log-log & Trimmed Mean (0.05)   & 8.52 & 8.53 & 11.07 & 11.36 & 6.43 & 6.18 & 4.80 \\ 
  Log-log & Trimmed Mean (0.10)   & 5.37 & 3.80 & 6.57 & 9.44 & 5.18 & 4.56 & 3.73 \\ \bottomrule
\end{tabular}
\caption{Relative improvement and gap metrics for the data-driven experiment in Section~\ref{subsec:R1_results_orangejuice}. \label{table:R1_DDE}}
\end{table}

From this table, we can see that for the semi-log demand model, the randomized robust approach leads to an improvement over the nominal approach, as the ordinary mean and the two trimmed means for $\RI_{\Nom \to \RR}$ are all greater than zero. Considering the ordinary mean, the randomized robust approach leads to an average improvement over nominal of 4.65\%; after removing the most extreme 5\% and 10\% of values, this average improvement is smaller but still large (1.76\% for the trimmed mean with $\alpha = 0.05$ and 1.40\% for the trimmed mean with $\alpha = 0.10$). Additionally, the randomized robust approach also leads to a greater improvement over nominal compared to the deterministic robust approach, as the mean and two trimmed means of $\RI_{\Nom \to \RR}$ are always greater than the corresponding means for $\RI_{\Nom \to \DR}$. In addition, the randomized robust approach is on average within about 3\% of the performance of the best method on the test set, compared to about 6\% for the nominal and deterministic robust approaches.

For the log-log demand model, the randomized robust approach still leads to an improvement over the nominal approach, as the ordinary mean of $\RI_{\Nom \to \RR}$ is 11.88\%. The deterministic robust approach, on the other hand, leads to a higher improvement over the nominal approach, with an ordinary mean of $\RI_{\Nom \to \DR}$ of 14.31\%. However, comparing the randomized robust and deterministic robust approaches using the two trimmed means of the $\RI$ metrics, we see that the trimmed mean with $\alpha = 0.05$ suggests that the two methods are similar (relative improvement of 8.52\% and 8.53\% for randomized and deterministic), whereas for $\alpha = 0.10$ the randomized robust approach is better (5.37\% for randomized vs 3.80\% for deterministic). These results suggest that while the deterministic robust approach can lead to very large gains over the nominal approach for some stores, leading to a large mean of the $\RI$ metric, the randomized robust approach does a better job of lifting the out-of-sample revenue for most stores and more uniformly improves on the nominal approach. In addition, when we consider the means of the gap metrics, we find that the deterministic and randomized robust approaches exhibit similar gap metrics. For example, using the ordinary mean, the deterministic robust approach is within 8.06\% of the best prescription, whereas the randomized robust approach is within 8.25\%, and both of these are better than the nominal approach (gap of 13.29\%). Using the trimmed mean with $\alpha = 0.10$, the randomized robust approach is slightly better than the deterministic robust approach and the nominal approach (5.18\% for deterministic robust and 9.44\% for nominal versus 4.56\% for randomized robust). 

It is also interesting to consider the hybrid method described earlier. With regard to the semi-log demand model, we find that the hybrid method improves on the deterministic robust approach, but is not as good as the randomized robust approach. For example, when considering the trimmed mean with $\alpha = 0.05$, the hybrid method has an average improvement of 1.20\% over nominal, which is higher than the deterministic robust approach (0.17\%) but smaller than the randomized robust approach (1.76\%). However, for the log-log demand model, the hybrid method actually outperforms both the deterministic and randomized robust approaches. For example, when considering the trimmed mean with $\alpha = 0.05$, the relative improvement of the hybrid approach is 11.07\% compared to 8.52\% and 8.53\% for the randomized robust and deterministic robust approaches respectively. The improvement over the pure robust approaches is also reflected in the gap metric, which is smaller for the hybrid method compared to the nominal, deterministic robust and randomized robust approaches, for all three types of means. These results suggest that from the perspective of out-of-sample performance, deterministic and randomized robust pricing can be synergistic and complementary. This contrasts with our theoretical results, which suggest that from the standpoint of worst-case performance, a randomized robust pricing strategy is always at least as good as a deterministic pricing strategy. 

We conclude this section by commenting on the limitations of this type of experiment. First, the revenue estimates that one obtains using the function $\bar{R}(\cdot)$ are inherently noisy, because each price vector only appears a small number of times for each store. For each store, the number of weeks of data available ranges from 87 to 121 with a mean of 116.3, while the number of unique price vectors ranges from 86 to 121 with a mean of 115.6, which implies that most price vectors only appear once, and thus for most stores and most price vectors, the function $\bar{R}(\cdot)$ reflects the revenue of a single week in which a price vector was offered. Second, as noted earlier, the test phase of Algorithm~\ref{alg:DDE} compares the randomized robust, deterministic robust and nominal approaches by considering a restricted pricing problem where one is restricted to those price vectors in $\Pcal_{\test}$. It is possible that by considering a different set of price vectors, one might see a different ordering of the methods. 

Notwithstanding these limitations, we believe that a strength of the approach is that it is entirely data-driven. We believe that this approach could be generalized to other contexts where the data corresponds to different historical decisions and their (noisy) objective values, and where the ground truth is unknown and inaccessible to us. The application of this data-driven evaluation approach in other settings is an interesting direction for future research.

\section{Conclusions}
\label{sec:conclusions}

In this paper, we considered the problem of designing randomized robust pricing strategies to maximize worst-case revenue. We presented idealized conditions under which such randomized pricing strategies fare no better than the deterministic robust pricing approach, and subsequently we developed solution methods for obtaining the randomized pricing strategies in different settings (when the price set is finite, and when the uncertainty set is either convex or discrete). We showed using both synthetic instances and real data instances that such randomized pricing strategies can lead to large improvements in worst-case revenue over deterministic robust price prescriptions, and using a data-driven experiment based on real data, we also demonstrate that randomized robust pricing can lead to large improvements in out-of-sample revenue over deterministic robust and non-robust pricing.

With regard to future research, there are several interesting directions to consider. First, a potential limitation of randomization in practice is that customers may engage in strategic behavior. For example, in a brick-and-mortar setting where prices are randomized by store or over time, strategic customers may travel to stores with lower prices, or they may wait for better prices before making their purchase. Although this type of customer behavior may not arise in every context (it may be costly for a customer to wait or to travel long distances to obtain a better price on a product), an interesting question is whether it is possible to design a randomized pricing strategy that accounts for this type of strategic behavior. 

Second, it would also be interesting to consider a version of the robust price optimization problem that incorporates contextual information. For example, in the ecommerce setting, different customers who log onto a retailer's website will differ in characteristics (age, web browser, operating system, etc.). This information could be used to craft a richer uncertainty set, and to motivate randomization strategies that randomize differently based on user characteristics. More generally, we hope that this work, which was inspired by the paper of \cite{wang2020randomized}, motivates further study in how randomization can be used to mitigate risk in revenue management applications.

\section*{Acknowledgments}
The authors thank the department editor Chung Piaw Teo, the associate editor and two anonymous reviewers for their comments and suggestions that have helped to significantly strengthen the paper.

\bibliographystyle{plainnat}
\bibliography{RRPO_literature_v3}

\ECSwitch

\section{Omitted proofs}

\subsection{Proof of Theorem~\ref{theorem:concave_revenue_function_set}}
\label{proof_theorem:concave_revenue_function_set}

To prove this, we prove that $Z^*_{\DR} \geq Z^*_{\RR}$. For any $R \in \Rcal$, and any distribution $F$ supported on $\Pcal$, we have
\begin{equation}
\int_{\Pcal} R(\pb)\; dF(\pb) \leq R \left( \int_{\Pcal} \pb \; dF(\pb) \right),
\end{equation}
which follows by Jensen's inequality and the concavity of $R$. This implies that for any $F \in \Fcal$, 
\begin{equation}
\inf_{R \in \Rcal} \int_{\Pcal} R(\pb)\; dF(\pb) \leq \inf_{R \in \Rcal} R \left( \int_{\Pcal} \pb \; dF(\pb) \right).
\end{equation}
Therefore, we have that 
\begin{align*}
Z^*_{\RR} & = \max_{F \in \Fcal} \left\{ \inf_{R \in \Rcal} \int_{\Pcal} R(\pb)\; dF(\pb) \right\} \\
& \leq \max_{F \in \Fcal} \left\{ \inf_{R \in \Rcal} R \left( \int_{\Pcal} \pb \; dF(\pb) \right) \right\} \\
& \leq \max_{\pb \in \Pcal} \left\{ \inf_{R \in \Rcal} R \left( \pb ) \right) \right\} \\
& = Z^*_{\DR},
\end{align*}
where the second inequality follows because $\Pcal$ is assumed to be convex, and thus for any $F \in \Fcal$, $\int_{\Pcal} \pb \; dF(\pb)$ is contained in $\Pcal$. \hfill \Halmos 

\subsection{Proof of Theorem~\ref{theorem:quasiconcave_quasiconvex}}
\label{proof_theorem:quasiconcave_quasiconvex}

We will establish this result by drawing a connection with game theory. Consider a two-player zero-sum game, where player 1 and 2 have strategy spaces $\Pcal$ and $\Ucal$, respectively, and payoff functions $g_1$ and $g_2$ given by 
\begin{align*}
g_1(\pb, \ub) & = R(\pb, \ub), \\
g_2(\pb, \ub) & = - R(\pb, \ub).
\end{align*}
Player 1 can choose the action $\pb \in \Pcal$ (a pure strategy) or some distribution $F$ over $\Pcal$ (a mixed strategy), and player 2 chooses $\ub$ or some distribution $Q$ over $\Ucal$. Observe that by assumption, $g_1$ and $g_2$ are both continuous. Additionally, $g_1$ is quasiconcave in $\pb$ by assumption, and for $g_2$, recall that the negative of a quasiconvex function is quasiconcave, so $g_2$ is also quasiconcave because $R(\pb,\ub)$ is quasiconvex in $\ub$. We now invoke the Debreu-Glicksberg-Fan theorem (see Theorem~1.2 in \citealt{fudenberg1991game}), which asserts that there exists a pure strategy Nash equilibrium for this game. In particular, there exists a $(\pb^*, \ub^*) \in \Pcal \times \Ucal$ such that
\begin{align*}
g_1(\pb^*, \ub^*) \geq g_1(\pb, \ub^*), \quad \forall\ \pb \in \Pcal, \\
g_2(\pb^*, \ub^*) \geq g_2(\pb^*, \ub), \quad \forall\ \ub \in \Ucal.
\end{align*}
Note that this is equivalent to $(\pb^*, \ub^*)$ being a saddle point of $R(\pb,\ub)$:
\begin{equation}
R(\pb, \ub^*) \leq R(\pb^*, \ub^*) \leq R(\pb^*, \ub), \quad \forall\ \pb \in \Pcal, \ub \in \Ucal.
\end{equation}
This implies that 
\begin{align*}
R(\pb^*, \ub^*) & \leq \min_{\ub \in \Ucal} R(\pb^*, \ub) \\
& \leq \max_{\pb \in \Pcal} \min_{\ub \in \Ucal} R(\pb, \ub),
\end{align*}
\begin{align*}
R(\pb^*, \ub^*) & \geq \max_{\pb \in \Pcal} R(\pb, \ub^*) \\
& \geq \min_{\ub \in \Ucal} \max_{\pb \in \Pcal} R(\pb, \ub),
\end{align*}
which together imply that 
\begin{equation}
\min_{\ub \in \Ucal} \max_{\pb \in \Pcal} R(\pb, \ub) \leq \max_{\pb \in \Pcal} \min_{\ub \in \Ucal} R(\pb, \ub). \label{eq:minmax_leq_maxmin}
\end{equation}
Observe now that 
\begin{align*}
Z^*_{\RR} & = \max_{F \in \Fcal} \min_{\ub \in \Ucal} \int_{\Pcal} R(\pb, \ub) \, dF(\pb) \\
& \leq \min_{\ub \in \Ucal} \max_{F \in \Fcal}  \int_{\Pcal} R(\pb, \ub) \, dF(\pb) \\
& = \min_{\ub \in \Ucal} \max_{\pb \in \Pcal}  R(\pb, \ub) \\
& \leq \max_{\pb \in \Pcal} \min_{\ub \in \Ucal} R(\pb, \ub) \\
& = Z^*_{\DR},
\end{align*}
and hence $Z^*_{\RR} = Z^*_{\DR}$. \Halmos

\subsection{Proof of Theorem~\ref{theorem:finite_Pcal}}
\label{proof_theorem:finite_Pcal}

%\begin{proofvvm}
To establish this result, observe that 
\begin{align*}
& \inf_{Q \in \Qcal} \max_{\pb \in \Pcal} \int_{\Ucal} R(\pb, \ub) \, dQ(\ub) \\
& = \inf_{Q \in \Qcal} \max_{\pib \in \Delta_{\Pcal}} \sum_{\pb \in \Pcal} \pi(\pb) \cdot \int_{\Ucal} R(\pb, \ub) \, dQ(\ub) \\
& = \max_{\pib \in \Delta_{\Pcal}} \inf_{Q \in \Qcal} \sum_{\pb \in \Pcal} \pi(\pb) \cdot \int_{\Ucal}  R(\pb, \ub) \, dQ(\ub) \\
& = \max_{\pib \in \Delta_{\Pcal}} \inf_{Q \in \Qcal} \int_{\Ucal} \sum_{\pb \in \Pcal} \pi(\pb) \cdot R(\pb, \ub) \, dQ(\ub) \\
& = \max_{\pib \in \Delta_{\Pcal}} \min_{\ub \in \Ucal} \sum_{\pb \in \Pcal} \pi(\pb) \cdot R(\pb, \ub) \\ 
& = Z^*_{\RR}.
\end{align*}
In the above, the steps are justified as follows. The first step follows because maximizing a function of $\pb$ over the finite set $\Pcal$ is the same as maximizing the expected value of that same function with respect to all probability mass functions $\pib$ supported on $\Pcal$. 

The second step follows by Sion's minimax theorem, because the quantity $\sum_{\pb \in \Pcal} \pi(\pb) \cdot \int_{\Ucal} R(\pb, \ub) \, dQ(\ub)$ is linear and therefore quasiconcave in $\pib$ for a fixed $Q$, and is linear and therefore quasiconvex in $Q$ for a fixed $\pib$; additionally, the set $\Delta_{\Pcal} = \{\pib \in \Rbb^{|\Pcal|} \mid \oneb^T \pib = 1, \pib \geq \zerob\}$ is a compact convex set, and $\Qcal$ is a convex set. Additionally, note that $\sum_{\pb \in \Pcal} \pi(\pb) \cdot \int_{\Ucal} R(\pb, \ub) \, dQ(\ub)$ is clearly continuous in $\pib$. It is also continuous in $Q$, because each term $\int_{\Ucal}  R(\pb, \ub) \, dQ(\ub)$ is continuous in $Q$ when $\Qcal$ is endowed with the topology of weak convergence, and there are finitely many such terms.

The third step follows by the linearity of integration. The fourth step follows by the fact that $\Qcal$ contains the Dirac delta distribution that places unit mass on $\ub$, for every $\ub \in \Ucal$. The final step just follows from the definition of $Z^*_{\RR}$. 

With this result in hand, observe that the existence of a $Q \in \Qcal$ such that for all $\pb \in \Pcal$, $\int_{\Ucal} R(\pb, \ub) dQ(\ub) \leq Z^*_{\DR}$ is equivalent to the existence of $Q \in \Qcal$ such that 
\begin{align*}
\max_{\pb \in \Pcal} \int_{\Ucal} R(\pb, \ub) dQ(\ub) \leq Z^*_{\DR},
\end{align*}
which is equivalent to 
\begin{align*}
\inf_{Q \in \Qcal} \max_{\pb \in \Pcal} \int_{\Ucal} R(\pb, \ub) dQ(\ub) \leq Z^*_{\DR}.
\end{align*}
Since the left-hand side of this inequality is equal to $Z^*_{\RR}$, the existence of the distribution $Q \in \Qcal$ as in the theorem statement is equivalent to $Z^*_{\RR} \leq Z^*_{\DR}$; since it is always the case that $Z^*_{\RR} \geq Z^*_{\DR}$, this is equivalent to the problem being randomization-proof. \hfill \Halmos 
%\end{proofvvm}

\subsection{Proof of Corollary~\ref{corollary:finite_Pcal_strong_duality}}
\label{proof_corollary:finite_Pcal_strong_duality}

Observe that since $\max_{\pb \in \Pcal} \min_{\ub \in \Ucal} R(\pb, \ub) \leq \min_{\ub \in \Ucal} \max_{\pb \in \Pcal}  R(\pb, \ub)$ always holds, equation~\eqref{eq:strong_duality} is equivalent to 
\begin{align*}
\max_{\pb \in \Pcal} \min_{\ub \in \Ucal} R(\pb, \ub) \geq \min_{\ub \in \Ucal} \max_{\pb \in \Pcal}  R(\pb, \ub),
\end{align*}
or equivalently,
\begin{align*}
Z^*_{\DR} \geq \min_{\ub \in \Ucal} \max_{\pb \in \Pcal}  R(\pb, \ub).
\end{align*}
Observe that the condition $\min_{\ub \in \Ucal} \max_{\pb \in \Pcal} R(\pb, \ub) \leq Z^*_{\DR}$ is exactly equivalent to the condition that there exists a $\ub \in \Ucal$ such that for all $\pb \in \Pcal$, $R(\pb, \ub) \leq Z^*_{\DR}$. 

To connect this to Theorem~\ref{theorem:finite_Pcal}, consider $Q = \delta_{\ub}$, where $\delta_{\ub}$ is the Dirac delta distribution centered at $\ub$. For any $\pb \in \Pcal$, $R(\pb, \ub) = \int_{\Ucal} R(\pb, \ub') \, dQ(\ub')$. Thus, for this choice of $Q$, it is the case that for all $\pb \in \Pcal$, $\int_{\Ucal} R(\pb, \ub') \, dQ(\ub') \leq Z^*_{\DR}$. By Theorem~\ref{theorem:finite_Pcal}, this is equivalent to randomization-proofness. Thus, it follows that the strong duality condition~\eqref{eq:strong_duality} implies that the RPO problem is randomization-proof. \hfill \Halmos

\subsection{Proof of Corollary~\ref{corollary:finite_Pcal_pDR}}
\label{proof_corollary:finite_Pcal_pDR}

\noindent \emph{Proof of a)}: Observe that by the definition of $\ub^*$, we have that $R(\pb^*_{\DR}, \ub^*) \leq R(\pb^*_{\DR}, \ub)$ for all $\ub \in \Ucal$. By the hypothesis that $\pb^*_{\DR}$ solves the nominal problem under $\ub^*$, that is, $\pb^*_{\DR} \in \arg \max R(\pb, \ub^*)$, we also have that $R(\pb^*_{\DR}, \ub^*) \geq R(\pb, \ub^*)$ for all $\pb \in \Pcal$. This automatically implies that $(\pb^*_{\DR}, \ub^*)$ is a saddle point of $R$, which implies that strong duality holds, i.e.,
\begin{equation}
\max_{\pb \in \Pcal} \min_{\ub \in \Ucal} R(\pb, \ub) = \min_{\ub \in \Ucal} \max_{\pb \in \Pcal}  R(\pb, \ub),
\end{equation}
and thus that the problem is randomization-proof. 

\noindent \emph{Proof of b)}: By using Theorem~\ref{theorem:finite_Pcal}, we know that there exists a distribution $Q$ such that for every $\pb \in \Pcal$, 
\begin{equation}
\int_{\Ucal} R(\pb, \ub) \, dQ(\ub) \leq Z^*_{\DR}. \label{eq:exists_Q_forall_p_leq_Z_DR}
\end{equation}
Recall that $Z^*_{\DR} = \min_{\ub \in \Ucal} R(\pb^*_{\DR}, \ub) = R(\pb^*_{\DR}, \ub^*)$. Since the minimizer of this problem is unique, we can examine the above condition when $\pb = \pb^*_{\DR}$: 
\begin{equation}
\int_{\Ucal} R(\pb^*_{\DR}, \ub) \, dQ(\ub) \leq \min_{\ub \in \Ucal} R(\pb^*_{\DR}, \ub). \label{eq:exists_Q_forall_p_leq_Z_DR}
\end{equation}
As a result, we must have that $Q = \delta_{\ub^*}$, as any other $Q$ would result in $\int_{\Ucal} R(\pb^*_{\DR}, \ub) \, dQ(\ub)$ being strictly more than $R(\pb^*_{\DR}, \ub^*)$, due to the uniqueness of $\ub^*$.

With this insight in hand, observe that for any $\pb \in \Pcal$, we have
\begin{align*}
\int_{\Ucal} R(\pb, \ub) \, dQ(\ub) = R(\pb, \ub^*).
\end{align*}
Using the property of $Q$ in \eqref{eq:exists_Q_forall_p_leq_Z_DR} at any arbitrary $\pb$ and the fact that $Z^*_{\DR} = R(\pb^*_{\DR}, \ub^*)$, we have
\begin{equation}
R(\pb, \ub^*) \leq R(\pb^*_{\DR}, \ub^*)
\end{equation}
for all $\pb$, which completes the proof. \Halmos

\subsection{Equivalence of strong duality and randomization-proofness when $\Pcal$ is finite}
\label{proof_proposition:equivalence_SD_randproof_finitePcal}

In this section, we develop the following result, which provides a condition under which strong duality and randomization-proofness are equivalent. % (one holds if and only if the other holds). 

\begin{proposition}
Suppose that $\Pcal$ is finite, $R(\pb, \ub)$ is convex in $\ub$ and $\Ucal$ is a convex set. Then the robust price optimization problem is randomization-proof if and only if strong duality holds. \label{proposition:equivalence_SD_randproof_finitePcal}
\end{proposition}

\begin{proofvvm}
Observe that 
\begin{align*}
Z^*_{\RR} & = \max_{\pib \in \Delta_{\Pcal}} \min_{\ub \in \Ucal} \sum_{\pb \in \Pcal} \pi_{\pb} R(\pb, \ub) \\
& = \min_{\ub \in \Ucal} \max_{\pib \in \Delta_{\Pcal}} \sum_{\pb \in \Pcal} \pi_{\pb} R(\pb, \ub) \\
& = \min_{\ub \in \Ucal} \max_{\pb \in \Pcal} R(\pb, \ub),
\end{align*}
where the second equality follows by Sion's minimax theorem, which follows because (1) $\sum_{\pb \in \Pcal} \pi_{\pb} R(\pb, \ub)$ is linear (and hence continuous and quasiconcave) in $\pib$, and convex (and hence continuous and quasiconvex) in $\ub$, (2) $\Delta_{\Pcal}$ is a compact convex set, and (3) $\Ucal$ is a convex set; and the third equality follows because the inner maximum is attained at a deterministic distribution that places all mass on a single $\pb$. Observe that this last min-max expression is an upper bound on $Z^*_{\DR} = \max_{\pb \in \Pcal} \min_{\ub \in \Ucal} R(\pb, \ub)$ by weak duality, and thus that the problem is exactly randomization-proof if and only if strong duality holds. \Halmos
\end{proofvvm}

\subsection{Proof of Proposition~\ref{proposition:loglog_extremal}}
\label{proof_proposition:loglog_extremal}
Let $(\mub^*, \pb^*)$ be an optimal solution of problem~\eqref{prob:loglog_maxmax}, which is guaranteed to exist because the objective function is continuous and the feasible set $\Delta_{[I]} \times \Pcal$ is compact. Observe that the objective function in \eqref{prob:loglog_maxmax} can be re-arranged as 
\begin{align*}
& \max_{\mub \in \Delta_{[I]}, \pb \in \Pcal} \left\{ \sum_{i=1}^I \mu_i (\alpha_i + \log p_i - \beta_i \log p_i + \sum_{j\neq i} \gamma_{i,j} \log p_j) - \sum_{i=1}^I \mu_i \log \mu_i \right\} \\
& = \max_{\mub \in \Delta_{[I]}} \max_{\pb \in \Pcal} \left\{ \sum_{i = 1}^I \mu_i \cdot \alpha_i + \sum_{i=1}^I  \left[ \mu_i \cdot (1 - \beta_i) + \sum_{j \neq i} \mu_j \gamma_{j,i} \right] \cdot \log p_i - \sum_{i=1}^I \mu_i \log \mu_i \right\} \\
& = \max_{\mub \in \Delta_{[I]}} \left[ \sum_{i = 1}^I \mu_i \cdot \alpha_i + \sum_{i=1}^I \max_{p_i \in \Pcal_i} \left\{  \left[ \mu_i \cdot (1 - \beta_i) + \sum_{j \neq i} \mu_j \gamma_{j,i} \right] \cdot \log p_i \right\} - \sum_{i=1}^I \mu_i \log \mu_i  \right]
\end{align*}
where the first step follows by algebra, and the second by the separability of the objective in $p_1,\dots, p_I$ and Assumption~\ref{assumption:Pcal_Cartesian_product} (since the price set is a Cartesian product and the objective is separable, each product's price can be optimized independently). Consider now holding $\mub$ fixed, and optimizing over $\pb$. With $\mub^*$ fixed, the optimal value of $p_i$ for the above objective depends on the sign of $(1 - \beta_i + \sum_{j \neq i} \gamma_{j,i})$. If this coefficient is positive, then since $\log p_i$ is increasing in $p_i$, it is optimal to set $p'_i = \max \Pcal_i $. If this coefficient is negative, then it is optimal to set $p'_i = \min \Pcal_i $. It thus follows that if $(\mub^*, \pb^*)$ is an optimal solution, then $(\mub^*, \pb')$ is also an optimal solution, where $\pb' \in \prod_{i=1}^I \{ \min \Pcal_i, \max \Pcal_i \}$.  \Halmos

\subsection{Proof of Theorem~\ref{theorem:complexity_convexUcal}}
\label{proof_theorem:complexity_convexUcal}

Let $t_{\floor}$ be any real number that satisfies
\begin{equation}
t_{\floor} < \min_{\pb \in \Pcal} \min_{\ub \in \Ucal} R(\pb, \ub). 
\end{equation}
Note that the right-hand side is well-defined: $R(\pb, \ub)$ is convex and therefore continuous in $\ub$, and $\Ucal$ is a compact set, and hence $\min_{\ub \in \Ucal} R(\pb, \ub)$ is finite by the extreme value theorem. Since $\Pcal$ is finite, the quantity $\min_{\pb \in \Pcal} \min_{\ub \in \Ucal} R(\pb, \ub)$ is therefore finite.

Observe that the RRPO problem $\max_{\pb \in \Pcal} \min_{\ub \in \Ucal} R(\pb, \ub)$ is equivalent to 
\begin{subequations}
\begin{alignat}{2}
\RR(\Pcal): \quad & \underset{t,\ub}{\text{minimize}} & \quad & t \\
& \text{subject to} & & t \geq R(\pb, \ub), \quad \forall \pb \in \Pcal, \\
& & & t \geq t_{\floor}, \\
& & & \ub \in \Ucal,
\end{alignat}
\end{subequations}
Note that this formulation is almost the same as formulation~\eqref{prob:finiteP_convexU_epigraph} in Section~\ref{subsec:finitePcal_convexUcal_algorithm}, with the exception of the constraint $t \geq t_{\floor}$; we shall discuss the reason for including this constraint momentarily, after problem~\eqref{prob:RRPO_discreteP_barPcal} is presented. Note that since $t_{\floor}$ by construction is always strictly less than $R(\pb, \ub)$ for any $\pb \in \Pcal$ and any $\ub \in \Ucal$, this constraint is vacuous and cannot become active/binding (i.e., hold at equality). 

Consider now the following generalization of the above problem, where we enforce the constraint $t \geq R(\pb, \ub)$ at a specific set of price vectors, $\bar{\Pcal} \subseteq \Pcal$: 
\begin{subequations}
\begin{alignat}{2}
\RR(\bar{\Pcal}): \quad & \underset{t,\ub}{\text{minimize}} & \quad & t \\
& \text{subject to} & & t \geq R(\pb, \ub), \quad \forall \pb \in \bar{\Pcal}, \\
& & & t \geq t_{\floor}, \\
& & & \ub \in \Ucal,
\end{alignat}
\label{prob:RRPO_discreteP_barPcal}%
\end{subequations}
Let us denote the optimal value by $Z^*_{\RR}(\bar{\Pcal})$, and the problem itself with this set of constraints by $\RR(\bar{\Pcal})$. Similarly to the previous formulation, which is $\RR(\Pcal)$ and has objective value $Z^*_{\RR}(\Pcal)$, the constraint $t \geq t_{\floor}$ in $\RR(\bar{\Pcal})$ cannot be active so long as $\bar{\Pcal}$ is nonempty. Here, we note that the reason for having the constraint $t \geq t_{\floor}$ is to avoid a pathological case where the problem becomes unbounded if $\bar{\Pcal}$ is empty; the inclusion of this constraint facilitates our analysis of the constructive procedure (Algorithm~\ref{algorithm:convexU_constructive_procedure}) that we will discuss shortly (see, in particular, the proof of Lemma~\ref{lemma:convexU_constructive_procedure_barPcal_never_empty}). 

With these definitions, to now appreciate our high-level proof strategy, observe that if we could succeed in showing that $Z^*_{\RR}(\Pcal) = Z^*_{\RR}(\Pcal')$ for a small set of prices $\Pcal' \subseteq \Pcal$ that satisfies $|\Pcal'| \leq M + 1$, then a similar application of Sion's minimax theorem used to derive the convex program in Section~\ref{subsec:finitePcal_convexUcal_algorithm} would allow us to show the following:
\begin{align}
& \max_{\pib \in \Delta_{\Pcal}} \min_{\ub \in \Ucal} \sum_{\pb \in \Pcal} \pi_{\pb} R(\pb, \ub) \label{eq:complexity_proof_chain_first} \\
& = Z^*_{\RR}(\Pcal) \\
& = Z^*_{\RR}(\Pcal') \\
& = \min_{\ub \in \Ucal} \max_{\pb \in \Pcal'} R(\pb, \ub) \\
& = \min_{\ub \in \Ucal} \max_{\pib \in \Delta_{\Pcal'}} \sum_{\pb \in \Pcal'} \pi_{\pb} R(\pb, \ub) \\
& = \max_{\pib \in \Delta_{\Pcal'}} \min_{\ub \in \Ucal} \sum_{\pb \in \Pcal'} \pi_{\pb} R(\pb, \ub), \label{eq:complexity_proof_chain_last}
\end{align}
which would imply that there exists a price distribution $\tilde{\pib}$ over $\Pcal'$ that achieves the same objective as the optimal solution of the RRPO problem when one considers all price vectors. In what follows, we will constructively establish that such a set $\Pcal'$ exists. 

Before proceeding, we first recall the concept of a support constraint from the paper of \cite{calafiore2005uncertain}.

\begin{definition}
Consider the convex optimization problem 
\begin{subequations}
\begin{alignat}{2}
\CONV: \quad & \underset{\xb \in \Rbb^n}{\text{minimize}} & \quad & \cb^T \xb \\
& \text{subject to} & & \xb \in \Xcal_i, \quad \forall i \in [m],
\end{alignat}
\end{subequations}
where $m$ and $n$ are positive integers, and each $\Xcal_i$ is a closed convex set. Consider the problem $\CONV_k$, defined as 
\begin{subequations}
\begin{alignat}{2}
\CONV_k: \quad & \underset{\xb \in \Rbb^n}{\text{minimize}} & \quad & \cb^T \xb \\
& \text{subject to} & & \xb \in \Xcal_i, \quad \forall i \in [m] \setminus \{ k \},
\end{alignat}
\end{subequations}
which is the problem obtained from removing the $k$th constraint $\xb \in \Xcal_k$ from $\CONV$. Let $\xb^*$ and $\xb^*_k$ be optimal solutions of $\CONV$ and $\CONV_k$, respectively. We say that constraint $k$ is a \emph{support constraint}, if $\cb^T \xb^* > \cb^T \xb^*_k$, i.e., if removing constraint $k$ from $\CONV$ causes a change to the optimal objective value of the problem.
\end{definition}

We now recall a critical result about support constraints from \cite{calafiore2005uncertain}. 
\begin{theorem}[Theorem 2 of \cite{calafiore2005uncertain}]
The number of support constraints for problem $\CONV$ is at most $n$. 
\end{theorem}

In our case, note that each inequality $t \geq R(\pb, \ub)$ in $\RR(\bar{\Pcal})$ can be equivalently regarded as the constraint $(t, \ub) \in \{ (t', \ub') \in \Rbb^{M+1} \mid t' \geq R(\pb, \ub') \}$, and the set in this latter definition is clearly a closed, convex set. 

With a slight abuse of terminology, we shall say that $\pb$ is a \emph{support constraint} for $\RR(\bar{\Pcal})$ if the constraint $(t, \ub) \in  \{ (t', \ub') \in \Rbb^{M+1} \mid t' \geq R(\pb, \ub') \}$ is a support constraint. Note that if $\pb$ is a support constraint, then the convex inequality $t \geq R(\pb, \ub)$ must hold at equality for any optimal solution of $\RR(\bar{\Pcal})$. (This is true, because if there were an optimal solution $(t, \ub)$ at which the constraint for $\pb$ was non-binding, then one could take the improved optimal solution $(t', \ub')$ of $\RR(\bar{\Pcal} \setminus \{\pb \})$ and take a convex combination with $(t, \ub)$ to obtain a better solution to $\RR(\bar{\Pcal})$.)

Additionally, given an optimal solution $(t^*,\ub^*)$, we shall say that $\pb$ is an \emph{active non-support constraint} for $\RR(\bar{\Pcal})$ if the constraint $t \geq R(\pb, \ub)$ is active at $(t^*,\ub^*)$, i.e., if $t^* = R(\pb, \ub^*)$, and $\pb$ is not a support constraint of $\RR(\bar{\Pcal})$. 

Lastly, given an optimal solution $(t^*,\ub^*)$, we shall say that $\pb$ is a \emph{non-active constraint} for $\RR(\bar{\Pcal})$ if the constraint $t \geq R(\pb,\ub)$ is not active at $(t^*,\ub^*)$, i.e., if $t^* > R(\pb, \ub^*)$. 

With these definitions, let us consider the following procedure, which we denote as Algorithm~\ref{algorithm:convexU_constructive_procedure}. In words, this algorithm starts with all of the price vectors in $\Pcal$, solves the randomized robust pricing problem with this set, and then decomposes the price vector set into sets corresponding to support constraints, active non-support constraints and non-active constraints. The non-active constraints are also discarded, and if there are any active non-support constraints, one price vector from that set of active non-support constraints is also discarded. The optimization problem is solved again with the remaining smaller set of constraints, and the remaining price vectors are then decomposed again. This procedure repeats until one has a set of price vectors for which all active constraints are support constraints, or equivalently, the set of active non-support constraints $\Ncal$ is empty. 

\begin{algorithm}
\begin{algorithmic}[1]
\STATE Initialize $\bar{\Pcal} = \Pcal$.
\STATE Solve $\RR(\bar{\Pcal})$ to obtain an optimal solution $(t^*, \ub^*)$.
\STATE Set $\Scal \subseteq \bar{\Pcal}$ to be the set of support constraints for $\RR(\bar{\Pcal})$.
\STATE Set $\Ncal \subseteq \bar{\Pcal}$ to be the set of active non-support constraints for $\RR(\bar{\Pcal})$ and $(t^*,\ub^*)$.
\STATE Set $\Ocal \subseteq \bar{\Pcal}$ to be the set of non-active constraints for $\RR(\bar{\Pcal})$ and $(t^*,\ub^*)$.
\WHILE{$|\Ncal| > 0$}
	\STATE Select any $\pb \in \Ncal$, and set $\bar{\Pcal} \gets \Ncal \setminus \{\pb\} \cup \Scal$. 
	\STATE Solve $\RR(\bar{\Pcal})$ to obtain an optimal solution $(t^*, \ub^*)$.
	\STATE Set $\Scal \subseteq \bar{\Pcal}$ to be the set of support constraints for $\RR(\bar{\Pcal})$.
	\STATE Set $\Ncal \subseteq \bar{\Pcal}$ to be the set of active non-support constraints for $\RR(\bar{\Pcal})$ and $(t^*,\ub^*)$
	\STATE Set $\Ocal \subseteq \bar{\Pcal}$ to be the set of non-active constraints for $\RR(\bar{\Pcal})$ and $(t^*,\ub^*)$
\ENDWHILE
\RETURN $\Scal^* \equiv \Scal$. 
\end{algorithmic}
\caption{Constructive procedure for reducing the price vector set $\Pcal$. \label{algorithm:convexU_constructive_procedure} }
\end{algorithm}

We now establish a few key properties of this procedure.
\begin{lemma}
Algorithm~\ref{algorithm:convexU_constructive_procedure} terminates in a finite number of iterations. \label{lemma:convexU_constructive_procedure_termination}
\end{lemma}
\begin{proofvvm}
Observe that after line 7 of Algorithm~\ref{algorithm:convexU_constructive_procedure}, $|\bar{\Pcal}| = |\Ncal | - 1 + |\Scal|$, whereas immediately before line 7, $|\bar{\Pcal}| = |\Ncal| + | \Scal | + |\Ocal|$, which is strictly higher than $|\Ncal| - 1 + |\Scal|$. Since the size of $\bar{\Pcal}$ is decreasing by at least 1 in each iteration, the initial set $\Pcal$ is finite, and $\Ncal$ is always updated to be a subset of $\bar{\Pcal}$, it must be that $\Ncal$ eventually becomes empty. \Halmos
\end{proofvvm}

%The second is that $\bar{\Pcal}$ never becomes empty during the execution of the procedure.

\begin{lemma}
The set $\bar{\Pcal}$ is non-empty throughout the execution of Algorithm~\ref{algorithm:convexU_constructive_procedure}. \label{lemma:convexU_constructive_procedure_barPcal_never_empty}
\end{lemma}
\begin{proofvvm}
Suppose for the sake of contradiction that $\bar{\Pcal}$ becomes empty. This would mean that in some iteration of the while loop, at line 7, one would have that $\Ncal \setminus \{\pb\} \cup \Scal$ is empty, which would imply that $\Ncal \setminus \{\pb\} = \emptyset$ and $\Scal = \emptyset$. Since $\pb$ is chosen from $\Ncal$, this would imply that $\Ncal = \{\pb \}$. This means that prior to line 7, $\RR(\bar{\Pcal})$ has only one active constraint, which is at $\pb$, and that constraint is not a support constraint. This, however, leads to a contradiction, because by removing this constraint, one would be free to set $t = t_{\floor}$, and recalling that $t_{\floor}$ is strictly smaller than $R(\pb,\ub)$ for all $\ub \in \Ucal$, this would mean we are able to improve the objective value. (In other words, $\pb$ has to be a support constraint.) \Halmos
\end{proofvvm}

%The third is that the procedure outputs a nonempty set of price vectors.

\begin{lemma}
Let $\Scal^*$ be the output of Algorithm~\ref{algorithm:convexU_constructive_procedure}. Then $\Scal^*$ is not empty, and $|\Scal^*| \leq M + 1$. \label{lemma:Scal_nonempty_Mplus1}
\end{lemma}
\begin{proofvvm}
Observe that Algorithm~\ref{algorithm:convexU_constructive_procedure} terminates by Lemma~\ref{lemma:convexU_constructive_procedure_termination}, and that by Lemma~\ref{lemma:convexU_constructive_procedure_barPcal_never_empty}, $\bar{\Pcal}$ never becomes empty during Algorithm~\ref{algorithm:convexU_constructive_procedure}. Consider now what happens in lines 9 and 10 of Algorithm~\ref{algorithm:convexU_constructive_procedure} immediately prior to termination. Since $\bar{\Pcal}$ is not empty, then $\RR(\bar{\Pcal})$ must have at least one active constraint from among the constraints in $\bar{\Pcal}$ (i.e., $\Scal \cup \Ncal$ is not empty); since the algorithm terminates and the loop exits, $\Ncal$ must be empty. This implies that $\Scal^*$ is not empty. 

To see that $|\Scal^*| \leq M + 1$, observe that $\Scal^*$ is the set of support constraints of $\RR(\bar{\Pcal})$ for some $\bar{\Pcal} \subseteq \Pcal$, which is a $(M + 1)$-dimensional convex optimization problem; by Theorem~2 of \cite{calafiore2005uncertain}, $|\Scal^*|$ can be at most $M + 1$. \Halmos
\end{proofvvm}

\begin{lemma}
Let $\Scal^*$ be the output of Algorithm~\ref{algorithm:convexU_constructive_procedure}. Then $\RR(\Scal^*) = \RR(\Pcal)$. \label{lemma:Scal_equals_Pcal}
\end{lemma}
\begin{proofvvm}
We shall establish this result by induction. Let $\bar{\Pcal} \subseteq \Pcal$ be any non-empty set of price vectors, and let $\Scal$, $\Ncal$ and $\Ocal$ be defined as in lines 9, 10, and 11 respectively. Let $(t^*, \ub^*)$ be any optimal solution of $\RR(\bar{\Pcal})$. We claim that 
\begin{equation}
Z^*_{\RR}(\bar{\Pcal}) \underset{(a)}{=} Z^*_{\RR}(\Scal \cup \Ncal) \underset{(b)}{=} Z^*_{\RR}(\Scal \cup \Ncal \setminus \{ \pb' \}), \label{eq:inductive_step}
\end{equation}
where $\pb'$ is any price vector in $\Ncal$. We now establish each of the two equalities, (a) and (b). \\

\noindent \emph{Proof of (a).} Observe that we trivially have $Z^*_{\RR}(\bar{\Pcal}) \geq Z^*_{\RR}(\Scal \cup \Ncal)$, because $\Scal \cup \Ncal$ is a subset of $\bar{\Pcal}$, and so the objective value of the minimization problem with $\Scal \cup \Ncal$ is either equal or lower than with $\bar{\Pcal}$. 

To see why $Z^*_{\RR}(\bar{\Pcal}) > Z^*_{\RR}(\Scal \cup \Ncal)$ cannot hold, suppose for the sake of contradiction that this was the case, and let $(t', \ub')$ be any optimal solution of $\RR(\Scal \cup \Ncal)$. Observe that by taking any convex combination of $(t^*, \ub^*)$ and $(t', \ub')$, we would obtain a new solution $(\tilde{t}, \tilde{\ub})$ that satisfies the constraint $t \geq R(\pb, \ub)$ for all $\pb \in \Scal \cup \Ncal$, because the function $R(\pb, \ub)$ is convex in $\ub$ and both $(t^*, \ub^*)$ and $(t', \ub')$ satisfy those constraints. Note also that while $(t', \ub')$ might not necessarily satisfy $t \geq R(\pb, \ub)$ for every $\pb \in \Ocal = \bar{\Pcal} \setminus (\Scal  \cup \Ncal)$, by a suitable choice of the weight on $(t^*, \ub^*)$, the convex combination $(\tilde{t}, \tilde{\ub})$ can be made to satisfy these because by assumption, $t^* > R(\pb, \ub^*)$ for all $\pb \in \Ocal$. (Because of the slack in these constraints for $(t^*, \ub^*)$, and because there are finitely many of these constraints, we can choose the weight in the convex combination so that none of these constraints are violated in the new solution $(\tilde{t}, \tilde{\ub})$.) Lastly, since the objective function is linear (in fact, it is exactly $t$), any convex combination of $t^*$ and $t'$ with a positive weight on $t'$ would ensure that $\tilde{t}$ is strictly smaller than $t^*$. Since we can construct a new solution $(\tilde{t}, \tilde{\ub})$ that is feasible for $\RR(\bar{\Pcal})$ but has an objective value smaller than $t^*$, we reach a contradiction, because $(t^*, \ub^*)$ was assumed to be optimal for $\RR(\bar{\Pcal})$. \\

\noindent \emph{Proof of (b).} We now establish the second equality. We shall do this by showing that $Z^*_{\RR}(\bar{\Pcal}) = Z^*_{\RR}(\Scal \cup \Ncal \setminus \{\pb'\})$; by then using (a), we will have that $Z^*_{\RR}(\Scal \cup \Ncal) = Z^*_{\RR}(\Scal \cup \Ncal \setminus \{ \pb' \})$. 

As in the case of (a), we again trivially have that $Z^*_{\RR}(\bar{\Pcal}) \geq Z^*_{\RR}(\Scal \cup \Ncal \setminus \{ \pb' \})$, because $\Scal \cup \Ncal \setminus \{\pb' \} \subseteq \bar{\Pcal}$.

%Again, we trivially have that $Z_{\RR}(\Scal \cup \Ncal) \geq Z_{\RR}(\Scal \cup \Ncal \setminus \{ \pb' \})$, because $\Scal \cup \Ncal \setminus \{\pb\} \subseteq \Scal \cup \Ncal$. 

To see why $Z^*_{\RR}(\bar{\Pcal}) > Z^*_{\RR}(\Scal \cup \Ncal \setminus \{ \pb' \})$ cannot hold, let us suppose that it does, and let $(t', \ub')$ be any optimal solution of $\RR(\Scal \cup \Ncal \setminus \{ \pb'\} )$. We thus have $t^* > t'$. Now, consider a new solution $(\tilde{t}, \tilde{\ub})$, which is obtained by taking a convex combination of $(t^*, \ub^*)$ and $(t', \ub')$. Observe that any convex combination of $(t^*, \ub^*)$ and $(t', \ub')$ satisfies $t \geq R(\pb, \ub)$ for all $\pb \in \Scal \cup \Ncal \setminus \{ \pb' \}$, because both $(t^*, \ub^*)$ and $(t', \ub')$ have to satisfy these constraints and $R(\pb, \ub)$ is convex in $\ub$. Observe also that because the constraints in $\Ocal$ are not binding at $(t^*, \ub^*)$ and $\Ocal$ is finite, we can choose the positive weight on $(t', \ub')$ in the convex combination so that $(\tilde{t}, \tilde{\ub})$ satisfies $t \geq R(\pb, \ub)$ for all $\pb \in \Ocal$ (this is similar to the argument in the case of equality (a)). Lastly, by choosing a suitable positive weight on $(t', \ub')$, the new solution $(\tilde{t}, \tilde{\ub})$ can also be made to have an objective value strictly lower than $t^*$. Thus, by choosing a suitable weight on $(t', \ub')$, we are able to construct the new solution $(\tilde{t}, \tilde{\ub})$ which is a feasible solution of $\RR(\bar{\Pcal} \setminus \{\pb' \})$, and has objective value $\tilde{t} < t^* = Z^*_{\RR}(\bar{\Pcal})$. This implies that $Z^*_{\RR}(\bar{\Pcal} \setminus \{ \pb' \}) < Z^*_{\RR}(\bar{\Pcal})$, which immediately leads to a contradiction, because it implies that $\pb'$ is a support constraint of $\RR(\bar{\Pcal})$. (By dropping the $\pb'$ constraint from $\RR(\bar{\Pcal})$, we get $Z^*_{\RR}(\bar{\Pcal} \setminus \{ \pb' \})$ which we have shown to be lower than $Z^*_{\RR}(\bar{\Pcal})$.) This contradicts the assumption that $\pb'$ is from $\Ncal$ and not a support constraint. \\

\noindent Having established equation~\eqref{eq:inductive_step}, suppose that the loop in Algorithm~\ref{algorithm:convexU_constructive_procedure} performs $k$ iterations, and let $\bar{\Pcal}_\ell$ be the set $\bar{\Pcal}$ after line 7 is performed in the $\ell$th iteration of the loop. Recall that $\Scal^*$ is the output of Algorithm~\ref{algorithm:convexU_constructive_procedure} in the hypothesis of the lemma. Then by applying equation~\eqref{eq:inductive_step} inductively, we have that
\begin{equation}
Z^*_{\RR}(\Pcal) = Z^*_{\RR}(\bar{\Pcal}_1) = Z^*_{\RR}(\bar{\Pcal}_2) = Z^*_{\RR}(\bar{\Pcal}_3) = \dots = Z^*_{\RR}(\bar{\Pcal}_k) = Z^*_{\RR}(\Scal^*),
\end{equation}
where each equality except for the last follows by using both (a) and (b) in \eqref{eq:inductive_step}, while the last equality $Z^*_{\RR}(\bar{\Pcal}_k) = Z^*_{\RR}(\Scal^*)$ follows by using (a) in \eqref{eq:inductive_step} together with the fact that $\Ncal = \emptyset$ (this is the case because it is the final iteration and the loop is exited). \Halmos
\end{proofvvm}

We now complete the proof of Theorem~\ref{theorem:complexity_convexUcal}. To establish the result, we apply Algorithm~\ref{algorithm:convexU_constructive_procedure}. The resulting $\Scal^*$ that we obtain is nonempty and has at most $M + 1$ elements (Lemma~\ref{lemma:Scal_nonempty_Mplus1}); moreover, we have that $Z^*_{\RR}(\Scal^*) = Z^*_{\RR}(\Pcal)$ (Lemma~\eqref{lemma:Scal_equals_Pcal}). By our earlier argument (equations~\eqref{eq:complexity_proof_chain_first} - \eqref{eq:complexity_proof_chain_last}), this implies that 
\begin{equation}
\max_{\pib \in \Delta_{\Scal^*}} \min_{\ub \in \Ucal} \sum_{\pb \in \Scal^*} \pi_{\pb} R(\pb, \ub) = \max_{\pib \in \Delta_{\Pcal}} \min_{\ub \in \Ucal} \sum_{\pb \in \Pcal} \pi_{\pb} R(\pb, \ub),
\end{equation}
which establishes the theorem. \Halmos

%%%%% 
% START OF EXAMPLE / COUNTEREXAMPLE SECTION FOR SECTION 4.
%%%%%

\section{Accompanying examples for Section~\ref{sec:benefits}}
\label{sec:examples_benefits}

In this section, we provide a number of examples that accompany our discussion of our theoretical results in Section~\ref{sec:benefits}. Section~\ref{subsec:examples_benefits_concave} provides examples to illustrate that certain conditions in Theorem~\ref{theorem:concave_revenue_function_set} cannot be relaxed, and an example to show that strong duality does not hold in the setting of this theorem. Section~\ref{subsec:examples_benefits_concave_applications} documents a number of examples of applications of Theorem~\ref{theorem:concave_revenue_function_set}. Section~\ref{subsec:examples_benefits_quasiconcave_quasiconvex_applications} documents a number of examples of applications of Theorem~\ref{theorem:quasiconcave_quasiconvex}. Section~\ref{counterexample_corollary:finite_Pcal_pDR} documents an example to show that uniqueness is necessary for part (b) of Corollary~\ref{corollary:finite_Pcal_pDR}, and also to show that randomization-proofness can hold without strong duality holding when $\Pcal$ is finite. Lastly, Section~\ref{subsec:examples_benefits_DR_RR_dual_all_different} documents an example that illustrates that $Z^*_{\DR}$, $Z^*_{\RR}$ and the dual bound $\min_{\ub \in \Ucal} \max_{\pb \in \Pcal} R(\pb,\ub)$ can in general all be different.

\subsection{Accompanying examples for Theorem~\ref{theorem:concave_revenue_function_set}}
\label{subsec:examples_benefits_concave}

As noted in Section~\ref{subsec:benefits_concave}, the condition that all functions in $\Rcal$ be concave cannot be relaxed, which we illustrate using the following example.

\begin{example}
Consider a single-product RPO problem, and suppose that the revenue function uncertainty set $\Rcal = \{R_1, R_2\}$, where $R_1(\cdot)$ and $R_2(\cdot)$ are defined as
\begin{align*}
R_1(p) & = p (10 - 2p), \\
R_2(p) & = p \cdot 10p^{-2} .
\end{align*}
Note that $R_1(p)$ is the revenue function corresponding to the linear demand function $d_1(p) = 10 - 2p$, while $R_2(p)$ is the revenue function corresponding to the log-log demand function $d_2(p) = \exp( \log(10) - 2 \log(p)) = 10p^{-2}$. Note also that $R_1(\cdot)$ is concave, while $R_2(\cdot)$ is convex. Suppose additionally that $\Pcal = [1,4]$.

We first calculate the optimal value of the DRPO problem. Observe that in the interval $[1,4]$, the only root of the equation $10 - 2p = 10p^{-2}$ is $p \approx p' = 1.137805\dots$. For $p < p'$, $d_2(p) > d_1(p)$, and for $p > p'$, $d_1(p) > d_2(p)$. Therefore, the optimal value of the of the DRPO problem can be calculated as
\begin{align*}
& \max_{p \in [1,4]} \min_{R \in \Rcal} R(p) \\
& = \max\{ \max_{p \in [1, p']} \min_{R \in \Rcal} R(p), \max_{p \in [p', 4]} \min_{R \in \Rcal} R(p) \} \\
& = \max\{ \max_{p \in [1, p']} p \cdot (10 - 2p), \max_{p \in [p', 4]}  p \cdot 10 p^{-2} \} \\
& = 10 p' - 2 p'^2 \\ 
& = 8.78885
\end{align*}
In the above, the first step follows because the best value of the worst-case revenue over $[1,4]$ is equivalent to taking the higher of the best worst-case revenue over either $[1,p']$ or $[p', 4]$. The second step follows because for every $p \in [1,p']$, $d_1(p) < d_2(p)$, and so $R_1(p) = p \cdot d_1(p)  < p \cdot d_2(p) = R_2(p)$; similarly, for every $p \in [p', 4]$, $d_1(p) > d_2(p)$, and so $R_2(p) < R_1(p)$. The third step follows by carrying out the maximization of each of the two functions from the prior step over its corresponding interval. 

Now, let us lower bound the optimal value of the RRPO problem. Consider a distribution $F$ that randomizes over prices in the following way:
\begin{equation}
p = \left\{ \begin{array}{ll} 1 & \text{with probability}\ 17/21, \\
2.5 & \text{with probability}\ 4/21. \end{array}  \right.
\end{equation}
The worst-case revenue for this distribution is 
\begin{align*}
& \min_{R \in \Rcal} \int_{1}^4 R(p) \, dF(p) \\
& = \min\left\{ \frac{17}{21} \cdot R_1(1) + \frac{4}{21} \cdot R_1(2.5),  \frac{17}{21} \cdot R_2(1) + \frac{4}{21} \cdot R_2(2.5) \right\} \\
& = \min\left\{ \frac{17}{21} \cdot 8 + \frac{4}{21} \cdot 12.5,  \frac{17}{21} \cdot 10 + \frac{4}{21} \cdot 4 \right\} \\
& = \min\left\{ \frac{62}{7}, \frac{62}{7} \right\} \\
& = \frac{62}{7} \\
& = 8.857143.
\end{align*}
This implies that $Z^*_{\RR} \geq 8.857143$, whereas $Z^*_{\DR} =  8.78885$, and thus $Z^*_{\RR} > Z^*_{\DR}$. \Halmos
\end{example}

The condition that $\Pcal$ be a convex set also cannot be relaxed in general, which we illustrate using the next example.

\begin{example}
\label{example:concave_revenue_fn_set_convexity_unrelaxable}
Consider again a single-product RPO problem. Suppose that $\Rcal = \{ R_1, R_2\}$, where $R_1(p) = p (10 - p)$, $R_2(p) = p (4 - 0.2p)$; $R_1$ and $R_2$ correspond to linear demand functions $d_1(p) = 10 - p$, $d_2(p) = 4 - 0.2p$. Suppose that $\Pcal = \{p_1, p_2\}$, where $p_1 = 5$, $p_2 = 10$. From this data, observe that:
\begin{align*}
R_1(p_1) & = 5 (10 - 5) = 25, \\
R_1(p_2) & = 10 (10 - 10) = 0, \\
R_2(p_1) & =5 (4 - 0.2 (5) ) = 15, \\
R_2(p_2) &= 10 (4 - 0.2 (10)) = 20.
\end{align*}
We first calculate the optimal value of the DRPO problem:
\begin{align*}
Z^*_{\DR} & = \max_{p \in \{p_1, p_2\}} \min \{ R_1(p), R_2(p) \} \\
& = \max\{ \min \{ 25, 15 \}, \min\{ 0, 20 \} \} \\
& = \max \{ 15, 0 \} \\
& = 15.
\end{align*}
For the RRPO problem, the optimal value is given by the following LP:
\begin{subequations}
\begin{alignat}{2}
& \underset{ \eta, \pib }{\text{maximize}} &\quad  & \eta \\
& \text{subject to} & & \eta \leq \pi_{p_1} \cdot p_1 \cdot (10 - p_1) + \pi_{p_2} \cdot p_2 \cdot (10 - p_2) \\
& & & \eta \leq \pi_{p_1}\cdot  p_1 \cdot(4 - 0.2 p_1) + \pi_{p_2} \cdot p_2 \cdot (4 - 0.2 p_2) \\
& & & \pi_{p_1} + \pi_{p_2} = 1\\
& & & \pi_{p_1}, \pi_{p_2} \geq 0.
\end{alignat}
\end{subequations}
The optimal distribution over $\Pcal = \{ p_1, p_2 \}$ is given by $\pi_{p_1} = 2/3$, $\pi_{p_2} = 1/3$, which leads to $Z^*_{\RR} = 50 / 3 = 16.6667$. Since this is higher than $Z^*_{\DR}$, we conclude that this particular instance is randomization-receptive. \Halmos
\end{example}

%\comebacktothis -- add the strong duality counterexample here??

Lastly, as noted in our discussion after Theorem~\ref{theorem:concave_revenue_function_set}, this result is not being driven by strong duality. The next example illustrates that a RPO problem can be randomization-proof without strong duality holding. (Note that this example is the same as our previous example, Example~\ref{example:concave_revenue_fn_set_convexity_unrelaxable}, with the price set $\{5, 10\}$ replaced by its convex hull, $[5,10]$.)

\begin{example}
\label{example:concave_strong_duality_fails}
Consider a single-product RPO problem, where:
\begin{align}
\Pcal & = [5, 10], \\
R_1(p) & = p(10 - p), \\
R_2(p) & = p(4 - 0.2p).
\end{align}
Observe that $\Pcal$ is convex and both $R_1$ and $R_2$ are concave quadratic functions, and so by Theorem~\ref{theorem:concave_revenue_function_set}, the problem is randomization-proof. However, observe that
\begin{align*}
& \max_{p \in [5,10]} \min \{ R_1(p), R_2(p) \} \\
& = 18.75,
\end{align*}
whereas
\begin{align*}
& \min \{ \max_{p \in [5,10]} R_1(p), \max_{p \in [5,10]} R_2(p) \} \\
& = \min \{ 25, 20 \} \\
& = 20,
\end{align*}
which exactly implies that $\max_{\pb \in \Pcal} \min_{R \in \Rcal} R(\pb) \neq \min_{R \in \Rcal} \max_{\pb \in \Pcal}  R(\pb)$. \Halmos
\end{example}

\subsection{Applications of Theorem~\ref{theorem:concave_revenue_function_set}}
% Theorem 1
\label{subsec:examples_benefits_concave_applications}

Theorem~\ref{theorem:concave_revenue_function_set} allows us to assert that a number of RPO problems are randomization-proof under certain conditions, which we illustrate through the following four examples. \\

\begin{example}
(Single-product pricing under linear demand). Suppose that $I = 1$, which corresponds to a single-product pricing problem. Let $\ub = (\alpha,\beta) \in \Rbb^2$ denote the vector of linear demand model parameters, and let $\Ucal \subseteq \Rbb^2$ be an uncertainty set of possible values of $(\alpha,\beta)$. Let $\Rcal = \{ R(\cdot, \ub) \mid \ub \in \Ucal \}$ be the set of revenue functions that arise from the uncertainty set $\Ucal$. Note that each revenue function is of the form $R(p) = \alpha p - \beta p^2$. Therefore, the condition that each $R \in \Rcal$ is concave implies that $R''(p) = - 2 \beta \leq 0$. Thus, if $\Ucal$ is such that $\inf \{ \beta \mid (\alpha, \beta) \in \Ucal \} \geq 0$, then the robust price optimization problem is randomization-proof. \Halmos \\

%Let $d_{\alpha,\beta}(p) = \alpha - \beta p$ be the demand function given $\alpha, \beta$. Let $\Ucal \subseteq \Rbb^2$ be an uncertainty set of possible values of $(\alpha, \beta)$; let $\Dcal = \{ d_{\alpha, \beta} \mid (\alpha, \beta) \in \Ucal\}$ be the demand function uncertainty set; and let $\Rcal = \{ R: p \mapsto p d(p) \mid d \in \Dcal\}$ be the revenue function uncertainty set. Note that for each $\alpha,\beta$, $R(p) = \alpha p - \beta p^2$. Therefore, the condition that each $R \in \Rcal$ is concave implies that $R''(p) = - 2 \beta \leq 0$. Thus, if $\Ucal$ is such that $\inf \{ \beta \mid (\alpha, \beta) \in \Ucal \} \geq 0$, then the robust price optimization problem is randomization-proof. 
\end{example}

\begin{example}
(Multi-product pricing under linear demand). In the more general multi-product pricing problem, let $\ub = (\alphab, \betab, \gammab) \in \Rbb^{I} \times \Rbb^{I} \times \Rbb^{I(I-1)}$ denote the vector of linear demand model parameters, and let $\Ucal$ be an arbitrary uncertainty set of these model parameter vectors. Let $\Rcal = \{ R(\cdot, \ub) \mid \ub \in \Ucal \}$ be the set of revenue functions that arise from the uncertainty set $\Ucal$. Observe that each revenue function $R(\cdot, \ub)$ is of the form
\begin{align*}
R(\pb) & = \sum_{i=1}^I p_i (\alpha_i - \beta_i p_i + \sum_{j \neq i} \gamma_{i,j} p_j ) \\
& = \alphab^T \pb - \pb^T \Mb_{\betab,\gammab} \pb,
\end{align*}
where $\Mb_{\betab,\gammab}$ is the matrix
\begin{align*}
\Mb_{\betab,\gammab} = \left[ \begin{array}{cccccc}
\beta_1 & -\gamma_{1,2} & -\gamma_{1,3} & \cdots & -\gamma_{1,I-1} & -\gamma_{1,I} \\
 -\gamma_{2,1} & \beta_2 & -\gamma_{2,3} & \cdots & -\gamma_{2,I-1} & -\gamma_{2,I} \\
\vdots & \vdots & \ddots & & \vdots & \vdots \\
 -\gamma_{I,1} & -\gamma_{I,2} & -\gamma_{I,3} & \cdots & -\gamma_{I,I-1} & \beta_I \end{array} \right].
\end{align*}
This implies that 
\begin{align*}
\nabla^2 R(\pb) = - 2 \Mb_{\betab,\gammab}.
\end{align*}
The function $R$ is therefore concave if the matrix $\Mb_{\betab,\gammab}$ is positive semidefinite. Therefore, if $\Ucal$ is such that $\inf \{ \lambda_{\min}( \Mb_{\betab,\gammab}) \mid (\alphab, \betab, \gammab) \in \Ucal \} \geq 0$, where $\lambda_{\min}(\Ab)$ denotes the minimum eigenvalue of a symmetric matrix $\Ab$, then the robust price optimization problem is randomization proof. \Halmos \\

%let $\db_{\alphab, \betab}$ be the vector-valued demand function. Let $\Ucal \subseteq \Rbb^{I} \times \Rbb^{I\times I}$, let $\Dcal = \{ \db_{\alphab, \betab} \mid (\alphab, \betab) \in \Ucal\}$, and let $\Rcal = \{ R(\pb): \pb \mapsto \sum_{i=1}^I p_i \cdot d_{i}(\pb) \mid \db(\cdot) = (d_1(\cdot), \dots, d_I(\cdot)) \in \Dcal \}$. Each revenue function $R \in \Rcal$ is of the form
%\begin{align*}
%R(\pb) & = \sum_{i=1}^I p_i (\alpha_i - \sum_{j=1}^I \beta_{i,j} p_j) \\
%& = \alphab^T \pb - \pb^T \betab \pb 
%\end{align*}
%which implies that 
%\begin{align*}
%\nabla^2 R(\pb) = - 2 \betab.
%\end{align*}
%The function $R$ is therefore concave if the matrix $\betab$ is positive semidefinite. Therefore, if $\Ucal$ is such that $\inf \{ \lambda_{\min}( \betab) \mid (\alphab, \betab) \in \Ucal \} \geq 0$, where $\lambda_{\min}(\Ab)$ denotes the minimum eigenvalue of a symmetric matrix $\Ab$, then the robust price optimization problem is randomization proof. 
 
\end{example}

\begin{example}
(Single-product pricing under semi-log demand). For the single-product pricing problem under semi-log demand, $d(p) = \exp(\alpha - \beta p)$ is the demand function given the parameter vector $\ub = (\alpha, \beta)$. Let $\Ucal \subseteq \Rbb^2$ be an uncertainty set of possible values of $(\alpha, \beta)$, and assume that $\beta$ is bounded away from zero, that is, $\inf \{ \beta \mid (\alpha,\beta) \in \Ucal\} \geq 0$. $\Rcal = \{ R(\cdot, \ub) \mid \ub \in \Ucal \}$ be the revenue function uncertainty set. For a given $R \in \Rcal$, its second derivative is $R''(p) = R''(p)  = \beta ( \beta p - 2 ) e^{\alpha - \beta p}$. 
%we calculate its second derivative as follows:
%\begin{align*}
%R(p) & = p e^{\alpha - \beta p } \\
%R'(p) & = (1 - \beta p ) e^{\alpha - \beta p} \\
%R''(p) & = \beta ( \beta p - 2 ) e^{\alpha - \beta p} 
%\end{align*}
Thus, for $R''(p)$ to be nonpositive, we need $\beta p - 2 \leq 0$ or equivalently $\beta p \leq 2$ (since $\beta$ is assumed to be nonnegative) for all $p \in \Pcal$ in order for $R(p)$ to be concave. Thus, if $\sup_{p \in \Pcal} \sup_{ (\alpha,\beta) \in \Ucal}  \{ \beta p \} \leq 2$ and $\inf \{ \beta \mid (\alpha,\beta) \in \Ucal\} \geq 0$, then the RPO problem is randomization-proof. \Halmos \\
\end{example}

\begin{example}
(Single-product pricing under log-log demand). For the single-product pricing problem under log-log demand, $d(p) = \exp(\alpha - \beta \log p) = e^{\alpha} \cdot p^{-\beta}$ is the demand function and $\ub =  (\alpha, \beta) \in \Rbb^2$ is the vector of uncertain demand model parameters. Let $\Ucal \subseteq \Rbb^2$ be an uncertainty set of possible values of $(\alpha, \beta)$, and assume that $\beta$ is bounded away from zero from below, that is, $\inf \{ \beta \mid (\alpha,\beta) \in \Ucal\} \geq 0$. Let $\Rcal = \{ R(\cdot, \ub) \mid \ub \in \Ucal \} $ be the revenue function uncertainty set. For a given $R \in \Rcal$, its second derivative is $R''(p) = e^{\alpha} \cdot (\beta - 1 ) (\beta) \cdot p^{- \beta - 1}$. 
%we calculate its second derivative as follows:
%\begin{align*}
%R(p) & = p \cdot e^{\alpha} \cdot p^{-\beta} \\
%& = e^{\alpha} \cdot p^{1 - \beta} \\
%R'(p) & = e^{\alpha} \cdot (1 - \beta) \cdot p^{- \beta} \\
%R''(p) & = e^{\alpha} \cdot (1 - \beta) (- \beta) \cdot p^{- \beta} \\ 
%& = e^{\alpha} \cdot (\beta - 1 ) (\beta) \cdot p^{- \beta}.
%\end{align*}
Thus, for $R''(p)$ to be nonpositive, we need $\beta - 1 \leq 0$, or equivalently $\beta \leq 1$. Thus, if $\sup_{ (\alpha, \beta) \in \Ucal} \beta \leq 1$ and $\inf_{(\alpha,\beta) \in \Ucal} \beta \geq 0$, then the RPO problem is randomization-proof. \Halmos
\end{example}

\subsection{Applications of Theorem~\ref{theorem:quasiconcave_quasiconvex}}
% Theorem 2 - quasiconvexity / quasiconcavity. 
\label{subsec:examples_benefits_quasiconcave_quasiconvex_applications}

Theorem~\ref{theorem:quasiconcave_quasiconvex} allows us to establish that a single product RPO problem under semi-log and log-log demand will be randomization-proof, even when the revenue function is not necessarily concave (and one cannot rely on Theorem~\ref{theorem:concave_revenue_function_set}).  \\

\begin{example}
Consider a single-product price optimization problem where the demand follows a semi-log model. The uncertain parameter is therefore $\ub = (\alpha, \beta)$. 

Observe that $R(p,\ub) = p e^{\alpha - \beta p}$ is convex in $\ub$. Thus, it is also quasiconvex in $\ub$ for a fixed $p$. Note also that the function $R$ is quasi-concave in $p$. To see this, observe that $\log R(p, \ub) = \log p + \alpha - \beta p$, 
%\begin{align*}
%\log R(p, \ub) & = \log [ p \cdot e^{\alpha - \beta p} ]\\
%& = \log p + \alpha - \beta p,
%\log R(p, \ub) & = \log p + \alpha - \beta p
%\end{align*}
which is concave in $p$; this means that $R$ is log-concave in $p$. Since any log-concave function is quasiconcave, it follows that $R$ is quasiconcave in $p$. 

Thus, if $\Pcal \subseteq \Rbb$ and $\Ucal \subseteq \Rbb^2$ are compact and convex, then Theorem~\ref{theorem:quasiconcave_quasiconvex} asserts that the RPO problem is randomization-proof. \Halmos \\
\end{example}

\begin{example}
Consider a single-product price optimization problem where the demand follows a log-log model. The uncertain parameter is $\ub = (\alpha, \beta)$, and $R(p, \ub) = p e^{\alpha - \beta \log p}$. Assume that $\Pcal \subseteq \Rbb$ is a compact convex set, and that $\min \{ p \mid p \in \Pcal \} > 0$. 

Observe that $R(p, \ub) = p \cdot e^{\alpha - \beta \log p}$ is convex in $\ub$, and therefore quasiconvex in $\ub$ for a fixed $p$. %Additionally, for any distribution $F$ over $\Pcal$, we have that the function $\int_{\Pcal} R(p, \ub) \, dF(p)$ is convex in $\ub$ and therefore also quasiconvex in $\ub$ for a fixed $F$. \comebacktothis -- quasiconcavity of integral is not necessary anymore. 
Additionally, with regard to quasiconcavity in $p$, observe that $\log R(p, \ub) = \log p + \alpha - \beta \log p = (1 - \beta) \log p + \alpha$, 
%\begin{align*}
%\log R(p, \ub) &=  \log [ p e^{\alpha - \beta \log p} ] \\
%& = \log p + \alpha - \beta \log p \\
%& = (1 - \beta) \log p + \alpha,
%\end{align*}
which means that $R$ is log-concave in $p$ whenever $1 - \beta > 0$ or equivalently $\beta < 1$. Therefore, $R$ will also be quasiconcave whenever $\beta < 1$. 

Thus, if $\Pcal \subseteq \Rbb$ and $\Ucal \subseteq \Rbb^2$ are compact and convex, and $\max \{ \beta \mid (\alpha, \beta) \in \Ucal \} < 1$, then Theorem~\ref{theorem:quasiconcave_quasiconvex} guarantees that the RPO problem is randomization-proof. \Halmos
\end{example}

\subsection{Example of necessity of uniqueness assumption in part (b) of Corollary~\ref{corollary:finite_Pcal_pDR}}
\label{counterexample_corollary:finite_Pcal_pDR}

The following example illustrates that the implication in part (b) of Corollary~\ref{corollary:finite_Pcal_pDR} may fail to hold if $\ub^*$ is not the unique solution of $\min_{\ub \in \Ucal} R(\pb^*_{\DR}, \ub)$. Additionally, this example also illustrates that even when $\Pcal$ is finite, the RPO problem can be randomization-proof in the absence of strong duality. 

\begin{example}
\label{example:finitePcal_uniqueness_ustar}

Consider a single product pricing instance, i.e., $I = 1$, which we define as follows. Let $\Pcal = \{ p_1, p_2, p_3 \}$ where $p_1 = 5$, $p_2 = 8$, $p_3  = 9$. Let the demand model $d$ be a linear demand model, so that the uncertain parameter $\ub = (\alpha, \beta)$ and $d(p, \ub) = \alpha - \beta p$. Finally, let $\Ucal = \{ (\alpha_1, \beta_1), (\alpha_2, \beta_2), (\alpha_3, \beta_3) \}$, where $(\alpha_1, \beta_1) = (10, 1)$, $(\alpha_2, \beta_2) = (3,0.1)$, $(\alpha_3, \beta_3) = (3.6, 0.2)$.

We first calculate $\min_{\ub \in \Ucal} R(p, \ub)$ for each $p \in \Pcal$. We have:
\begin{itemize}
\item For $p_1 = 5$: $p_1 (\alpha_2 - \beta_2 p_1) = 12.5 <  p_1 (\alpha_3 - \beta_3 p_1)  = 13 < p_1(\alpha_1 -  \beta_1 p_1) = 25$. Hence, $\min_{\ub \in \Ucal} R(p_1, \ub) = \min\{12.5, 13, 25\} = 12.5$. 
\item For $p_2 = 8$: $p_2 (\alpha_1 - \beta_1 p_2) = p_2 (\alpha_3 - \beta_3 p_2) = 16 < p_2 (\alpha_2 - \beta_2 p_2) = 17.6$. Hence, $\min_{\ub \in \Ucal} R(p_2, \ub) = \min\{16, 16, 17.6\} = 16$, and note also that the minimizing $\ub$ is not unique (the minimum is attained at both $(\alpha_1, \beta_1)$ and $(\alpha_3, \beta_3)$).
\item For $p_3 = 9$: $p_3 (\alpha_1 - \beta_1 p_3) = 9 < p_3 (\alpha_3 - \beta_3 p_3) = 16.2 < p_3 (\alpha_2 - \beta_2 p_3) = 18.9$. Hence, $\min_{\ub \in \Ucal} R(p_3, \ub) = \min \{9, 16.2, 18.9\} = 9$. 
\end{itemize}
From this, we can see that the optimal deterministic robust price is $p^*_{\DR}  = p_2 = 8$ and the optimal deterministic robust objective value is $Z^*_{\DR} =16$. At $p = 8$, we can see that $\arg \min_{\ub \in \Ucal} R(p_2, \ub) = \{ (\alpha_1, \beta_1), (\alpha_3, \beta_3)\}$. %{\color{blue}In other words, the minimizing $\ub^*$ for this $p$ is not unique.}

Let us now consider the RRPO problem. When we write the problem $\max_{\pib \in \Delta_{\Pcal}} \min_{\ub \in \Ucal} \sum_{\pb \in \Pcal} \pi_{\pb} R(\pb, \ub)$ as a linear program, we get the following problem: 
\begin{subequations}
\begin{alignat}{2}
& \underset{\pib, t}{\text{maximize}} & \quad & t \\
& \text{subject to} & & t \leq \pi_{p_1} \cdot p_1 (\alpha_1 - \beta_1 p_1) + \pi_{p_2} \cdot p_2 (\alpha_1 - \beta_1 p_2) + \pi_{p_3} \cdot p_3 (\alpha_1 - \beta_1 p_3), \\
& & & t \leq \pi_{p_1} \cdot p_1 (\alpha_2 - \beta_2 p_1) + \pi_{p_2} \cdot p_2 (\alpha_2 - \beta_2 p_2) + \pi_{p_3} \cdot p_3 (\alpha_2 - \beta_2 p_3), \\
& & & t \leq \pi_{p_1} \cdot p_1 (\alpha_3 - \beta_3 p_1) + \pi_{p_2} \cdot p_2 (\alpha_3 - \beta_3 p_2) + \pi_{p_3} \cdot p_3 (\alpha_3 - \beta_3 p_3), \\
& & & \pi_{p_1} + \pi_{p_2} + \pi_{p_3} = 1, \\
& & & \pi_{p_1}, \pi_{p_2}, \pi_{p_3} \geq 0,
\end{alignat}
\end{subequations}
or equivalently,
\begin{subequations}
\begin{alignat}{2}
& \underset{\pib, t}{\text{maximize}} & \quad & t \\
& \text{subject to} & & t \leq 25 \pi_{p_1} + 16 \pi_{p_2} + 9 \pi_{p_3}  \\
& & & t \leq 12.5 \pi_{p_1}  + 17.6 \pi_{p_2}  + 18.9 \pi_{p_3} , \\
& & & t \leq 13 \pi_{p_1} + 16 \pi_{p_2} + 16.2 \pi_{p_3}, \\
& & & \pi_{p_1} + \pi_{p_2} + \pi_{p_3} = 1, \\
& & & \pi_{p_1}, \pi_{p_2}, \pi_{p_3} \geq 0,
\end{alignat}
\end{subequations}
for which the optimal objective value is $Z^*_{\RR} = 16$, which is the same as $Z^*_{\DR}$. Thus, we can see that if the problem is randomization-proof, but the uniqueness condition on $\min_{\ub \in \Ucal} R(\pb^*_{\DR}, \ub)$ is not met, then $p^*_{\DR}$ is not necessarily a maximizer of the nominal problem at any $\ub^*$. (For $\ub^* = (\alpha_1,\beta_1)$, the $p$ that maximizes $R(\cdot, \ub^*)$ is $p_1$, but for $\ub^* = (\alpha_3, \beta_3)$, the $p$ that maximizes $R(\cdot, \ub^*)$ is $p_3$.)

%Thus, if the uniqueness condition on $\min_{\ub \in \Ucal} R(\pb^*_{\DR}, \ub)$ is relaxed, then it is possible for the problem to be randomization proof. 

Additionally, observe that in this instance, strong duality does not hold. If we compute $\min_{\ub \in \Ucal} \max_{p \in \Pcal} R(p, \ub)$, we get:
\begin{align*}
& \min_{\ub \in \Ucal} \max_{p \in \Pcal} R(p, \ub) \\
& = \min \{ \max_{p \in \Pcal} R(p, \ub_1), \max_{p \in \Pcal} R(p, \ub_2), \max_{p \in \Pcal} R(p, \ub_3) \} \\
& = \min\{ \max\{25, 16, 9\}, \max\{12.5, 17.6, 18.9\}, \max\{13, 16, 16.2\} \} \\
& = \min\{25, 18.9, 16.2\} \\
& = 16.2,
\end{align*}
which is strictly higher than $Z^*_{\RR} = 16$. \Halmos
\end{example}

\subsection{Example to illustrate that $Z^*_{\DR}$, $Z^*_{\RR}$ and dual objective can be distinct when $\Pcal$ is finite}
\label{subsec:examples_benefits_DR_RR_dual_all_different}

In this section, we revisit Example~\ref{example:concave_revenue_fn_set_convexity_unrelaxable} from Section~\ref{subsec:examples_benefits_concave}, which can be viewed as an instance of a finite $\Pcal$ robust price optimization problem.

\begin{example}
Consider a single-product RPO problem, and suppose that $\Pcal = \{p_1, p_2\}$ where $p_1 = 5$, $p_2 = 10$. Consider $\ub = (u_1)$, where $\Ucal = \{ 0, 1 \}$, and let $R(\pb, \ub)$ be defined as
\begin{equation}
R(p,\ub) = u_1 \cdot p(10 - p) + (1 - u_1) \cdot p (4 - 0.2p).
\end{equation}
Note that this is the same instance as Example~\ref{example:concave_revenue_fn_set_convexity_unrelaxable}, with the only difference being that we have re-defined the uncertain parameter $\ub$ to be a binary variable that toggles between the two revenue functions $p(10 - p)$ and $p(4 - 0.2p)$. For this instance, $Z^*_{\DR} = 15$ and $Z^*_{\RR} = 50/3 = 16.6667$, as before, and the dual objective $\min_{\ub \in \Ucal} \max_{\pb \in \Pcal} R(\pb, \ub)$ is 
\begin{align}
\min_{\ub \in \Ucal} \max_{\pb \in \Pcal} R(\pb, \ub) & = \min\{ \max\{25, 0\}, \max\{15, 20\} \} \\
& = \min\{ 25, 20 \} \\
& = 20,
\end{align}
which is strictly larger than $Z^*_{\RR}$. \Halmos
\end{example}

Note that in the previous example, the difference between $\min_{\ub \in \Ucal} \max_{\pb \in \Pcal} R(\pb, \ub)$ and $Z^*_{\RR}$ is driven by $\Ucal$ being a discrete set. If we replaced $\Ucal$ with its convex hull, then $Z^*_{\RR}$ and the dual objective would be identical, which is consistent with Proposition~\ref{proof_proposition:equivalence_SD_randproof_finitePcal} (when $\Ucal$ is convex and $R(\pb, \ub)$ is convex in $\ub$, then $Z^*_{\RR}$ and the dual objective are equal, and strong duality is equivalent to randomization-proofness). 

\begin{example}
Consider the same single-product RPO problem as in the preceding example, but let $\Ucal = [0,1]$. In this case, both $Z^*_{\RR}$ and $Z^*_{\DR}$ are unchanged because the objective of the inner minimization over $\ub \in \Ucal$ is unchanged when we replace $\Ucal = \{0,1\}$ with $\Ucal = [0,1]$. For the dual objective, we get 
\begin{align}
\min_{\ub \in \Ucal} \max_{\pb \in \Pcal} R(\pb, \ub) & = \min_{0 \leq u_1 \leq 1} \max\{ u_1 \cdot 25 + (1-u_1) 15, (1 - u_1) \cdot 20 \} \\
& =  \min_{0 \leq u_1 \leq 1} \{ 10u_1 + 15, 20 - 20u_1 \} \\
& = 50 / 3,
\end{align}
which is the same as $Z^*_{\RR}$, as predicted by Proposition~\ref{proof_proposition:equivalence_SD_randproof_finitePcal}. \Halmos
\end{example}

%%%% 
% END OF EXAMPLE / COUNTEREXAMPLE SECTION FOR SECTION $.
%%%%

%\end{proofvvm}

%%%% 
% START OF MIECP EXPLICIT FORMULATION FOR SECTION 5
%%%%

\section{Mixed-integer exponential cone formulations for Section~\ref{sec:finitePcal_convexUcal} }
\label{sec:miecp_explicit_models}

In this section, we describe how the two mixed-integer convex programs in Section~\ref{sec:finitePcal_convexUcal}, for the semi-log demand model (problem~\eqref{prob:semilog_MIECP} in Section~\ref{subsec:finitePcal_convexUcal_semilog}) and the log-log demand model (problem~\eqref{prob:loglog_MIECP} in Section~\ref{subsec:finitePcal_convexUcal_loglog}) can be reformulated as mixed-integer exponential cone programs. 

We first recall that the standard exponential cone is defined as
\begin{equation}
\Kcal_{\exp} = \{ (r, 0, t) \in \mathbb{R}^3 \mid r \geq 0, t \leq 0\} \cup \{ (r,s,t) \in \mathbb{R}^3 \mid s > 0, r \geq s \exp( t / s) \}.
\end{equation}
We also recall that the function $f(\mu) = \mu \log \mu$ is convex, for $\mu > 0$. To model this function using the exponential cone, let us introduce the epigraph variable $q$. We can then model the epigraph of this function as follows:
\begin{align}
& q \geq \mu \log \mu \\
& \Leftrightarrow \frac{q}{\mu} \geq \log \mu \\
& \Leftrightarrow e^{q / \mu} \geq \mu \\
& \Leftrightarrow 1 \geq \mu e^{- q / \mu} \\
& \Leftrightarrow (1, \mu, -q) \in \Kcal_{\exp}.
\end{align}

Therefore, the semi-log demand formulation~\eqref{prob:semilog_MIECP} can be written using the auxiliary epigraph variables $q_1, \dots, q_I$ as
\begin{subequations}
\begin{alignat}{2}
& \underset{\mub, \qb, \wb, \xb}{\text{maximize}} & \quad & \sum_{i=1}^I \mu_i \alpha_i + \sum_{i=1}^I \sum_{t \in \Pcal_i} \log t \cdot w_{i,i,t} - \sum_{i=1}^I \beta_i \cdot \sum_{t \in \Pcal_i} t \cdot w_{i,i,t} + \sum_{i=1}^I \sum_{j\neq i} \gamma_{i,j} \cdot (\sum_{t \in \Pcal_j} t \cdot w_{i,j,t}) - \sum_{i=1}^I q_i \\
& \text{subject to} & & (1, \mu_i, -q_i) \in \Kcal_{\exp}, \quad \forall \ i \in [I], \\
& & & \text{constraints~\eqref{prob:semilog_MIECP_w_sums_to_mu} -- \eqref{prob:semilog_MIECP_mu_nonnegative}}. 
\end{alignat}
\end{subequations}

In a similar fashion, the log-log demand formulation~\eqref{prob:loglog_MIECP} can be written as 
\begin{subequations}
\begin{alignat}{2}
& \underset{\mub, \qb, \wb, \xb}{\text{maximize}} & \quad & \sum_{i=1}^I \mu_i \alpha_i + \sum_{i=1}^I \sum_{t \in \Pcal_i} \log t \cdot w_{i,i,t} - \sum_{i=1}^I \beta_i \cdot \sum_{t \in \Pcal_i} \log t \cdot w_{i,i,t} + \sum_{i=1}^I \sum_{j\neq i} \gamma_{i,j} \sum_{t \in \Pcal_j} \log t \cdot w_{i,j,t} - \sum_{i=1}^I q_i \\ 
& \text{subject to} & & (1, \mu_i, -q_i) \in \Kcal_{\exp}, \quad \forall \ i \in [I], \\
& & & \text{constraints~\eqref{prob:semilog_MIECP_w_sums_to_mu} -- \eqref{prob:semilog_MIECP_mu_nonnegative}}. 
\end{alignat}
\end{subequations}

%%%%
% END OF MIECP EXPLICIT FORMULATIONS FOR SECTION 5.
%%%%

\section{Deterministic robust price optimization for finite $\Pcal$, convex $\Ucal$ under the semi-log and log-log demand models}
\label{sec:convexUcal_DRPO}

In this section, we describe how to formulate the DRPO problem as a mixed-integer exponential cone program for the semi-log and log-log demand models. In both cases, we assume that $\Ucal$ is a convex uncertainty set, and that Assumption~\ref{assumption:Pcal_Cartesian_product} on the structure of $\Pcal$ holds. 

\subsection{Semi-log model}
\label{subsec:convexUcal_DRPO_semilog}

For the semi-log demand model, we can write the DRPO problem as
\begin{align}
& \max_{\pb \in \Pcal} \min_{\ub \in \Ucal} R(\pb, \ub) \nonumber \\
& = \max_{\pb \in \Pcal} \min_{\ub \in \Ucal} \sum_{i=1}^I p_i \cdot e^{\alpha_i - \beta_i p_i + \sum_{j \neq i} \gamma_{i,j} p_j }. \label{prob:DRPO_semilog_abstract}
\end{align}
To accomplish our reformulation, we will make use of the fact that the optimal solution set of the DRPO problem is unchanged upon log-transformation, that is, 
\begin{align*}
\arg \max_{\pb \in \Pcal} \min_{\ub \in \Ucal} R(\pb, \ub) = \arg \max_{\pb \in \Pcal} \min_{\ub \in \Ucal} \log R(\pb, \ub).
\end{align*}
Thus, instead of problem~\eqref{prob:DRPO_semilog_abstract}, we can focus on the following problem:
\begin{align*}
& \max_{\pb \in \Pcal} \min_{\ub \in \Ucal} \log \left( \sum_{i=1}^I p_i \cdot e^{\alpha_i - \beta_i p_i + \sum_{j \neq i} \gamma_{i,j} p_j } \right) \\
& = \max_{\pb \in \Pcal} \min_{\ub \in \Ucal} \log \left( \sum_{i=1}^I e^{\log p_i + \alpha_i - \beta_i p_i + \sum_{j \neq i} \gamma_{i,j} p_j } \right)
\end{align*}
Here, we can again use the log-sum-exp biconjugate technique to reformulate the objective function in the following way:
\begin{align*}
& \log \left( \sum_{i=1}^I e^{\log p_i + \alpha_i - \beta_i p_i + \sum_{j \neq i} \gamma_{i,j} p_j } \right) \\
& = \max_{\mub \in \Delta_{[I]}} \left\{ \sum_{i=1}^I \mu_i \cdot (\log p_i + \alpha_i - \beta_i p_i + \sum_{j \neq i} \gamma_{i,j} p_j) - \sum_{i=1}^I \mu_i \log \mu_i \right\}.
\end{align*}
Thus, the overall problem becomes the following max-min-max problem:
\begin{align*}
\max_{\pb \in \Pcal} \min_{\ub \in \Ucal} \max_{\mub \in \Delta_{[I]}} \left\{ \sum_{i=1}^I \mu_i \cdot (\log p_i + \alpha_i - \beta_i p_i + \sum_{j \neq i} \gamma_{i,j} p_j) - \sum_{i=1}^I \mu_i \log \mu_i \right\}.
\end{align*}
Here, we observe that the objective function is linear in $\ub = (\alphab, \betab, \gammab)$, and is concave in $\mub$; additionally, the feasible region of $\ub$ is assumed to be convex, and the feasible region of $\mub$ is convex and compact (being just the $(|I|-1)$-dimensional unit simplex). Therefore, we can use Sion's minimax theorem to interchange the minimization over $\ub$ and the maximization over $\mub$, which gives us 
\begin{align*}
& \max_{\pb \in \Pcal} \min_{\ub \in \Ucal} \max_{\mub \in \Delta_{[I]}} \left\{ \sum_{i=1}^I \mu_i \cdot (\log p_i + \alpha_i - \beta_i p_i + \sum_{j \neq i} \gamma_{i,j} p_j) - \sum_{i=1}^I \mu_i \log \mu_i \right\} \\
& = \max_{\pb \in \Pcal} \max_{\mub \in \Delta_{[I]}}  \min_{\ub \in \Ucal} \left\{ \sum_{i=1}^I \mu_i \cdot (\log p_i + \alpha_i - \beta_i p_i + \sum_{j \neq i} \gamma_{i,j} p_j) - \sum_{i=1}^I \mu_i \log \mu_i \right\} \\
& = \max_{\pb \in \Pcal, \mub \in \Delta_{[I]}} \min_{\ub \in \Ucal} \left\{ \sum_{i=1}^I \mu_i \cdot (\log p_i + \alpha_i - \beta_i p_i + \sum_{j \neq i} \gamma_{i,j} p_j) - \sum_{i=1}^I \mu_i \log \mu_i \right\}
\end{align*}
Under Assumption~\ref{assumption:Pcal_Cartesian_product}, this final problem can then be reformulated as robust mixed-integer exponential cone program, just as in Section~\ref{subsec:finitePcal_convexUcal_semilog}. We introduce the same binary decision variable $x_{i,t}$ which is 1 if product $i$ is offered at price $t \in \Pcal_i$, and 0 otherwise, and we use $w_{i,j,t}$ to denote the linearization of $\mu_i \cdot x_{j,t}$ for $i, j \in [I]$, $t \in \Pcal_j$. This gives rise to the following program:
\begin{subequations}
\begin{alignat}{2}
& \underset{\mub, \wb, \xb}{\text{maximize}} & \quad & \min_{\ub \in \Ucal} \left\{ \sum_{i=1}^I \mu_i \alpha_i + \sum_{i=1}^I  \left( \sum_{t \in \Pcal_i} \log t \cdot w_{i,i,t} - \beta_i \cdot \sum_{t \in \Pcal_i} t \cdot w_{i,i,t} + \sum_{j\neq i} \gamma_{i,j} \sum_{t \in \Pcal_j} t \cdot w_{i,j,t} \right) - \sum_{i=1}^I \mu_i \log \mu_i  \right\} \\ 
& \text{subject to} & &  \sum_{t \in \Pcal_j} w_{i,j,t} = \mu_i, \quad \forall \ i \in [I], \ j \in [I], \label{prob:DRPO_semilog_MIECP_w_sums_to_mu} \\
& & & \sum_{i=1}^I w_{i,j,t} = x_{j,t}, \quad \forall \ j \in [I], \ t \in \Pcal_j, \label{prob:DRPO_semilog_MIECP_w_sums_to_x} \\
& & & \sum_{i=1}^I \mu_i = 1, \\
& & & \sum_{t \in \Pcal_i} x_{i,t} = 1, \quad \forall \ i \in [I], \\
& & & w_{i,j,t} \geq 0, \quad \forall \ i \in [I],\ j \in [I], t \in \Pcal_j, \label{prob:DRPO_semilog_MIECP_w_nonnegative} \\
& & & x_{i,t} \in \{0,1\}, \quad \forall \ i \in [I], \ t \in \Pcal_i, \\
& & & \mu_i \geq 0, \quad \forall \ i \in [I].
\end{alignat}
\label{prob:DRPO_semilog_MIECP}%
\end{subequations}
Note that the feasible region of this problem is identical to that of problem~\eqref{prob:semilog_MIECP}, which appeared in our discussion of the separation problem for the RRPO problem when $\Ucal$ is convex and $\Pcal$ is finite. The difference here is that the objective is now a robust objective; it is the worst-case value of the objective of problem~\eqref{prob:semilog_MIECP}, taken over the convex uncertainty set $\Ucal$. Depending on the structure of $\Ucal$, the overall problem can remain in the mixed-integer convex program problem class. For example, if $\Ucal$ is a polyhedron, then one can use LP duality to reformulate the robust problem exactly by introducing additional variables and constraints, as is normally done in robust optimization \citep{bertsimas2004price,ben2000robust,bertsimas2011theory}. Similarly, if $\Ucal$ is a second-order cone representable set, then one can again use conic duality to reformulate the problem. Alternatively, one can also consider solving the problem using a cutting plane method/delayed constraint generation approach, whereby one reformulates the program in epigraph form and then solves the inner minimization over $\ub$ to identify new cuts to add \citep{bertsimas2016reformulation}.

\subsection{Log-log model} 
\label{subsec:convexUcal_DRPO_loglog}
For the log-log demand model, we can write the DRPO problem as
\begin{align*}
& \max_{\pb \in \Pcal} \min_{\ub \in \Ucal} R(\pb, \ub) \\
& = \max_{\pb \in \Pcal} \min_{\ub \in \Ucal} \sum_{i=1}^I p_i \cdot e^{\alpha_i - \beta_i \log p_i + \sum_{j \neq i} \gamma_{i,j} \log p_j } \\
& = \max_{\pb \in \Pcal} \min_{\ub \in \Ucal} \sum_{i=1}^I e^{\alpha_i + (1 - \beta_i) \log p_i + \sum_{j \neq i} \gamma_{i,j} \log p_j } 
\end{align*}
Again, as with the semi-log model, solving the above problem is equivalent to solving the same problem with a log-transformed objective. Taking this log-transformed problem as our starting point, replacing the log-sum-exp function with its biconjugate and applying Sion's minimax theorem gives us:
\begin{align*}
& \max_{\pb \in \Pcal} \min_{\ub \in \Ucal} \log \left( \sum_{i=1}^I e^{\alpha_i + (1 - \beta_i) \log p_i + \sum_{j \neq i} \gamma_{i,j} \log p_j } \right) \\
& = \max_{\pb \in \Pcal} \min_{\ub \in \Ucal} \max_{\mub \in \Delta_{[I]}} \left\{ \sum_{i=1}^I \left[ \alpha_i \mu_i + (1 - \beta_i) \mu_i \cdot \log p_i + \sum_{j \neq i} \gamma_{i,j} \mu_i \cdot \log p_j \right] - \sum_{i=1}^I \mu_i \log \mu_i  \right\} \\
& = \max_{\pb \in \Pcal, \mub \in \Delta_{[I]}} \min_{\ub \in \Ucal} \left\{ \sum_{i=1}^I \left[ \alpha_i \mu_i + (1 - \beta_i) \mu_i \cdot \log p_i + \sum_{j \neq i} \gamma_{i,j} \mu_i \cdot \log p_j \right] - \sum_{i=1}^I \mu_i \log \mu_i \right\}.
\end{align*}
Under Assumption~\ref{assumption:Pcal_Cartesian_product}, this last problem can be re-written as the following robust version of problem~\eqref{prob:loglog_MIECP}, with the decision variables defined identically:
\begin{subequations}
\begin{alignat}{2}
& \underset{\mub, \wb, \xb}{\text{maximize}} & & \min_{\ub \in \Ucal} \left\{  \sum_{i=1}^I \mu_i \alpha_i + \sum_{i=1}^I \left(\sum_{t \in \Pcal_i} \log t \cdot w_{i,i,t} - \beta_i \cdot \sum_{t \in \Pcal_i} \log t \cdot w_{i,i,t} + \sum_{j\neq i} \gamma_{i,j} \sum_{t \in \Pcal_j} \log t \cdot w_{i,j,t} \right) - \sum_{i=1}^I \mu_i \log \mu_i  \right\} \\
& \text{subject to} & \quad & \text{constraints~\eqref{prob:semilog_MIECP_w_sums_to_mu} -- \eqref{prob:semilog_MIECP_mu_nonnegative}.}
%\text{constraints~\eqref{prob:loglog_MIECP_w_sums_to_mu} - \eqref{prob:loglog_MIECP_mu_nonnegative}} 
\end{alignat}
\label{prob:DRPO_loglog_MIECP}%
\end{subequations}
Again, this problem has exactly the same feasible region as the log-log separation problem~\eqref{prob:loglog_MIECP} and the semi-log separation problem~\eqref{prob:semilog_MIECP}. Additionally, just as with the deterministic robust problem~\eqref{prob:DRPO_semilog_MIECP} for the semi-log model, this problem can be further reformulated by exploiting the structure of $\Ucal$, or otherwise one can design a cutting plane method that generates violated uncertain parameter vectors $\ub \in \Ucal$ on the fly.

%%%%
% START OF CONSTRAINED DISTRIBUTION SECTION
%%%%

\section{Incorporating constraints on price distributions}
\label{sec:constrained}

In this section, we present an extension of our framework to consider the possibility of constrained distributions over price vectors in the finite $\Pcal$, convex $\Ucal$ case. In Section~\ref{subsec:constrained_model}, we show how a particular type of constrained RRPO problem can be reformulated as a large-scale convex program. In Section~\ref{subsec:constrained_results}, we present results of a small numerical experiment to show that incorporating constraints on the price distribution produced by the RRPO approach does not result in a large deterioration in performance compared to the DRPO approach.

\subsection{Randomized pricing with constraints}
\label{subsec:constrained_model}

Suppose that we consider the following more general instance of the RRPO problem:
\begin{equation}
\textsc{RRPO-C}: \quad \max_{\pib \in \tilde{\Delta}_{\Pcal} } \min_{\ub \in \Ucal} \sum_{\pb \in \Pcal} \pi_{\pb} R(\pb, \ub)
\end{equation}
where $\tilde{\Delta}_{\Pcal}$ is the following subset of the probability simplex $\Delta_{\Pcal}$: 
%\begin{equation}
%\tilde{\Delta}_{\Pcal} = \left\{ \pib \in \Delta_{\Pcal} \mid \sum_{\pb \in \Pcal} g_k(\pb) \pi_{\pb} = c_k, \ \forall k \in [K], \right\}.
%\end{equation}
 \begin{equation}
 \tilde{\Delta}_{\Pcal} = \left\{ \pib \in \Delta_{\Pcal} \ \vline \ \begin{array}{l} \sum_{\pb \in \Pcal} g_k(\pb) \pi_{\pb} = c_k, \ \forall k \in [K_1], \\
 \sum_{\pb \in \Pcal} h_k(\pb) \pi_{\pb} \leq d_k, \ \forall k \in [K_2], \end{array} \right\},
 \end{equation}
 where $K_1$ and $K_2$ are nonnegative integers.  In words, $\tilde{\Delta}_{\Pcal}$ is the set of all probability distributions that satisfy a set of \emph{moment conditions} with respect to a set of moment functions $\{ g_k(\cdot) \}_{k=1}^{K_1}$, $\{ h_k(\cdot) \}_{k=1}^{K_2}$; in particular, the expected value of $g_k(\tilde{\pb})$, where $\tilde{\pb}$ is the random price vector that follows the distribution $\pib$, must be equal to $c_k$, for $K_1$ different moment functions. Similarly, the expected value of $h_k(\tilde{\pb})$ must be less than or equal to $d_k$ for $K_2$ different moment functions. 

Let us now see how to reformulate this as a tractable problem. Let us assume that $\Ucal$ is convex and that $R(\pb, \ub)$ is convex in $\ub$. In this case, Sion's minimax theorem again lets us assert
\begin{align}
& \max_{\pib \in \tilde{\Delta}_{\Pcal} } \min_{\ub \in \Ucal} \sum_{\pb \in \Pcal} \pi_{\pb} R(\pb, \ub)\\
& = \min_{\ub \in \Ucal} \max_{\pib \in \tilde{\Delta}_{\Pcal} } \sum_{\pb \in \Pcal} \pi_{\pb} R(\pb, \ub), \label{eq:min_ub_max_pib_tildeDelta}
\end{align}
where we are also using the fact that $\tilde{\Delta}_{\Pcal}$ is convex, as it is a polyhedron, and is compact, as it is a closed subset of the compact set $\Delta_{\Pcal}$. 

Now, let us introduce a $K_1$-dimensional Lagrange multiplier vector $\lambdab$ and a $K_2$-dimensional nonpositive Lagrange multiplier vector $\rhob$. Observe that the inner maximization over $\pib$ can be re-written equivalently as 
\begin{align}
& \max_{\pib \in \tilde{\Delta}_{\Pcal} } \sum_{\pb \in \Pcal} \pi_{\pb} R(\pb, \ub)  \\
& = \max_{\pib \in \Delta_{\Pcal} } \min_{\lambdab \in \Rbb^{K_1}, \rhob \leq \zerob} \left\{ \sum_{\pb \in \Pcal} \pi_{\pb} R(\pb, \ub) + \sum_{k=1}^{K_1} \lambda_k \cdot \left(\sum_{\pb \in \Pcal} g_k(\pb) \pi_{\pb} - c_k \right)  + \sum_{k=1}^{K_2} \rho_k \cdot \left( \sum_{\pb \in \Pcal} h_k(\pb) \pi_{\pb} - d_k \right) \right\}  \\
& = \max_{\pib \in \Delta_{\Pcal} } \min_{\lambdab \in \Rbb^{K_1}, \rhob \leq \zerob} \left\{ \sum_{\pb \in \Pcal} \pi_{\pb} \left( R(\pb, \ub) + \sum_{k=1}^{K_1} \lambda_k g_k(\pb) + \sum_{k=1}^{K_2} \rho_k h_k(\pb) \right) - \sum_{k=1}^{K_1} \lambda_k c_k - \sum_{k=1}^{K_2} \rho_k d_k  \right\}  \\
%& = \max_{\pib \in \Delta_{\Pcal} } \min_{\lambdab \in \Rbb^K} \left\{ \sum_{\pb \in \Pcal} \pi_{\pb} \left[ R(\pb, \ub)  \ri+ \sum_{k=1}^K \lambda_k (   \right\}  \\
& = \min_{\lambdab \in \Rbb^{K_1}, \rhob \leq \zerob} \max_{\pib \in \Delta_{\Pcal} } \left\{ \sum_{\pb \in \Pcal} \pi_{\pb} \left( R(\pb, \ub) + \sum_{k=1}^{K_1} \lambda_k g_k(\pb) + \sum_{k=1}^{K_2} \rho_k h_k(\pb) \right) - \sum_{k=1}^{K_1} \lambda_k c_k - \sum_{k=1}^{K_2} \rho_k d_k  \right\}  \\
& = \min_{\lambdab \in \Rbb^{K_1}, \rhob \leq \zerob} \left\{ \max_{\pib \in \Delta_{\Pcal} } \left\{ \sum_{\pb \in \Pcal} \pi_{\pb} \left( R(\pb, \ub) + \sum_{k=1}^{K_1} \lambda_k g_k(\pb) + \sum_{k=1}^{K_2} \rho_k h_k(\pb) \right) \right\} - \sum_{k=1}^{K_1} \lambda_k c_k  - \sum_{k=1}^{K_2} \rho_k d_k \right\}  \\
& = \min_{\lambdab \in \Rbb^{K_1}, \rhob \leq \zerob} \left\{ \max_{\pb \in \Pcal } \left( R(\pb, \ub) + \sum_{k=1}^{K_1} \lambda_k g_k(\pb) + \sum_{k=1}^{K_2} \rho_k h_k(\pb) \right) - \sum_{k=1}^{K_1} \lambda_k c_k  - \sum_{k=1}^{K_2} \rho_k d_k \right\} 
\label{eq:min_lambda_max_pb}
\end{align}
where the interchange of maximization over $\pib \in \Delta_{\Pcal}$ and minimization over $\lambdab$ and $\rhob$ follows by linear programming strong duality, and the final equality follows because the optimal solution of the inner maximization over $\pib \in \Delta_{\Pcal}$ occurs when all of the probability mass is placed at a single $\pb \in \Pcal$. Combining~\eqref{eq:min_lambda_max_pb} with \eqref{eq:min_ub_max_pib_tildeDelta}, we obtain 
\begin{align}
& \min_{\ub \in \Ucal} \max_{\pib \in \tilde{\Delta}_{\Pcal} } \sum_{\pb \in \Pcal} \pi_{\pb} R(\pb, \ub)\\
& = \min_{\ub \in \Ucal} \min_{\lambdab \in \Rbb^{K_1}, \rhob \leq \zerob} \left\{ \max_{ \pb \in \Pcal } \left( R(\pb, \ub) + \sum_{k=1}^{K_1} \lambda_k g_k(\pb) + \sum_{k=1}^{K_2} \rho_k h_k(\pb) \right)  - \sum_{k=1}^{K_1} \lambda_k c_k - \sum_{k=1}^{K_2} \rho_k d_k \right\} \\
& = \min_{\ub \in \Ucal, \lambdab \in \Rbb^{K_1}, \rhob \leq \zerob} \left\{ \max_{ \pb \in \Pcal } \left( R(\pb, \ub) + \sum_{k=1}^{K_1} \lambda_k g_k(\pb) + \sum_{k=1}^{K_2} \rho_k h_k(\pb) \right)  - \sum_{k=1}^{K_1} \lambda_k c_k - \sum_{k=1}^{K_2} \rho_k d_k  \right\} 
\end{align}
This final problem can be re-written as the following convex optimization problem, where we introduce an epigraph variable $t$ to model the innermost maximum over $\pb$:
\begin{subequations}
\begin{alignat}{2}
& \underset{\lambdab, \rhob, t, \ub}{\text{minimize}} & \quad & t - \sum_{k=1}^{K_1} \lambda_k c_k - \sum_{k=1}^{K_2} \rho_k d_k \\
& \text{subject to} &  & t \geq R(\pb, \ub) + \sum_{k=1}^{K_1} \lambda_k g_k(\pb) + \sum_{k=1}^{K_2} \rho_k h_k(\pb), \quad \forall \ \pb \in \Pcal, \\ 
& & & \ub \in \Ucal, \\
& & & \lambdab \in \Rbb^{K_1}, \\
& & & \rhob \leq \zerob, \\
& & & t \in \Rbb.
\end{alignat}
\end{subequations}
Like the unconstrained problem~\eqref{prob:finiteP_convexU_epigraph}, this problem can again potentially be solved via constraint generation. The challenge in solving this problem now becomes how one solves the modified separation problem
\begin{equation}
\max_{\pb \in \Pcal} \left\{ R(\pb, \ub) + \sum_{k=1}^{K_1} \lambda_k g_k(\pb)  + \sum_{k=1}^{K_2} \rho_k h_k(\pb) \right\}.
\end{equation}
This separation problem is more challenging than the separation problem in the unconstrained distribution setting, because of the presence of the $\sum_{k=1}^{K_1} \lambda_k g_k(\pb)$ and $\sum_{k=1}^{K_2} \rho_k h_k(\pb)$ terms in the objective function. In particular, because of these terms and the fact that the $\lambda_k$ and $\rho_k$ variables can be negative, one cannot straightforwardly apply the biconjugate reformulation technique for the case of log-log and semi-log demand model. However, one can solve the problem heuristically, using the same random improvement approach that we use in Section~\ref{subsec:results_discreteUcal_loglog_semilog}. We adopt this approach in our numerical experiment in Section~\ref{subsec:constrained_results}. 

With regard to the moment functions $g_k(\cdot)$ and $h_k(\cdot)$, one natural choice of moment functions correspond to the first and second moments of each product's random price. In particular, letting $\tilde{\pb}$ be the random price vector, one may wish to enforce that the distribution $\pib$ over $\Pcal$ is such that 
\begin{align}
\Ebb_{\pib} [ \tilde{p}_i ] & = \mu_i, \quad \forall \ i \in [I], \label{eq:moment_first}\\
\Ebb_{\pib} [ \tilde{p}^2_i ] & \leq \mu_i^2 + \sigma_i^2, \quad \forall \ i \in [I], \label{eq:moment_second}
\end{align}
where $\mu_i$ and $\sigma_i$ are user-specified values for the desired mean and standard deviation of each product's random price. This corresponds to specifying $2I$ moment functions, where $g_i(\pb) = p_i$ for $1 \leq i \leq I$, and $h_{i}(\pb) = p^2_i$ for $1 \leq i \leq I$. The resulting instance of the RRPO problem is one where we seek a distribution $\pib$ that maximizes the worst-case revenue, with the additional restriction that the mean of each product's random price be equal to a predefined value $\mu_i$, and the standard deviation is at most $\sigma_i$. By appropriately setting $\mu_i$ and $\sigma_i$, one can obtain randomized pricing schemes with controlled behavior, e.g., where the average price of a product is fixed to a certain value and to ensure that the variability around that value is not too large. As another example, one could set $\mu_i$ and $\sigma_i$ so that the resulting price distribution $\pib$ matches historical price data for a collection of products, i.e., on average, $\pib$ results in the same average prices and the same variability in those prices as in the empirical distribution of $\pb$ in a historical data set. 

We close this section by noting that moment constraints have a rich history in the distributionally robust optimization (DRO) literature; see, for example, \cite{zymler2013distributionally}, \cite{rujeerapaiboon2016robust} and the references therein. We remark here that our use of moment constraints differs from how it is typically used in DRO. In DRO, one seeks to solve a stochastic program where the distribution of a parameter is uncertain, and one constructs an ambiguity set of possible distributions that are consistent with known information and data through moment constraints (e.g., the first and second moments of the uncertain distribution should be equal to certain known values); thus, moment constraints serve to limit how conservative the decision maker is. This contrasts with our use of moment constraints, where the moment constraints are on the firm's decision (the distribution of the price vector $\tilde{\pb}$), and their purpose is to facilitate the implementation of the randomized pricing strategy.  

\subsection{Numerical experiment with \orangeJuice data set}
\label{subsec:constrained_results}

In this section, we demonstrate the value of constrained price distributions using the \orangeJuice problem instances studied in Section~\ref{subsec:results_orangejuice}. We consider the same set of robust price optimization problems as in Section~\ref{subsec:results_orangejuice}, with the modification that $\Pcal_i$ includes $\hat{\mu}_i$, which is the target mean for product $i$; the reason for this will be discussed shortly.

For each value of $\theta$, we solve the constrained randomized problem, with the moment constraints~\eqref{eq:moment_first} and \eqref{eq:moment_second}. We specifically impose that 
\begin{align}
\Ebb_{\pib} [ \tilde{p}_i ] & = \hat{\mu}_i, \quad \forall \ i \in [I], \\
\Ebb_{\pib} [ \tilde{p}^2_i ] & \leq \hat{\mu}^2_i + \alpha^2 \hat{\sigma}_i^2, \quad \forall \ i \in [I],
\end{align}
where $\hat{\mu}_i$ and $\hat{\sigma}_i$ are the mean and standard deviation, respectively, of the price of product $i$, calculated with respect to the optimal unconstrained randomized robust price distribution $\pib'$. In other words, we seek a distribution $\pib$ that has the same marginal means as $\pib'$, but has less variance in each product's price than $\pib'$. We test values of $\alpha \in \{0.5, 0.6, 0.7, 0.8, 0.9\}$. Note that when $\alpha = 1$, the constrained RRPO problem is equivalent to the unconstrained RRPO problem. Additionally, we note that for any $\alpha \geq 0$, the set of distributions $\tilde{\Delta}_{\Pcal}$ defined according to the moment constraints is non-empty, because $\Pcal_i$ includes $\hat{\mu}_i$ (as noted earlier).

%where both are calculated over all stores and weeks in the data set; and $\alpha \in (0,1)$ is a parameter that controls the variability. Note that under this formulation, the mean price of product $i$ is constrained to be $\hat{\mu}_i$, while the standard deviation of the price of product $i$ is constrained to be $\alpha \hat{\sigma}_i$. By varying $\alpha$, we consider price distributions with different amounts of dispersion. We consider values of $\alpha \in \{0.1, 0.2, 0.5, 1.0\}$. Note that because $\Pcal_i$ includes $\hat{\mu}_i$, and because $\Pcal_i$ includes values smaller and larger than $\hat{\mu}_i$, the constrained RRPO problem is always feasible. 

Tables~\ref{table:constrained_loglog} and \ref{table:constrained_semilog} display the results for the log-log and semi-log demand models. From these tables, we can see that in general, when the uncertainty set is small, constraining the distribution using the moment constraints described above results in weaker randomized pricing schemes than the unconstrained approach, and the resulting distributions do not improve on the deterministic robust price solution. For example, with the log-log demand model and $\theta = 0.1$, all of the constrained distributions with $\alpha \leq 0.9$ have a negative $\RI$ metric, indicating that their worst-case revenue is worse than the deterministic robust price solution. However, as the uncertainty set gets larger, there are values of $\alpha < 1$ for which the constrained RRPO solution yields a positive $\RI$. For example, with the log-log demand model and $\theta = 0.5$, the constrained RRPO solution with $\alpha$ set to 0.8 or 0.9 outperforms the deterministic robust solution. As $\theta$ gets larger, it is possible to find distributions with less variability (i.e., with smaller $\alpha$) that outperform the deterministic robust solution. For example, with the log-log demand model and $\theta = 1.0$, the constrained RRPO solution with $\alpha = 0.6$ continues to give a positive $\RI$. These results suggest that by considering a constrained form of the RRPO problem, it is possible to trade-off some of the worst-case performance of the unconstrained RRPO solution to obtain price distributions that are more structured (in the examples shown above, ones with less variability than the unconstrained RRPO solution). Of course, the moment constraints that we have considered are only one way in which one can constrain the price distribution. An interesting direction for future research is to consider other types of constraints that could lead to structured price distributions that continue to provide a benefit over the deterministic robust scheme. %An interesting follow-on question is whether it is possible to allow the mean $\hat{\mu}_i$ to vary, as it might be the case that there are less variable price distributions that are centered elsewhere (other than $(\hat\mu_1,\dots,\hat\mu_I)$) within the convex hull of $\Pcal$ that continue to outperform the deterministic robust solution. From a methodological standpoint, it 

\begin{table}
\centering
\begin{tabular}{lllr} \toprule
$\theta$ & Method $m$ & $Z^*_{m}$ & RI (\%) \\ \midrule
0.10 & RRPO-C, $\alpha = 0.5$ & 237386.69 & -58.05 \\ 
  0.10 & RRPO-C, $\alpha = 0.6$ & 262704.91 & -53.57 \\ 
  0.10 & RRPO-C, $\alpha = 0.7$ & 285955.77 & -49.47 \\ 
  0.10 & RRPO-C, $\alpha = 0.8$ & 313211.80 & -44.65 \\ 
  0.10 & RRPO-C, $\alpha = 0.9$ & 415820.50 & -26.52 \\ 
  0.10 & RRPO & 722647.22 & 27.71 \\ \midrule
  0.50 & RRPO-C, $\alpha = 0.5$ & 162117.27 & -30.54 \\ 
  0.50 & RRPO-C, $\alpha = 0.6$ & 190673.79 & -18.30 \\ 
  0.50 & RRPO-C, $\alpha = 0.7$ & 224890.46 & -3.64 \\ 
  0.50 & RRPO-C, $\alpha = 0.8$ & 251484.85 & 7.75 \\ 
  0.50 & RRPO-C, $\alpha = 0.9$ & 286534.25 & 22.77 \\ 
  0.50 & RRPO & 342614.34 & 46.80 \\ \midrule
  0.80 & RRPO-C, $\alpha = 0.5$ & 131098.05 & -19.21 \\ 
  0.80 & RRPO-C, $\alpha = 0.6$ & 157189.92 & -3.13 \\ 
  0.80 & RRPO-C, $\alpha = 0.7$ & 175564.41 & 8.19 \\ 
  0.80 & RRPO-C, $\alpha = 0.8$ & 200577.13 & 23.60 \\ 
  0.80 & RRPO-C, $\alpha = 0.9$ & 228397.76 & 40.75 \\ 
  0.80 & RRPO & 260049.66 & 60.25 \\ \midrule
  1.00 & RRPO-C, $\alpha = 0.5$ & 111909.24 & -12.72 \\ 
  1.00 & RRPO-C, $\alpha = 0.6$ & 131534.71 & 2.58 \\ 
  1.00 & RRPO-C, $\alpha = 0.7$ & 149969.36 & 16.96 \\ 
  1.00 & RRPO-C, $\alpha = 0.8$ & 172717.05 & 34.70 \\ 
  1.00 & RRPO-C, $\alpha = 0.9$ & 191272.60 & 49.17 \\  
  1.00 & RRPO & 217580.86 & 69.69 \\ \midrule
 1.50 & RRPO-C, $\alpha = 0.5$ & 78581.62 & 3.54 \\ 
  1.50 & RRPO-C, $\alpha = 0.6$ & 92204.17 & 21.48 \\ 
  1.50 & RRPO-C, $\alpha = 0.7$ & 103533.52 & 36.41 \\ 
  1.50 & RRPO-C, $\alpha = 0.8$ & 116139.87 & 53.02 \\ 
  1.50 & RRPO-C, $\alpha = 0.9$ & 127482.81 & 67.97 \\ 
%  1.50 & RRPO-C, $\alpha = 1$ & 141926.68 & 87.00 \\ 
  1.50 & RRPO & 142307.66 & 87.50 \\ \midrule
  2.00 & RRPO-C, $\alpha = 0.5$ & 53377.06 & 8.23 \\ 
  2.00 & RRPO-C, $\alpha = 0.6$ & 61205.36 & 24.10 \\ 
  2.00 & RRPO-C, $\alpha = 0.7$ & 70425.84 & 42.80 \\ 
  2.00 & RRPO-C, $\alpha = 0.8$ & 78122.02 & 58.40 \\ 
  2.00 & RRPO-C, $\alpha = 0.9$ & 87114.25 & 76.63 \\ 
%  2.00 & RRPO-C, $\alpha = 1$ & 94625.20 & 91.86 \\ 
  2.00 & RRPO & 94847.37 & 92.31 \\ \bottomrule
\end{tabular}
\caption{Results of comparison between RRPO, constrained RRPO and DRPO on \orangeJuice instances with log-log demand model. \label{table:constrained_loglog}}
\end{table}

\begin{table}
\centering
\begin{tabular}{lllr} \toprule
$\theta$ & Method $m$ & $Z^*_{m}$ & RI (\%) \\ \midrule
0.10 & RRPO-C, $\alpha = 0.5$ & 153375.29 & -47.20 \\ 
  0.10 & RRPO-C, $\alpha = 0.6$ & 147513.70 & -49.22 \\ 
  0.10 & RRPO-C, $\alpha = 0.7$ & 173292.78 & -40.34 \\ 
  0.10 & RRPO-C, $\alpha = 0.8$ & 171790.80 & -40.86 \\ 
  0.10 & RRPO-C, $\alpha = 0.9$ & 205353.74 & -29.30 \\ 
%  0.10 & RRPO-C, $\alpha = 1$ & 198023.35 & -31.83 \\ 
  0.10 & RRPO & 342357.06 & 17.86 \\ \midrule
  0.50 & RRPO-C, $\alpha = 0.5$ & 117253.37 & -20.64 \\ 
  0.50 & RRPO-C, $\alpha = 0.6$ & 131225.11 & -11.18 \\ 
  0.50 & RRPO-C, $\alpha = 0.7$ & 133415.93 & -9.70 \\ 
  0.50 & RRPO-C, $\alpha = 0.8$ & 155148.87 & 5.01 \\ 
  0.50 & RRPO-C, $\alpha = 0.9$ & 171173.27 & 15.85 \\ 
%  0.50 & RRPO-C, $\alpha = 1$ & 191568.59 & 29.66 \\ 
  0.50 & RRPO & 197517.06 & 33.68 \\ \midrule
  0.80 & RRPO-C, $\alpha = 0.5$ & 87201.93 & -17.53 \\ 
  0.80 & RRPO-C, $\alpha = 0.6$ & 98469.28 & -6.87 \\ 
  0.80 & RRPO-C, $\alpha = 0.7$ & 104328.93 & -1.33 \\ 
  0.80 & RRPO-C, $\alpha = 0.8$ & 111791.18 & 5.73 \\ 
  0.80 & RRPO-C, $\alpha = 0.9$ & 122249.99 & 15.62 \\ 
%  0.80 & RRPO-C, $\alpha = 1$ & 146523.03 & 38.58 \\ 
  0.80 & RRPO & 149709.04 & 41.59 \\ \midrule
  1.00 & RRPO-C, $\alpha = 0.5$ & 74421.62 & -14.44 \\ 
  1.00 & RRPO-C, $\alpha = 0.6$ & 77826.59 & -10.52 \\ 
  1.00 & RRPO-C, $\alpha = 0.7$ & 90332.98 & 3.86 \\ 
  1.00 & RRPO-C, $\alpha = 0.8$ & 103167.94 & 18.61 \\ 
  1.00 & RRPO-C, $\alpha = 0.9$ & 107505.80 & 23.60 \\ 
%  1.00 & RRPO-C, $\alpha = 1$ & 122010.14 & 40.28 \\ 
  1.00 & RRPO & 125987.02 & 44.85 \\ \midrule
  1.50 & RRPO-C, $\alpha = 0.5$ & 53040.97 & -6.08 \\ 
  1.50 & RRPO-C, $\alpha = 0.6$ & 56740.90 & 0.47 \\ 
  1.50 & RRPO-C, $\alpha = 0.7$ & 62651.62 & 10.94 \\ 
  1.50 & RRPO-C, $\alpha = 0.8$ & 66842.72 & 18.36 \\ 
  1.50 & RRPO-C, $\alpha = 0.9$ & 73892.42 & 30.84 \\ 
%  1.50 & RRPO-C, $\alpha = 1$ & 81976.27 & 45.16 \\ 
  1.50 & RRPO & 82880.96 & 46.76 \\ \midrule
  2.00 & RRPO-C, $\alpha = 0.5$ & 35416.94 & -4.70 \\ 
  2.00 & RRPO-C, $\alpha = 0.6$ & 38862.06 & 4.57 \\ 
  2.00 & RRPO-C, $\alpha = 0.7$ & 41834.74 & 12.57 \\ 
  2.00 & RRPO-C, $\alpha = 0.8$ & 45750.84 & 23.10 \\ 
  2.00 & RRPO-C, $\alpha = 0.9$ & 50007.72 & 34.56 \\ 
%  2.00 & RRPO-C, $\alpha = 1$ & 54356.68 & 46.26 \\ 
  2.00 & RRPO & 54665.15 & 47.09 \\ \bottomrule
\end{tabular}
\caption{Results of comparison between RRPO, constrained RRPO and DRPO on \orangeJuice instances with semi-log demand model. \label{table:constrained_semilog}}
\end{table}

%%%%
% END OF CONSTRAINED DISTRIBUTION SECTION
%%%%

%%%% 
% START OF DISCRETIZATION SECTION
%%%%

\clearpage
\pagebreak

\section{Approximation of convex price sets by discrete price sets}
\label{sec:approximation_convex_discrete}

In this section, we show that under certain conditions, the optimal randomized robust pricing scheme over a convex price set $\Pcal$ can be well-approximated by the optimal randomized robust pricing scheme over a finite subset of $\Pcal$. 

Let $Z^*_{\RR}(\Pcal')$ denote the optimal objective value of the RRPO problem with respect to the price set $\Pcal'$, which could be a finite set or an uncountable set:
\begin{equation}
Z^*_{\RR}(\Pcal') = \sup_{F \in \Fcal_{\Pcal'}} \inf_{\ub \in \Ucal} \int R(\pb, \ub) \, dF(\pb).
\end{equation}
We additionally make the following assumption regarding the revenue function $R$. 
\begin{assumption}
There exists $L > 0$ such that for all $\ub \in \Ucal$, $R(\cdot, \ub)$ is Lipschitz continuous with constant $L$, i.e.,
\begin{equation}
| R(\pb, \ub) - R(\pb', \ub)| \leq L \| \pb - \pb' \|_{\infty},
\end{equation}
for all $\pb, \pb' \in \Pcal$. \label{assumption:R_lipschitz}
\end{assumption}

Our main result can be stated as follows.
\begin{theorem}
Suppose that $\Pcal$ and $\Ucal$ are nonempty compact convex sets. Suppose that $R(\pb,\ub)$ is convex in $\ub$. Suppose also that Assumption~\ref{assumption:R_lipschitz} holds. 

Then, for any $\epsilon > 0$ there exists a finite subset $\hat{\Pcal} \subseteq \Pcal$, such that 
\begin{equation}
Z^*_{\RR}(\Pcal) - Z^*_{\RR}(\hat{\Pcal}) \leq \epsilon.
\end{equation} \label{theorem:discrete_vs_continuous_Lipschitz}
\end{theorem}

To establish the result, we need to establish several auxiliary results. The first result states that if we are given a finite $\delta$-net of $\Pcal$, the nominal problem for a given $\ub$ with the exact set $\Pcal$ and the $\delta$-net of $\Pcal$ will be close in their optimal objective values.

\begin{lemma}
Let $\Pcal \subseteq \Rbb^I$ be a nonempty set of price vectors, and let $\hat{\Pcal} \subseteq \Pcal$ be a finite set such that for every $\pb \in \Pcal$, there exists $\pb' \in \hat{\Pcal}$ such that $\| \pb - \pb' \|_{\infty} \leq \delta$. Suppose also that Assumption~\ref{assumption:R_lipschitz} holds. Then 
\begin{equation}
\sup_{\pb \in \Pcal} R(\pb, \ub) - \max_{\pb \in \hat{\Pcal}} R(\pb, \ub) \leq L \delta
\end{equation}
for all $\ub \in \Ucal$. \label{lemma:max_over_p_Ldelta}
\end{lemma}
\begin{proofvvm}
Fix any $\epsilon > 0$ and any $\ub \in \Ucal$. Let $\pb^*$ be such that $R(\pb^*, \ub) \geq \sup_{\pb \in \Pcal} R(\pb, \ub) - \epsilon$ or equivalently, $R(\pb^*, \ub) + \epsilon \geq \sup_{\pb \in \Pcal} R(\pb, \ub)$, which is guaranteed to exist because of the definition of the supremum. 

Let $\hat{\pb} \in \arg \max_{\pb \in \hat{\Pcal}} R(\pb, \ub)$. Note that $\hat{\pb}$ is well-defined, because $\hat{\Pcal}$ is a finite set. 
%Note that both sets of maximizers are well-defined, because $R$ is assumed to be Lipschitz continuous and $\Pcal$ is compact, and so by the extreme value theorem the maximum of $R(\cdot, \ub)$ is attained, whereas $\hat{\Pcal}$ is finite, and so the maximum of $R(\cdot, \ub)$ over $\hat{\Pcal}$ must also be attained. 

Let $\tilde{\pb} \in \hat{\Pcal}$ be a price vector such that $\| \tilde{\pb} - \pb^* \|_{\infty} \leq \delta$, which is guaranteed to exist by the hypotheses of the lemma. We now have
\begin{align*}
& \sup_{\pb \in \Pcal} R(\pb, \ub) - \max_{\pb \in \hat{\Pcal}} R(\pb, \ub) \\
& = R(\pb^*, \ub) + \epsilon - R(\hat{\pb}, \ub) \\
& = R(\pb^*, \ub) - R(\tilde{\pb}, \ub) + R(\tilde{\pb}, \ub) - R(\hat{\pb}, \ub) + \epsilon\\
& \leq | R(\pb^*, \ub) - R(\tilde{\pb}, \ub) | + 0 + \epsilon \\
& \leq L \delta + \epsilon,
\end{align*}
where the first step follows by the definitions of $\pb^*$ and $\hat{\pb}$; the second step follows by algebra; the third step follows by basic properties of absolute values, and also the fact that $\tilde{\pb}$ is a feasible but not necessarily optimal solution to $\max_{\pb \in \hat{\Pcal}} R(\pb, \ub)$, and hence $R(\tilde{\pb}, \ub) \leq R(\hat{\pb}, \ub)$ (or equivalently, $R(\tilde{\pb}, \ub) - R(\hat{\pb}, \ub) \leq 0$); and the last step follows by Assumption~\ref{assumption:R_lipschitz}. 

Since $\epsilon > 0$ is arbitrary, it must therefore be that $\sup_{\pb \in \Pcal} R(\pb, \ub) - \max_{\pb \in \hat{\Pcal}} R(\pb, \ub) \leq L \delta$, as required. \Halmos
\end{proofvvm}

The second result says that given a $\delta$-net of $\Pcal$, the objective value of a certain min-max problem, where the outer minimization is over $\ub$ and the inner maximization is over $\Pcal$, will not change too much if the inner maximization is performed over the $\delta$-net.

\begin{lemma}
Let $\Pcal \subseteq \Rbb^I$ be a nonempty set of price vectors, and let $\hat{\Pcal} \subseteq \Pcal$ be a finite set such that for every $\pb \in \Pcal$, there exists $\pb' \in \hat{\Pcal}$ such that $\| \pb - \pb' \|_{\infty} \leq \delta$. Suppose also that Assumption~\ref{assumption:R_lipschitz} holds. Then 
\begin{equation}
\inf_{\ub \in \Ucal} \sup_{\pb \in \Pcal}  R(\pb, \ub) - \inf_{\ub \in \Ucal} \max_{\pb \in \hat{\Pcal}} R(\pb, \ub) \leq L \delta.
\end{equation} \label{lemma:inf_sup_leq_L_delta}
\end{lemma}

\begin{proofvvm}
Fix any $\epsilon > 0$. Let $\hat{\ub} \in \Ucal$ be an uncertain parameter vector such that $\max_{\pb \in \hat{\Pcal}} R(\pb, \hat{\ub}) \leq \inf_{\ub \in \Ucal} \max_{\pb \in \hat{\Pcal}} R(\pb, \ub) + \epsilon$. Such a $\hat{\ub}$ is guaranteed to exist by the definition of the infimum. 

We now have:
\begin{align}
& \inf_{\ub \in \Ucal} \sup_{\pb \in \Pcal} R(\pb, \ub) - \inf_{\ub \in \Ucal} \max_{\pb \in \hat{\Pcal}} R(\pb, \ub) \\
& \leq \inf_{\ub \in \Ucal} \sup_{\pb \in \Pcal} R(\pb, \ub) - \max_{\pb \in \hat{\Pcal}} R(\pb, \hat{\ub}) + \epsilon \\
& \leq \sup_{\pb \in \Pcal} R(\pb, \hat{\ub}) - \max_{\pb \in \hat{\Pcal}} R(\pb, \hat{\ub}) + \epsilon \\
& \leq L \delta + \epsilon,
\end{align}
where the first step follows by the definition of $\hat{\ub}$ as an $\epsilon$-approximate solution to the discretized problem $\inf_{\ub \in \Ucal} \max_{\pb \in \hat{\Pcal}} R(\pb, \ub)$; the second step follows because $\hat{\ub}$ is a feasible solution to the continuous (non-discretized) problem $\inf_{\ub \in \Ucal} \sup_{\pb \in \Pcal} R(\pb, \ub)$; and the final step by Lemma~\ref{lemma:max_over_p_Ldelta}. Since $\epsilon$ was arbitrary, it must be that $\inf_{\ub \in \Ucal} \sup_{\pb \in \Pcal} R(\pb, \ub) - \inf_{\ub \in \Ucal} \max_{\pb \in \hat{\Pcal}} R(\pb, \ub)  \leq L \delta$ as required. \Halmos
\end{proofvvm}

Finally, we prove Theorem~\ref{theorem:discrete_vs_continuous_Lipschitz}. 

\begin{proof}{Proof of Theorem~\ref{theorem:discrete_vs_continuous_Lipschitz}:}
Let $\delta = \epsilon / L$. Since $\Pcal$ is compact, there exists a finite set $\hat{\Pcal}$ which is a $\delta$-net, that is, for every $\pb \in \Pcal$, there exists $\pb' \in \hat{\Pcal}$ such that 
\begin{equation}
\| \pb - \pb' \|_{\infty} \leq \delta.
\end{equation}

We now have
\begin{align}
& Z^*_{\RR}(\Pcal) - Z^*_{\RR}(\hat{\Pcal}) \\
& = \sup_{F \in \Fcal(\Pcal)} \inf_{\ub \in \Ucal} \int_{\Pcal} R(\pb, \ub) \, dF(\pb) -  \sup_{F \in \Fcal(\hat{\Pcal})} \inf_{\ub \in \Ucal} \int_{\Pcal} R(\pb, \ub) \, dF(\pb) \\
& = \inf_{\ub \in \Ucal} \sup_{F \in \Fcal(\Pcal)} \int_{\Pcal} R(\pb, \ub) \, dF(\pb) -  \inf_{\ub \in \Ucal} \sup_{F \in \Fcal(\hat{\Pcal})} \int_{\Pcal} R(\pb, \ub) \, dF(\pb) \\
& = \inf_{\ub \in \Ucal} \sup_{\pb \in \Pcal} R(\pb, \ub)  -  \inf_{\ub \in \Ucal} \max_{\pb \in \hat{\Pcal}} R(\pb, \ub)  \\ 
& \leq L \delta \\
& = \epsilon,
\end{align}
where the first step follows by the definition of $Z^*_{\RR}(\Pcal')$, the second by applying Sion's minimax theorem, the third by the equivalence of $\sup_{F \in \Fcal(\Pcal)} \int R(\pb, \ub) \, dF(\pb)$ and $\sup_{\pb \in \Pcal} R(\pb, \ub)$, the fourth by Lemma~\ref{lemma:inf_sup_leq_L_delta}, and the final step by the definition of $\delta$. (Note that in the second step, Sion's minimax theorem applies as $\int_\Pcal R(\pb, \ub) \, dF(\pb)$ is convex in $\ub$ (and hence both quasiconvex and continuous in $\ub$), and linear in $F$ (and hence both quasiconcave and continuous in the topology of weak convergence in $F$), and $\Fcal(\Pcal)$ and $\Ucal$ are both convex, and $\Ucal$ is additionally compact.) \Halmos \\
\end{proof}

We now pause to make a few comments about Theorem~\ref{theorem:discrete_vs_continuous_Lipschitz}. First, the Lipschitz continuity assumption on the revenue function $R$ is not particularly restrictive, and is satisfied by many popular demand models. In particular, when $\Pcal$ and $\Ucal$ are compact sets and the gradient $\nabla_\pb R(\pb, \ub)$ with respect to $\pb$ is continuous in both $\pb$ and $\ub$, then $R$ will satisfy the requirement of the theorem for a Lipschitz constant $L$ that is defined as
\begin{equation}
L = \max_{\ub \in \Ucal} \max_{\pb \in \Pcal} \| \nabla_{\pb} R(\pb, \ub) \|_1.
\end{equation}
For the semi-log and log-log demand models, the revenue function $R$ has a gradient $\nabla_{\pb} R(\pb, \ub)$ that is a continuous function of $\pb$ and $\ub$.

Second, this theorem implies that when the revenue function is Lipschitz continuous with the same constant across every $\ub$, discretization will perform well. In particular, for a price set of the form $\Pcal = [p_{\min,i}, p_{\max,i}]^I$, a uniform discretization of the form $\hat{\Pcal} = \{ \delta \lceil p_{\min,i} / \delta \rceil, \delta (\lceil p_{\min,i} / \delta \rceil + 1), \delta (\lceil p_{\min,i} / \delta \rceil + 2), \dots, \delta \lfloor p_{\max,i} / \delta \rfloor \}^I$ for some $\delta > 0$ will yield the $\delta$-net in the $\ell_{\infty}$ norm required in the proof of Theorem~\ref{theorem:discrete_vs_continuous_Lipschitz}. By setting $\delta$ to be suitably small, we can thus ensure that the discretized RRPO problem returns a randomized pricing scheme that is arbitrarily close to the optimal randomized pricing scheme over the continuous price set $\Pcal$. This result thus provides a theoretical justification for our focus on discrete price sets in Sections~\ref{sec:finitePcal_convexUcal} and \ref{sec:finitePcal_finiteUcal} of the paper.

\section{Solution method for finite $\Pcal$, finite $\Ucal$}
\label{sec:finitePcal_finiteUcal}

The second solution approach we consider is for the case where both $\Pcal$ and $\Ucal$ are finite sets. In particular, we assume that the uncertainty set $\Ucal$ is a binary representable set. For fixed positive integers $m$ and $n$, we let $\Ucal$ be defined as 
\begin{equation}
\Ucal = \{ \ub = \Fb \zb \mid \Ab \zb \leq \bb, \zb \in \{0,1\}^n \},
\end{equation}
where $\bb$ is a $m$-dimensional real vector, $\Ab$ is a $m$-by-$n$ real matrix and $\Fb$ is a $M$-by-$n$ real matrix, where $M$ is the dimension of the uncertain parameter vector $\ub$. We note here that any finite set can be represented in this way: if $\Ucal = \{\ub_1, \dots, \ub_H \}$, then one can define $\zb$ as $\zb = (z_1,\dots, z_H) \in \{0,1\}^H$, use the linear inequality system $\Ab \zb \leq \bb$ to represent the constraints $\oneb^T \zb \leq 1$, $- \oneb^T \zb \leq -1$, and define $\Fb = [ \ub_1 \ \dots \ \ub_H ]$ (i.e., the matrix obtained from stacking the vectors $\ub_1,\dots, \ub_H$ as columns). 

From a modeling standpoint, a discrete set $\Ucal$ can be used to model the firm's uncertainty in different ways. One could use this type of uncertainty set to model a finite set of scenarios in terms of the customer population's demand response to the firm's prices. Another possibility is that the firm faces a collection of different customer segments, where each one corresponds to a particular type of demand model, and at most one segment will respond to the firm's products; in this setting, the firm could use the uncertainty set to model its uncertainty as to which segment will respond. Given the factor structure $\ub = \Fb \zb$ within $\Ucal$, one could also consider a setting where the firm estimates a richer demand model where the elasticity parameters (represented by $\ub$) are dependent on macroeconomic factors (represented by $\zb$) that are uncertain and represented in a discrete way; for example, one could use binary variables $z_1,z_2,z_3$ to indicate three different unemployment scenarios, and $z_4,z_5,z_6$ to indicate three different interest rate scenarios. The firm would then obtain a randomized pricing strategy that accounts for uncertainty in these macroeconomic factors. 

Recall that when $\Pcal$ is finite, then the RRPO problem is
\begin{align}
Z^*_{\RR} = \max_{\pib \in \Delta_{\Pcal}} \min_{\ub \in \Ucal} \sum_{\pb \in \Pcal} \pi_{\pb} R(\pb, \ub). \label{prob:RRPO_DCG_primal}
\end{align}
% \begin{subequations}
% \begin{alignat}{2}
% & \text{maximize} & \quad &  \min_{\ub \in \Ucal} \sum_{\pb \in \Pcal} \pi_{\pb} R(\pb, \ub)
% \end{alignat}
% \end{subequations}
We can transform this problem into a dual problem where the outer problem is to optimize a distribution over uncertain parameter vectors, and the inner problem is to optimize over the price vector, as follows:
\begin{align}
Z^*_{\RR} & = \max_{\pib \in \Delta_{\Pcal}} \min_{\ub \in \Ucal} \sum_{\pb \in \Pcal} \pi_{\pb} R(\pb, \ub) \\
& = \max_{\pib \in \Delta_{\Pcal}} \min_{\lambdab \in \Delta_{\Ucal}} \sum_{\pb \in \Pcal} \sum_{\ub \in \Ucal} \pi_{\pb} \lambda_{\ub} R(\pb, \ub) \\
& = \min_{\lambdab \in \Delta_{\Ucal}} \max_{\pib \in \Delta_{\Pcal}}  \sum_{\pb \in \Pcal} \sum_{\ub \in \Ucal} \pi_{\pb} \lambda_{\ub} R(\pb, \ub) \\
& = \min_{\lambdab \in \Delta_{\Ucal}} \max_{\pb \in \Pcal} \sum_{\ub \in \Ucal} \lambda_{\ub} R(\pb, \ub), \label{prob:RRPO_DCG_dual}
\end{align}
where the first equality follows because minimization of a function of $\ub$ over the finite set $\Ucal$ is the same as minimizing the expected value of that function over all probability mass functions supported on $\Ucal$; the second equality follows by linear programming duality; and the final equality follows because maximization of a function of $\pb$ over $\Pcal$ is the same as maximizing the expected value of that function over all probability mass functions supported on $\Pcal$. We refer to problem~\eqref{prob:RRPO_DCG_primal} as the \emph{primal} problem and \eqref{prob:RRPO_DCG_dual} as the \emph{dual} problem. 

Consider now the \emph{restricted primal problem}, where we replace $\Pcal$ with a subset $\hat{\Pcal} \subseteq \Pcal$ in problem~\eqref{prob:RRPO_DCG_primal}, and the \emph{restricted dual problem}, where we replace $\Ucal$ with a subset $\hat{\Ucal} \subseteq \Ucal$ in problem~\eqref{prob:RRPO_DCG_dual}. Let us denote the objective values of these two problems with $Z_{P, \hat{\Pcal}}$ and $Z_{D, \hat{\Ucal}}$, respectively. These two problems are:
\begin{align}
Z_{P, \hat{\Pcal}} = & \max_{\pib \in \Delta_{\hat{\Pcal}}} \min_{\ub \in \Ucal} \sum_{\pb \in \hat{\Pcal}} \pi_{\pb} R(\pb, \ub), \label{prob:restricted_primal} \\
Z_{D, \hat{\Ucal}} = & \min_{\lambdab \in \Delta_{\hat{\Ucal}}}  \max_{\pb \in \Pcal} \sum_{\ub \in \hat{\Ucal}} \lambda_{\ub} R(\pb, \ub). \label{prob:restricted_dual} 
\end{align}
Observe that $Z_{P,\hat{\Pcal}}$ and $Z_{D, \hat{\Ucal}}$ bound $Z^*_{\RR}$ from below and above, that is,
\begin{equation*}
Z_{P, \hat{\Pcal}} \leq Z^*_{\RR} \leq Z_{D, \hat{\Ucal}}.
\end{equation*}
In the above, the justification for the first inequality is because maximizing over distributions supported on the smaller set of price vectors $\hat{\Pcal}$ cannot result in a higher worst-case objective than solving the full primal problem with $\Pcal$, which gives the value $Z^*_{\RR}$. The second inequality similarly follows because minimizing over distributions supported on the smaller set of uncertainty realizations $\hat{\Ucal}$ cannot result in a lower worst-case objective than solving the full dual problem with $\Ucal$, which gives $Z^*_{\RR}$. 

The idea of double column generation is as follows. Let us pick some subset of price vectors $\hat{\Pcal} \subseteq \Pcal$ and some subset of uncertainty realizations $\hat{\Ucal} \subseteq \Ucal$. Observe that the restricted primal problem~\eqref{prob:restricted_primal} for $\hat{\Pcal}$ can be written in epigraph form as
\begin{subequations}
\begin{alignat}{2}
& \underset{\pib, t}{\text{maximize}} & \quad & t \\
& \text{subject to} & & t \leq \sum_{\pb \in \hat{\Pcal}} \pi_{\pb} R(\pb, \ub), \quad \forall\ \ub \in \Ucal, \\
& & & \sum_{\pb \in \hat{\Pcal}} \pi_{\pb} = 1, \\
& & & \pi_{\pb} \geq 0, \quad \forall \ \pb \in \hat{\Pcal}.
\end{alignat}
\label{prob:restricted_primal_epigraph}%
\end{subequations}
This problem has a huge number of constraints (one for each $\ub \in \Ucal$). However, we can solve it using delayed constraint generation, starting from the set $\hat{\Ucal}$. Upon solving it in this way, at termination we will have a subset $\Ucal'$ of uncertainty realizations from $\Ucal$ that were found during the constraint generation process. We update $\hat{\Ucal}$ to be equal to $\Ucal'$. 

With this (updated) subset $\hat{\Ucal}$ in hand, we now solve the restricted dual problem~\eqref{prob:restricted_dual} for $\hat{\Ucal}$. This problem can be written in epigraph form as 
\begin{subequations}
\begin{alignat}{2}
& \underset{\lambdab, \rho}{\text{minimize}} & \quad & \rho \\
& \text{subject to} & & \rho \geq \sum_{\ub \in \hat{\Ucal}} \lambda_{\ub} R(\pb, \ub), \quad \forall\ \pb \in \Pcal, \\
& & & \sum_{\ub \in \hat{\Ucal}} \lambda_{\ub} = 1, \\
& & & \lambda_{\ub} \geq 0, \quad \forall \ \ub \in \hat{\Ucal}.
\end{alignat}
\label{prob:restricted_dual_epigraph}%
\end{subequations}
This problem also has a huge number of constraints, but again we can solve it using delayed constraint generation, with the initial subset of price vectors set to $\hat{\Pcal}$. At termination, we will have a new subset $\Pcal'$ of price vectors, which will contain the original set of price vectors in $\hat{\Pcal}$. We then update $\hat{\Pcal}$ to $\Pcal'$, and go back to solving the restricted primal problem. The process then repeats: after solving the restricted primal, we will have a new (bigger) $\hat{\Ucal}$; we then solve the restricted dual, after which we have a new (bigger) $\hat{\Pcal}$; we then go back to the restricted primal, and so on. After each iteration of solving the restricted primal and restricted dual, the set $\hat{\Pcal}$ expands and the set $\hat{\Ucal}$ expands. Thus, the bounds $Z_{P, \hat{\Pcal}}$ and $Z_{D,\hat{\Ucal}}$ get closer and closer to $Z^*_{\RR}$. The algorithm can then be terminated either when $Z_{P, \hat{\Pcal}} = Z_{D, \hat{\Ucal}}$, which would imply that both restricted primal and restricted dual objective values exactly coincide with $Z^*_{\RR}$; or otherwise, one can terminate when $Z_{D, \hat{\Ucal}} - Z_{P,\hat{\Pcal}} < \epsilon$, where $\epsilon > 0$ is a user specified tolerance. 

The overall algorithmic approach is formalized as Algorithm~\ref{alg:doubleCG}. This algorithm invokes two procedures, \PrimalCG (Algorithm~\ref{alg:primalCG}) and \DualCG (Algorithm~\ref{alg:dualCG}), which are delayed constraint generation algorithms for solving the restricted primal and dual problems respectively. We note that Algorithm~\ref{alg:doubleCG} is an adaptation of the double column generation algorithm of \cite{wang2020randomized} for the randomized robust assortment optimization problem, which is itself adapted from the double column generation algorithm of \cite{delage2022value} for solving mixed-integer distributionally robust optimization problems. The proof of correctness of this procedure follows similarly to \cite{delage2022value}, and is omitted for brevity. The novelty in our approach lies in how we handle the separation problems which are at the heart of \PrimalCG and \DualCG, which we discuss next.  Sections~\ref{subsec:finitePcal_finiteUcal_linear}, \ref{subsec:finitePcal_finiteUcal_semilog} and \ref{subsec:finitePcal_finiteUcal_loglog} provide the details for the linear, semi-log and log-log demand models. Section~\ref{subsec:finitePcal_finiteUcal_complexity} discusses the complexity of randomization and provides a result that is analogous to Theorem~\ref{theorem:complexity_convexUcal} for the case of a finite $\Pcal$ and finite $\Ucal$. 

\begin{algorithm}
\caption{Double column generation method for solving the finite $\Pcal$, finite $\Ucal$ RRPO problem. \label{alg:doubleCG}}
\begin{algorithmic}[1]
\STATE Initialize $\hat{\Pcal}$ to be a non-empty subset of $\Pcal$, and $\hat{\Ucal}$ to be a non-empty subset of $\Ucal$.
\STATE Set $\LB \gets -\infty$, $\UB \gets +\infty$
\REPEAT
	\STATE Run $\PrimalCG(\hat{\Pcal}, \hat{\Ucal})$ to solve the restricted primal problem with $\hat{\Pcal}$ and with $\hat{\Ucal}$ as the initial uncertainty set. Let the objective value be $Z_{P, \hat{\Pcal}}$ and the new uncertainty set be $\Ucal'$.
	\STATE Set $\hat{\Ucal} \gets \Ucal'$.
	\STATE Set  $\LB \gets Z_{P, \hat{\Pcal}}$.
	\STATE Run $\DualCG(\hat{\Pcal}, \hat{\Ucal})$ to solve the restricted dual problem with $\hat{\Ucal}$ and with $\hat{\Pcal}$ as the initial price vector set. Let the objective value be $Z_{D, \hat{\Ucal}}$ and the new price vector set be $\Pcal'$.
	\STATE Set $\hat{\Pcal} \gets \Pcal'$.
	\STATE Set $\UB \gets Z_{D, \hat{\Ucal}}$. 
\UNTIL{ $\UB - \LB < \epsilon$}
\end{algorithmic}
\end{algorithm}

\begin{algorithm}
\caption{\PrimalCG procedure.\label{alg:primalCG}}
\begin{algorithmic}[1]
\STATE Initialize $\Ucal' \gets \hat{\Ucal}$
\REPEAT
	\STATE Solve the doubly restricted primal problem:
	\begin{subequations}
	\begin{alignat}{2}
	& \underset{\pib, t}{\text{maximize}} & \quad & t \\
	& \text{subject to} & & t \leq \sum_{\pb \in \hat{\Pcal}} \pi_{\pb} R(\pb, \ub), \quad \forall \ub \in \Ucal', \\
	& & & \sum_{\pb \in \hat{\Pcal}} \pi_{\pb} = 1, \\
	& & & \pi_{\pb} \geq 0, \quad \forall \ \pb \in \hat{\Pcal}. 
	\end{alignat}
	\label{prob:primal_DR_epigraph}
	\end{subequations}
	Let $(\pib, t^*)$ be the optimal solution of the doubly restricted problem. 
	\STATE Solve the primal separation problem:
	\begin{align}
	 \min_{\ub \in \Ucal} \sum_{\pb \in \hat{\Pcal}} \pi_{\pb} \cdot R(\pb, \ub). \label{prob:primal_separation}
	\end{align}
	Let $t'$ and $\ub^*$ be the optimal objective value and solution of this separation problem.
	\IF{$t^* > t'$} 
		\STATE Set $\Ucal' \gets \Ucal' \cup \{ \ub^* \}$
	\ENDIF
\UNTIL{ $t^* \leq t'$ }
\STATE Set $Z_{P, \hat{\Pcal}} \gets t^*$
\RETURN $(Z_{P, \hat{\Pcal}}, \Ucal')$. 
\end{algorithmic}
\end{algorithm}

\begin{algorithm}
\caption{\DualCG procedure. \label{alg:dualCG}}
\begin{algorithmic}[1]
\STATE Initialize $\Pcal' \gets \hat{\Pcal}$
\REPEAT
	\STATE Solve the doubly restricted dual problem:
	\begin{subequations}
	\begin{alignat}{2}
	& \underset{\lambdab, \rho}{\text{minimize}} & \quad & \rho \\
	& \text{subject to} & & \rho \geq \sum_{\ub \in \hat{\Ucal}} \lambda_{\ub} R(\pb, \ub), \quad \forall \pb \in \Pcal', \\
	& & & \sum_{\ub \in \hat{\Ucal}} \lambda_{\ub} = 1, \\
	& & & \lambda_{\ub} \geq 0, \quad \forall \ \ub \in \hat{\Ucal}. 
	\end{alignat}
	\label{prob:dual_DR_epigraph}
	\end{subequations}
	Let $(\lambdab, \rho^*)$ be the optimal solution of the doubly restricted problem. 
	\STATE Solve the dual separation problem:
	\begin{align}
	 \max_{\pb \in \Pcal} \sum_{\ub \in \hat{\Ucal}} \lambda_{\ub} \cdot R(\pb, \ub) \label{prob:dual_separation}
	\end{align}
	Let $\rho'$ and $\pb^*$ be the optimal objective value and solution of this separation problem.
	\IF{$\rho^* < \rho'$}
		\STATE Set $\Pcal' \gets \Pcal' \cup \{ \pb^* \}$
	\ENDIF
\UNTIL{ $\rho^* \geq \rho'$ }
\STATE Set $Z_{D, \hat{\Ucal}} \gets \rho^*$
\RETURN $(Z_{D, \hat{\Ucal}}, \Pcal')$. 
\end{algorithmic}
\end{algorithm}

Note that the doubly restricted primal and dual problems~\eqref{prob:primal_DR_epigraph} and \eqref{prob:dual_DR_epigraph} solved in \PrimalCG and \DualCG are both linear programs, and can be thus be solved easily. The principal difficulty in these procedures comes from the primal and dual separation problems~\eqref{prob:primal_separation} and \eqref{prob:dual_separation}, which require optimizing over an uncertain parameter vector $\ub \in \Ucal$ and a price vector $\pb \in \Pcal$ respectively. In the following sections, we discuss how these two separation problems can be tackled for the linear, semi-log and log-log demand models. Note that in all three sections, we continue to make Assumption~\ref{assumption:Pcal_Cartesian_product}, which states that $\Pcal$ can be written as the Cartesian product of finite sets of prices for each product, i.e., $\Pcal = \Pcal_1 \times \dots \times \Pcal_I$, where $\Pcal_1, \dots, \Pcal_I$ are finite sets.

\subsection{Primal and dual subproblems for linear demand model}
\label{subsec:finitePcal_finiteUcal_linear}

For the linear demand model, the primal separation problem is 
\begin{align}
\min_{\ub \in \Ucal} \sum_{\pb \in \hat{\Pcal}} \pi_{\pb} \cdot \left[ \sum_{i=1}^I p_i \cdot (\alpha_i - \beta_i p_i + \sum_{j \neq i} \gamma_{i,j} p_j ) \right].
\end{align}
Note that this objective function is linear in $\ub = (\alphab, \betab, \gammab)$. Therefore, the whole problem can be expressed as 
\begin{subequations}
\begin{alignat}{2}
& \text{minimize} & \quad & \sum_{\pb \in \hat{\Pcal}} \pi_{\pb} \cdot \left[ \sum_{i=1}^I p_i \cdot (\alpha_i - \beta_i p_i + \sum_{j \neq i} \gamma_{i,j} p_j ) \right] \\
& \text{subject to} & & \ub = \Fb \zb, \\
& & & \Ab \zb \leq \bb, \\
& & & \zb \in \{0,1\}^n,
\end{alignat}
\end{subequations}
which is a mixed-integer linear program. 

The dual separation problem is
\begin{align}
\max_{\pb \in \Pcal} \sum_{\ub \in \hat{\Ucal}} \lambda_{\ub} \cdot \left[ \sum_{i=1}^I p_i (\alpha_i - \beta_i p_i + \sum_{j \neq i} \gamma_{i,j} p_j ) \right].
\end{align}
By introducing the same binary variables as in the separation problem~\eqref{prob:linear_MILP} (the linear demand separation problem for the convex $\Ucal$ setting), we obtain the following mixed-integer linear program:
\begin{subequations}
\begin{alignat}{2}
& \underset{\xb, \yb}{\text{maximize}} & &  \sum_{\ub \in \hat{\Ucal}} \lambda_{\ub} \cdot \left[ \sum_{i=1}^I \sum_{t \in \Pcal_i} \alpha_i \cdot t \cdot x_{i,t} - \sum_{i=1}^I \sum_{t \in \Pcal_i} t^2 \cdot \beta_i \cdot x_{i,t} + \sum_{i=1}^I \sum_{ j \neq i} \sum_{t_1 \in \Pcal_i} \sum_{t_2 \in \Pcal_j} \gamma_{i,j} \cdot t_1 \cdot t_2 \cdot y_{i,j,t_1,t_2}  \right] \\
& \text{subject to} & \quad & \sum_{t \in \Pcal_i} x_{i,t} = 1, \quad \forall \ i \in [I], \\
& & & \sum_{t_2 \in \Pcal_j} y_{i,j,t_1, t_2} = x_{i,t_1}, \quad \forall \ i, j \in [I], j \neq i, t_2 \in \Pcal_j, \\
& & & \sum_{t_1 \in \Pcal_i} y_{i,j,t_1, t_2} = x_{j,t_1}, \quad \forall \ i, j \in [I], j \neq i, t_1 \in \Pcal_i, \\
& & & x_{i,t} \in \{0,1\}, \quad \forall \ i \in [I], \ t \in \Pcal_i, \\
& & & y_{i,j,t_1,t_2} \in \{0,1\}, \quad \forall \ i, j\in [I], i \neq j, t_1 \in \Pcal_i, t_2 \in \Pcal_j.
\end{alignat}
\end{subequations}
Importantly, note that the size of this problem does not scale with the number of uncertainty realizations inside $\hat{\Ucal}$; the form of this problem is equivalent to problem~\eqref{prob:linear_MILP} where $\ub$ is replaced with $\sum_{\ub \in \hat{\Ucal}} \lambda_{\ub} \cdot \ub$ (the ``average'' uncertain demand parameter). As we will see in the next couple of sections, the same will not be true for the semi-log and log-log demand models.

\subsection{Primal and dual subproblems for semi-log demand model}
\label{subsec:finitePcal_finiteUcal_semilog}

For the semi-log demand model, the primal separation problem is 
\begin{align}
\min_{\ub \in \Ucal} \sum_{\pb \in \hat{\Pcal}} \pi_{\pb} \cdot \left[ \sum_{i=1}^I p_i \cdot e^{\alpha_i - \beta_i p_i + \sum_{j \neq i} \gamma_{i,j} p_j } \right].
\end{align}
Note that this objective function is convex in $\ub = (\alphab, \betab, \gammab)$, because the weights $\pi_{\pb}$ and $p_i$ for a given $\pb \in \Pcal$ and $i \in [I]$ are nonnegative, and because the function $e^{\alpha_i - \beta_i p_i + \sum_{j \neq i} \gamma_{i,j} p_j }$ is convex in $\ub = (\alphab, \betab, \gammab)$. Thus, the whole problem can be expressed as 
\begin{subequations}
\begin{alignat}{2}
& \underset{\ub, \zb}{\text{minimize}} & \quad & \sum_{\pb \in \hat{\Pcal}} \pi_{\pb} \cdot \left[ \sum_{i=1}^I p_i \cdot e^{\alpha_i - \beta_i p_i + \sum_{j \neq i} \gamma_{i,j} p_j } \right] \\
& \text{subject to} & & \ub = \Fb \zb, \\
& & &  \Ab \zb \leq \bb, \\
& & &  \zb \in \{0,1\}^n,
\end{alignat}
\end{subequations}
which can be re-written as a mixed-integer exponential cone program.

The dual separation problem is 
\begin{align}
\max_{\pb \in \Pcal} \sum_{\ub \in \hat{\Ucal}} \lambda_{\ub} \cdot \left[ \sum_{i=1}^I p_i \cdot e^{\alpha_i - \beta_i p_i + \sum_{j \neq i} \gamma_{i,j} p_j }  \right].
\end{align}
The objective function of this problem is in general not concave in $\pb$. However, just as in Section~\ref{subsec:finitePcal_convexUcal_semilog}, the related problem of optimizing the logarithm of this objective, which is 
\begin{align}
& \max_{\pb \in \Pcal} \log \left[ \sum_{\ub \in \hat{\Ucal}} \lambda_{\ub} \cdot \left[ \sum_{i=1}^I p_i \cdot e^{\alpha_i - \beta_i p_i + \sum_{j \neq i} \gamma_{i,j} p_j }  \right] \right] \\
& = \max_{\pb \in \Pcal} \log \left[ \sum_{\ub \in \hat{\Ucal}} \sum_{i=1}^I e^{\log \lambda_{\ub} + \log p_i +  \alpha_i - \beta_i p_i + \sum_{j \neq i} \gamma_{i,j} p_j }  \right] \label{eq:finitePcal_finiteUcal_semilog_biconjugate_precursor}
\end{align}
can be reformulated as a mixed-integer exponential cone program using the same biconjugate-based technique in Section~\ref{subsec:finitePcal_convexUcal_semilog}. In particular, when Assumption~\ref{assumption:Pcal_Cartesian_product} holds, then problem~\eqref{eq:finitePcal_finiteUcal_semilog_biconjugate_precursor} is equivalent to
\clearpage
\begin{subequations}
\begin{alignat}{2}
& \underset{\wb, \xb, \mub}{\text{maximize}} & \quad & \sum_{\ub \in \hat{\Ucal}} \sum_{i = 1}^I \mu_{\ub, i} \cdot 
\left( \log \lambda_{\ub} + \alpha_i \right) \nonumber \\
& & & + \sum_{\ub \in \hat{\Ucal}} \sum_{i = 1}^I\sum_{t \in \Pcal_i} \log t w_{\ub,i,i,t} \nonumber \\
& & & + \sum_{\ub \in \hat{\Ucal}} \sum_{i = 1}^I \sum_{t \in \Pcal_i} (- \beta_i) \cdot t \cdot w_{\ub,i,i,t} \nonumber \\
& & & + \sum_{\ub \in \hat{\Ucal}} \sum_{i = 1}^I \sum_{j \neq i} \gamma_{i,j} \cdot \sum_{t \in \Pcal_j} t \cdot w_{\ub,i,j,t} \nonumber \\
& & & - \sum_{\ub \in \hat{\Ucal}} \sum_{i=1}^I \mu_{\ub,i} \log \mu_{\ub,i} \\
& \text{subject to} & & \sum_{\ub \in \hat{\Ucal}} \sum_{i = 1}^I \mu_{\ub,i} = 1, \\
& & & \sum_{t \in \Pcal_j} w_{\ub,i,j,t} = \mu_{\ub,i}, \quad \forall \ \ub \in \hat{\Ucal}, \ i, j \in [I], \\
& & & \sum_{\ub \in \hat{\Ucal}} \sum_{i = 1}^I w_{\ub, i, j, t} = x_{j,t}, \quad \forall \ j \in [I], \ t \in \Pcal_j, \\
& & & \sum_{t \in \Pcal_j} x_{j,t} = 1, \quad \forall j \in [I], \\
& & & w_{\ub,i,j,t} \geq 0, \quad \forall \ \ub \in \hat{\Ucal}, \ i, j \in [I], t \in \Pcal_j, \\
& & & \mu_{\ub,i} \geq 0, \quad \forall \ \ub \in \hat{\Ucal}, \ i \in [I], \\
& & & x_{j,t} \in \{0,1\}, \quad \forall \ j \in [I], \ t \in \Pcal_j,
\end{alignat}
\label{prob:finitePcal_finiteUcal_semilog_dual_separation}%
\end{subequations}
where $x_{j,t}$ is a binary decision variable that is 1 if product $j$ is offered at price $t \in \Pcal_j$, and 0 otherwise; $\mu_{\ub,i}$ is a nonnegative decision variable introduced as part of the biconjugate-based reformulation; and $w_{\ub,i,j,t}$ is a decision variable that represents the linearization of $\mu_{\ub,i} \cdot x_{j,t}$ for all $\ub \in \hat{\Ucal}$, $i,j \in[I]$, and $t \in \Pcal_j$.

As with problem~\eqref{prob:semilog_MIECP}, this problem can be expressed as a mixed-integer exponential cone program. One notable difference between formulation~\eqref{prob:finitePcal_finiteUcal_semilog_dual_separation} and formulation~\eqref{prob:semilog_MIECP} from earlier is that the number of decision variables and constraints is larger because the decision variable $\mu_{\ub,i}$ is introduced for every combination of an uncertainty realization in $\hat{\Ucal}$ and each product $i$; thus, $\mub$ represents a probability mass function over the set $\hat{\Ucal} \times [I]$.

% )_{\ub \in \hat{\Ucal}, i \in [I]}$ is the

\subsection{Primal and dual subproblems for log-log demand model}
\label{subsec:finitePcal_finiteUcal_loglog}

For the log-log demand model, the primal separation problem is 
\begin{align}
\min_{\ub \in \Ucal} \sum_{\pb \in \hat{\Pcal}} \pi_{\pb} \cdot \sum_{i=1}^I p_i \cdot e^{ \alpha_i - \beta_i \log p_i + \sum_{j \neq i} \gamma_{i,j} \log p_j }.
\end{align}
Note that the objective function is convex in $\ub = (\alphab, \betab, \gammab)$; it is the nonnegative weighted combination of terms of the form $e^{ \alpha_i - \beta_i \log p_i + \sum_{j \neq i} \gamma_{i,j} \log p_j }$, each of which are convex in $(\alphab, \betab, \gammab)$. Thus, the overall problem, which can be stated as 
\begin{subequations}
\begin{alignat}{2}
& \underset{\ub, \zb}{\text{minimize}} & & \sum_{\pb \in \hat{\Pcal}} \pi_{\pb} \cdot \sum_{i=1}^I p_i \cdot e^{ \alpha_i - \beta_i \log p_i + \sum_{j \neq i} \gamma_{i,j} \log p_j } \\ 
& \text{subject to} & \quad & \ub = \Fb \zb, \\
& & &  \Ab \zb \leq \bb, \\
& & &  \zb \in \{0,1\}^n,
\end{alignat}
\end{subequations}
is a mixed-integer convex program, and can be expressed as a mixed-integer exponential cone program. 

The dual separation problem is 
\begin{align}
\max_{\pb \in \Pcal} \sum_{\ub \in \hat{\Ucal}} \lambda_{\ub} \cdot \left[ \sum_{i=1}^I p_i \cdot e^{ \alpha_i - \beta_i \log p_i + \sum_{j \neq i} \gamma_{i,j} \log p_j } \right]. 
\end{align}
The objective function of this problem is in general not concave in $\pb$. However, following the same method as in Section~\ref{subsec:finitePcal_convexUcal_loglog}, we can reformulate the related problem of maximizing the logarithm, which is 
\begin{align}
& \max_{\pb \in \Pcal} \log \left( \sum_{\ub \in \hat{\Ucal}} \lambda_{\ub} \cdot \left[ \sum_{i=1}^I p_i \cdot e^{ \alpha_i - \beta_i \log p_i + \sum_{j \neq i} \gamma_{i,j} \log p_j } \right]  \right) \\
& = \max_{\pb \in \Pcal} \log \left( \sum_{\ub \in \hat{\Ucal}} \sum_{i=1}^I e^{ \log \lambda_{\ub} + \log p_i + \alpha_i - \beta_i \log p_i + \sum_{j \neq i} \gamma_{i,j} \log p_j }  \right)
\end{align}
as a mixed-integer exponential cone program. Under Assumption~\ref{assumption:Pcal_Cartesian_product}, the resulting formulation is 
\begin{subequations}
\begin{alignat}{2}
& \underset{\xb, \wb, \mub}{\text{maximize}} & & \sum_{\ub \in \hat{\Ucal}} \sum_{i = 1}^I \mu_{\ub, i} \cdot 
\left( \log \lambda_{\ub} + \alpha_i \right) \nonumber \\
& & & + \sum_{\ub \in \hat{\Ucal}} \sum_{i = 1}^I\sum_{t \in \Pcal_i} (1 - \beta_i) \log t \cdot w_{\ub,i,i,t} \nonumber  \\
& & & + \sum_{\ub \in \hat{\Ucal}} \sum_{i = 1}^I \sum_{j \neq i} \gamma_{i,j} \cdot \sum_{t \in \Pcal_j} \log t \cdot w_{\ub,i,j,t} \nonumber  \\
& & & - \sum_{\ub \in \hat{\Ucal}} \sum_{i=1}^I \mu_{\ub,i} \log \mu_{\ub,i} \\
& \text{subject to} & \quad & \sum_{\ub \in \hat{\Ucal}} \sum_{i = 1}^I \mu_{\ub,i} = 1, \\
& & & \sum_{t \in \Pcal_j} w_{\ub,i,j,t} = \mu_{\ub,i}, \quad \forall \ \ub \in \hat{\Ucal}, \ i, j\in [I], \\
& & & \sum_{\ub \in \hat{\Ucal}} \sum_{i = 1}^I w_{\ub, i, j, t} = x_{j,t}, \quad \forall \ j \in [I], \ t \in \Pcal_j, \\
& & & \sum_{t \in \Pcal_j} x_{j,t} = 1, \quad \forall j \in [I], \\
& & & w_{\ub,i,j,t} \geq 0, \quad \forall \ \ub \in \hat{\Ucal}, \ i, j \in [I], t \in \Pcal_j, \\
& & & \mu_{\ub,i} \geq 0, \quad \forall \ \ub \in \hat{\Ucal}, \ i \in [I], \\
& & & x_{j,t} \in \{0,1\}, \quad \forall \ j \in [I], \ t \in \Pcal_j,
\end{alignat}
\label{prob:finitePcal_finiteUcal_loglog_dual_separation}
\end{subequations}
where the decision variables have the same meaning as those in formulation~\eqref{prob:finitePcal_finiteUcal_semilog_dual_separation}.

\subsection{Complexity of randomization for finite $\Pcal$, finite $\Ucal$}
\label{subsec:finitePcal_finiteUcal_complexity}

In the case that $\Pcal$ and $\Ucal$ are both finite, it is possible to characterize the support size more directly than in the case where $\Ucal$ is convex, by using standard results in linear programming. In particular, we have the following result. 

\begin{proposition}
Suppose that the uncertainty set $\Ucal$ is finite, and that $R(\pb, \ub) \geq 0$ for all $\pb \in \Pcal$ and $\ub \in \Ucal$. The RRPO problem $\max_{\pib \in \Delta_{\Pcal}} \min_{\ub \in \Ucal} \sum_{\pb \in \Pcal} \pi_{\pb} R(\pb, \ub)$ has an optimal solution $\pib^*$ with a support of size at most $|\Ucal| + 1$. \label{proposition:complexity_discreteUcal}
\end{proposition}

\begin{proofvvm}
The RRPO problem can be equivalently written as 
\begin{subequations}
\begin{alignat}{2}
& \underset{\sbee, t, \pib}{\text{maximize}} & \quad & t \\
& \text{subject to} & & t + s_{\ub} = \sum_{\pb \in \Pcal} \pi_{\pb} R(\pb, \ub), \quad \forall \ \ub \in \Ucal, \\
& & & \sum_{\pb \in \Pcal} \pi_{\pb} = 1, \\ 
& & & \pi_{\pb} \geq 0, \quad \forall\ \pb \in \Pcal, \\
& & & s_{\ub} \geq 0, \quad \forall\ \ub \in \Ucal, \\
& & & t \geq 0, 
\end{alignat}
\end{subequations}
where $s_{\ub}$ can be regarded as a slack variable for the constraint $t - \sum_{\pb \in \Pcal} \pi_{\pb} R(\pb, \ub) \leq 0$ (or equivalently, $t \leq \sum_{\pb \in \Pcal} \pi_{\pb} R(\pb, \ub)$). Observe that this is a linear program in standard form, with finitely many equality constraints and finitely many variables. Note that the feasible region is non-empty, because given any $\pib \in \Delta_{\Pcal}$, we can select $t$ and $\sbee = (s_{\ub})_{\ub \in \Ucal}$ to obtain a feasible solution to the LP above; note that additionally the problem is bounded, because $\{ R(\pb, \ub) \mid \pb \in \Pcal, \ub \in \Ucal\}$ is a finite set. Hence, it follows that there exists an optimal solution. Since the problem is in standard form, there is at least one extreme point (Corollary 2.2 of \citealt{bertsimas1997introduction}). By Theorem 2.7 of \cite{bertsimas1997introduction}, there must exist an extreme point solution that is optimal. Since extreme points correspond to basic feasible solutions, this optimal solution must be such that at most $|\Ucal|+1$ coordinates of $(\sbee, t, \pib)$ are non-zero. For the resulting distribution $\pib$, at most $|\Ucal|+1$ prices therefore have non-zero probability. \Halmos 
\end{proofvvm}

%{\color{blue}
%\section{Constraining the price distribution for finite $\Pcal$, convex $\Ucal$}
%
%
%
%\comebacktothis
%
%}

\section{Additional numerical results}

%%%%
% START OF DISCRETE UCAL DISCRETE PCAL SYNTHETIC EXPERIMENTS
%%%%

\subsection{Experiments with discrete $\Ucal$ and log-log and semi-log demand models} 
\label{subsec:results_discreteUcal_loglog_semilog}

In this set of experiments, we consider linear, log-log and semi-log demand models, where uncertainty is modeled through a discrete uncertainty set. We specifically consider a discrete budget uncertainty set $\Ucal$ here:
\begin{equation}
\Ucal= \{\ub = \ub_0 - (\ub_0-\bar{\ub})\circ\xib - (\ub_0-\underline{\ub})\circ\etab \mid \eb^\top \xib+\eb^\top\etab \leq \Gamma, \xib+\etab \leq \oneb, \xib,\etab\in\{0,1\}^{I+I^2}\}, \label{eq:discrete_budget_Ucal}
\end{equation}
where $\underline{\ub}$ and $\bar{\ub}$ are respectively the component-wise lower and upper bounds of $\ub$, $\ub_0$ is the nominal value of the uncertain parameter vector $\ub = (\alphab, \betab,\gammab)$, $\circ$ denotes the element-wise product of two vectors, $\Gamma$ is the budget of uncertainty and $I + I + I(I-1) = I + I^2$ is the total number of demand model parameters. Under the budget uncertainty set $\Ucal$, up to $\Gamma$ parameters can attain their lower bounds or upper bounds, whereas the remaining parameters can only attain their nominal values. We shall assume that the lower bound vector $\underline{\ub}$ and upper bound vector $\overline{\ub}
$ are defined as $\underline{\ub}=0.7\ub_0$ and $\overline{\ub}=1.3\ub_0$, where $\ub_0$ is the vector of nominal parameters.

For each of the three demand models (linear, semi-log and log-log), we vary the number of products $I$ in $\{5, 10, 15\}$. For each value of $I$, we generate 24 random instances, where the values of $\alphab, \betab, \gammab$ are independently randomly generated as follows:
\begin{enumerate}
\item \emph{Linear demand}. Each $\alpha_i \sim \Uniform(100, 200)$, $\beta_i \sim \Uniform(5,15)$, $\gamma_{i,j} \sim \Uniform(-0.1,+0.1)$.
\item \emph{Semi-log demand}. Each $\alpha_i \sim \Uniform(8,10)$, $\beta_i \sim \Uniform(1.5,2)$, $\gamma_{i,j} \sim \Uniform(-0.5, +0.5)$.
\item \emph{Log-log demand}. Each $\alpha_i \sim \Uniform(10,14)$, $\beta_i \sim \Uniform(1.5,2)$, $\gamma_{i,j} \sim \Uniform(-0.8, +0.8)$. 
\end{enumerate}
We set the price set of each $i \in [I]$ as $\Pcal_i = \{1,2,3,4,5\}$. 

For each instance, we solve the nominal problem, the DRPO problem and the RRPO problem. For both RRPO and DRPO, we test a different collection of $\Gamma$ values for the uncertainty set depending on the value of $I$. 

To solve the RRPO problem for each instance, we execute the double column generation algorithm described in Section~\ref{sec:finitePcal_finiteUcal}. For the linear demand model, we solve both primal and dual separation problems as mixed-integer programs in Gurobi.  

For the log-log and semi-log demand models, we proceed differently. In our preliminary experimentation with the restricted dual problem, we observed that exactly solving the dual separation problem \eqref{prob:finitePcal_finiteUcal_semilog_dual_separation} (for semi-log demand) or \eqref{prob:finitePcal_finiteUcal_loglog_dual_separation} (for log-log demand) via Mosek takes quite a long time. Therefore, to reduce the computation time of RRPO with discrete $\Ucal$, we instead use a random improvement heuristic to obtain the solution of dual separation problem. Specifically, we randomly select a price vector $\pb^0$ as a starting point. We start with changing the price of product $i=1$ and keeping the prices of all other products unchanged, to search for a price vector $\pb^1$ that makes the objective value of the dual separation problem the largest. Then based on the current price vector $\pb^1$, we change the price of product  $i=2$ and keep the prices of all other products unchanged, to search for a better price vector $\pb^2$. We repeat this for all of the products, yielding the price vector $\pb^{I}$. We repeat this procedure with 100 random starting points, and retain the best solution over these 100 repetitions. We note that in the heuristic form of the double column generation method, we continue to solve the primal separation problem (over $\ub \in \Ucal$) exactly using Mosek. The lower bound given by the restricted primal, $Z_{P, \hat{\Pcal}}$, therefore remains a valid lower bound on the true optimal RRPO objective, and will correspond to the worst-case expected revenue of a possibly suboptimal price distribution $\pib$, at every iteration of the double column generation method. At termination, the objective value returned by the algorithm will therefore be a lower bound on the true optimal RRPO objective.

Although this approximate method cannot guarantee that the overall double column generation procedure converges to a provably optimal solution, our preliminary experimentation with small instances suggests that it obtains the exact solution of RRPO that one would obtain if the dual separation problem were solved to provable optimality. In Section~\ref{subsec:extra_results_exact_vs_heuristic} of the ecompanion, we undertake an in-depth comparison of this heuristic form of double column generation against the exact double column generation approach, where both separation problems are solved to provable optimality using Mosek. When one imposes a time limit on the exact approach, we find that the heuristic approach returns solutions of comparable quality (in cases where the exact approach terminates with the optimal solution), or better (in cases where the exact approach exhausts the time limit and does not return an optimal solution), and in general, the heuristic approach terminates in significantly less time than the exact approach.

With regard to the DRPO problem for each log-log and semi-log instance, we note that we do not have a solution algorithm or formulation to solve it exactly. Therefore, we again use the same random improvement heuristic to obtain an approximate solution of DRPO with these demand models. We randomly pick a starting price vector, and change the price of one product at a time to improve the worst case objective value until we no longer get an improvement. We repeat this procedure 50 times and select the best resulting price vector from these 50 repetitions as the approximate solution of DRPO. We note that we use a smaller number of repetitions because each repetition involves solving worst-case problem over $\ub \in \Ucal$ repeatedly in order to evaluate the robust objective of each candidate price vector; this contributes to a large overall computation time for this approach. We denote this approximate solution by $\hat{\pb}_{\DR}$. With regard to the DRPO problem for linear demand, we observe that the objective function of DRPO is linear in the uncertain parameter vector $\ub$, and that the description of the set polyhedron~\eqref{eq:discrete_budget_Ucal} is integral (i.e., extreme points of this polyhedron naturally correspond to $\xib, \etab \in \{0,1\}^{2I + I^2}$). Therefore, DRPO can be solved exactly by relaxing the requirement $\xib, \etab \in \{0,1\}^{2I + I^2}$ in the uncertainty set~\eqref{eq:discrete_budget_Ucal}, and reformulating the worst-case objective using LP duality, leading to a mixed-integer linear program. 

Lastly, for the nominal problem for each instance, we either formulate it as a mixed-integer linear program (in the linear demand case) or use the biconjugate technique to formulate it as a mixed-integer exponential cone program (in the log-log and semi-log demand cases).

We report the same metrics as in Section~\ref{subsec:results_convexUcal}, with several modifications for the log-log and semi-log demand models. For those models, we use $\hat{Z}_{\DR}$ and $R(\hat{\pb}_{\DR}, \ub_0)$ to denote the approximate objective value of the DRPO problem and revenue of the approximate DRPO solution $\hat{\pb}_{\DR}$ under the nominal parameter vector $\ub_0$. We also use $\hat{Z}_{\RR}$ to denote the approximate objective of the RRPO problem, and $\Ebb[ R(\hat{\pb}_{\RR}, \ub_0)]$ denotes the expected revenue of the resulting randomized price vector $\hat{\pb}_{\RR}$ under the nominal demand model parameter vector $\ub_0$. Finally, the approximate relative improvement of RRPO over DRPO is defined as $\hat{\RI} = (\hat{Z}_{\RR} - \hat{Z}_{\DR} )/\hat{Z}_{\DR} \times 100\%$.

Table~\ref{table:result_discreteU_linear} shows the results for the linear demand model. Here, we interestingly find that the vast majority of instances are randomization-proof, i.e., the average $\RI$ is below 1\%, if not exactly 0\%. We note here that we tested other families of instances where $(\alphab, \betab, \gammab)$ and $\Pcal_i$ are generated differently, but in virtually every case we found that the relative improvement of randomized over deterministic robust pricing was very small. These results, together with those for the convex $\Ucal$ case, suggest that randomized pricing is of limited benefit compared to deterministic pricing for the uncertain linear demand model case. % \comebacktothis

Tables~\ref{table:result_discreteU_semi_log} and \ref{table:result_discreteU_log_log} show how the results vary for different values of discrete uncertainty budget $\Gamma$ for semi-log and log-log demand.\footnote{We note here that for the log-log model, we encountered one instance ($I = 10$, $\Gamma = 44$) where $\hat{Z}_{\DR}$ was higher than $\hat{Z}_{\RR}$; in general, when solved to perfect optimality, one should see $Z^*_{\RR}$ should be higher than $Z^*_{\DR}$. We have verified that the reason for this anomaly was a numerical error in the solution of the worst-case subproblem in Mosek within the DRPO random improvement heuristic. This instance is omitted in our calculation of $\hat{Z}_{\DR}$, $\hat{\RI}$ and $\Ebb[R(\hat{\pb}_{\DR}, \ub_0)]$, and the affected entries are indicated by * in Table~\ref{table:result_discreteU_log_log}. }  We can see that, in most of the cases we test, the randomized robust pricing strategy provides a substantial benefit over the deterministic robust price solution. The percentage improvement given by randomization ranges from $0\%$ to as much as $488.59\%$ for semi-log instances, and from $0\%$ to $175.18\%$ for log-log instances. Similar to the cases with convex $\Ucal$, both $\hat{Z}_{\RR}$ and $ \hat{Z}_{\DR}$ decrease as the uncertainty set becomes larger. While the RI metric generally decreases as $\Gamma$ increases, in some instances it can be increasing in $\Gamma$ at small values of $\Gamma$ (this is visible in the average results metrics for $I = 10$ with semi-log demand). When $\Gamma$ is large enough, the RI metric often becomes very small or even zero. This makes sense when interpreted through Corollary~\ref{corollary:finite_Pcal_pDR}. Specifically, when nature is able to make a large number of demand model parameters take their worst values, it is likely that at the $\ub^*$ at which the optimal objective of DRPO is attained is such that the price vector for the nominal problem with $\ub^*$ coincides with the optimal price vector for DRPO. %(As an extreme case, suppose that $\Gamma$ were set to $I + I^2$; because the prices are nonnegative, it is optimal for the worst case problem to set $\etab = \oneb$ and $\xib = \zerob$, and it is clear that the optimal price vector for DRPO will be the same as the nominal problem with $\underline{\ub}$ as the vector of demand model parameters.) 
Thus, by Corollary~\ref{corollary:finite_Pcal_pDR}, the problem will be randomization-proof.   %Besides, by the comparison between $\Ebb[R(\pb^*_{\RR},\ub_0)]$ and $R(\pb^*_{\DR},\ub_0)$, we can see that the nominal revenue of RRPO solution is larger than that of DRPO solution in nearly $1/3$ of the instances.

With regard to the computation time, the computation time for both RRPO and DRPO increases with $I$. Interestingly, the computation time required by RRPO does not necessarily increase as the discrete uncertainty budget $\Gamma$ increases; in some cases, when $\Gamma$ is large, the RRPO solution degenerates to the DRPO solution, allowing the double column generation algorithm to terminate quickly. By comparing $t_{\RR}$ and $t_{\DR}$, we can see that RRPO in general takes less time than DRPO. The computation time of the RRPO problem in semi-log instances is no more than approximately two minutes on average ($I=15$, $\Gamma=36$), while in log-log instances, solving RRPO requires no more than 1.5 minutes on average   ($I=15$, $\Gamma=18$). Lastly, for linear demand, the computation time for RRPO is extremely small, requiring no more than a few seconds on average. 

% We cannot compare the performance among different $I$ because these are totally different instances. The scale of revenue does not necessaily increase with $I$.

\begin{table}[ht]
	\centering
	\small
	\begin{tabular}{llllllllllllll}
		\toprule
		$I$&$\Gamma$ & $t_{\RR}$ & $Z_{\RR}^*$ & $\Ebb[R(\pb^*_{\RR},\ub_0)]$ & $t_{\DR}$ & $Z^*_{\DR}$  & RI(\%) & $R(\pb^*_{\DR},\ub_0)$ & $t_{\Nom}$ & $Z_{\Nom}^*$ & $Z_{\Nom,\WC}$ \\ 
		\midrule
		5 &   3 & 8.76 & 1724.38 & 2458.17 & 0.43 & 1723.46  & 0.06 & 2462.66 & 0.32 & 2473.30 & 1719.29 \\     
		5 &   6 & 0.27 & 1316.39 & 2382.10 & 0.04 & 1313.68  & 0.23 & 2383.35 & -- & --  & 1261.84 \\     
		5 &   9 & 0.27 & 1159.33 & 2294.44 & 0.03 & 1157.98  & 0.12 & 2293.84 & --& -- & 1039.24 \\     
		5 &  12 & 0.19 & 1126.85 & 2259.05 & 0.03 & 1126.85  & 0.00 & 2259.05 & -- & -- & 987.29 \\     
		5 &  18 & 0.20 & 1124.87 & 2256.23 & 0.03 & 1124.87  & 0.00 & 2256.23 & -- & --  & 984.10 \\     
		5 &  24 & 0.17 & 1123.78 & 2256.23 & 0.03 & 1123.78  & 0.00 & 2256.23 & -- & --  & 982.18 \\ 
		\midrule
		10 &   5 & 0.53 & 3717.70 & 4990.89 & 0.11 & 3714.08  & 0.10 & 4997.68 & 0.09 & 5009.03 & 3711.50 \\    
		10 &   7 & 0.59 & 3315.59 & 4972.42 & 0.11 & 3312.91  & 0.08 & 4973.49 & -- & -- & 3297.07 \\    
		10 &   9 & 0.70 & 2982.02 & 4941.26 & 0.11 & 2979.68  & 0.08 & 4944.40 &-- & -- & 2942.80 \\    
		10 &  14 & 1.35 & 2589.70 & 4782.85 & 0.12 & 2583.93  & 0.23 & 4791.83 & -- & --  & 2437.09 \\    
		10 &  19 & 0.84 & 2372.37 & 4662.16 & 0.11 & 2371.21  & 0.05 & 4657.17 & --& --  & 2133.38 \\    
		10 &  26 & 0.56 & 2342.57 & 4636.61 & 0.10 & 2342.57  & 0.00 & 4636.61 & -- & --  & 2086.03 \\    
		10 &  33 & 0.57 & 2339.18 & 4633.87 & 0.11 & 2339.18  & 0.00 & 4636.61 & -- & -- & 2081.46 \\    
		10 &  44 & 0.57 & 2334.90 & 4627.78 & 0.10 & 2334.90  &0.00 & 4627.78 &-- & -- & 2075.15 \\
		\midrule
		15 &   6 & 0.92 & 5920.62 & 7495.77 & 0.22 & 5918.51  & 0.04 & 7506.26 & 0.17 & 7513.95  & 5917.60 \\    
		15 &  12 & 1.79 & 4706.04 & 7435.56 & 0.25 & 4701.16  & 0.11 & 7442.32 & --& --  & 4668.21 \\    
		15 &  18 & 3.16 & 4051.33 & 7250.57 & 0.27 & 4045.35  & 0.15 & 7248.00 & -- & --  & 3889.25 \\    
		15 &  24 & 4.50 & 3722.40 & 7071.15 & 0.29 & 3713.86  & 0.23 & 7068.49 & --& --  & 3428.59 \\    
		15 &  36 & 1.12 & 3500.73 & 6889.94 & 0.21 & 3500.73  & 0.00 & 6889.47 & -- & --  & 3112.83 \\    
		%15 &  48 & 1.20 & 3494.32 & 6886.66 & 0.21 & 3494.32  & 0.00& 6889.47 & -- & -- & 3104.39 \\    
		%15 &  60 & 1.16 & 3488.86 & 6882.94 & 0.22 & 3488.86  & 0.00& 6882.94 & --& --  & 3096.57 \\    
		%15 &  72 & 1.24 & 3484.07 & 6880.49 & 0.22 & 3484.06  & 0.00& 6880.38 & -- & --  & 3089.30 \\ 
		\bottomrule
	\end{tabular}
	\caption{Results for linear instances with discrete $\Ucal$. \label{table:result_discreteU_linear}}
\end{table}

\begin{table}[ht]
	\centering
	\small
		\begin{tabular}{llllllllllll}
			\toprule
			$I$&$\Gamma$ & $t_{\RR}$ & $\hat{Z}_{\RR}$ & $\Ebb[R(\hat{\pb}_{\RR},\ub_0)]$ & $t_{\DR}$ & $\hat{Z}_{\DR}$ & $\hat{\RI}$ (\%) & $R(\hat{\pb}_{\DR},\ub_0)$ & $t_{\Nom}$ & $Z_{\Nom}^*$ & $Z_{\Nom,\WC}$ \\ 
			\midrule
			5 &   3 & 13.23 & \numvvm{8649.86} & \numvvm{157098.87} & 20.35 & \numvvm{5561.00} & 63.01 & \numvvm{216781.36} & 0.93 & \numvvm{227228.32} & \numvvm{5040.46} \\ 
    		5 &   6 & 0.43 & \numvvm{3060.07} & \numvvm{181189.81} & 20.52 & \numvvm{1918.06} & 66.89 & \numvvm{219660.30} & -- & -- & \numvvm{1776.72} \\ 
    		5 &   9 & 0.55 & \numvvm{1840.13} & \numvvm{208530.96} & 22.96 & \numvvm{1639.60} & 18.42 & \numvvm{222140.73} & -- &-- & \numvvm{1625.07} \\ 
    		5 &  12 & 0.63 & \numvvm{1647.64} & \numvvm{221571.95} & 25.53 & \numvvm{1599.35} & 4.39 & \numvvm{224743.54} & -- & -- & \numvvm{1599.51} \\ 
    		5 &  18 & 0.27 & \numvvm{1590.56} & \numvvm{227101.80} & 23.87 & \numvvm{1589.58} & 0.11 & \numvvm{227228.15} & -- & -- & \numvvm{1589.58} \\ 
    		5 &  24 & 0.08 & \numvvm{1588.80} & \numvvm{227227.34} & 13.53 & \numvvm{1588.80} & 0.00 & \numvvm{227228.15} & -- & -- & \numvvm{1588.80} \\ 
			\midrule
			10 &   5 & 1.59 & \numvvm{1207444.93} & \numvvm{53046716.46} & 117.69 & \numvvm{494733.98} & 173.47 & \numvvm{67299115.62} & 0.15 & \numvvm{71616384.73} & \numvvm{447565.73} \\ 
   			10 &   7 & 2.59 & \numvvm{649903.90} & \numvvm{41041125.23} & 119.28 & \numvvm{186678.65} & 233.11 & \numvvm{56898320.22} & -- & -- & \numvvm{188458.06} \\ 
   			10 &   9 & 3.49 & \numvvm{432040.81} & \numvvm{36404899.11} & 121.43 & \numvvm{108902.10} & 250.12 & \numvvm{54074501.31} & -- & -- & \numvvm{104687.21} \\ 
   			10 &  14 & 6.44 & \numvvm{180180.92} & \numvvm{36662587.82} & 180.63 & \numvvm{66896.02} & 161.08 & \numvvm{69555655.70} & -- & -- & \numvvm{67766.09} \\ 
   			10 &  19 & 10.23 & \numvvm{96340.91} & \numvvm{52870810.71} & 362.29 & \numvvm{63587.96} & 67.80 & \numvvm{71477486.58} & -- & -- & \numvvm{63625.57} \\ 
  			 10 &  26 & 12.14 & \numvvm{66891.41} & \numvvm{64013288.45} & 559.63 & \numvvm{63244.97} & 19.70 & \numvvm{71477486.58} & -- & -- & \numvvm{63329.47} \\ 
   			10 &  33 & 11.37 & \numvvm{63736.41} & \numvvm{71112478.62} & 699.06 & \numvvm{63287.93} & 6.47 & \numvvm{71540861.08} & -- & -- & \numvvm{63291.95} \\ 
   			10 &  44 & 7.84 & \numvvm{63306.22} & \numvvm{71494123.33} & 799.81 & \numvvm{63299.46} & 0.09 & \numvvm{71557767.25} & -- & -- & \numvvm{63278.65} \\
			\midrule
			15 &   6 & 37.63 & \numvvm{355373227.60} & \numvvm{3220896064.00} & 331.13 & \numvvm{15852841.99} & 488.59 & \numvvm{4518485256.00} & 0.44 & \numvvm{5053760085.00} & \numvvm{14200239.04} \\ 
   			15 &  12 & 18.86 & \numvvm{8685177.19} & \numvvm{3236504376.00} & 344.94 & \numvvm{2151155.02} & 471.14 & \numvvm{3509674763.00} & -- & -- & \numvvm{1773100.56} \\ 
   			15 &  18 & 33.08 & \numvvm{3106034.33} & \numvvm{3079453445.00} & 653.64 & \numvvm{1177177.17} & 323.61 & \numvvm{4913148447.00} & -- & -- & \numvvm{1155687.93} \\ 
   			15 &  24 & 48.52 & \numvvm{1748546.17} & \numvvm{3671378650.00} & 1525.41 & \numvvm{1115540.50} & 165.75 & \numvvm{5032597478.00} & -- & -- & \numvvm{1110257.31} \\ 
   			15 &  36 & 101.38 & \numvvm{1211100.07} & \numvvm{4819415642.00} & 3400.17 & \numvvm{1099215.62} & 38.58 & \numvvm{5024624319.00} & -- & -- & \numvvm{1100760.78} \\ 
   			%15 &  48 & 129.43 & \numvvm{1116581.47} & \numvvm{4992959039.00} & 5150.42 & \numvvm{1100761.52} & 6.99 & \numvvm{5036957282.00} & -- & -- & \numvvm{1099800.15} \\ 
   			%15 &  60 & 136.72 & \numvvm{1102953.48} & \numvvm{5036076696.00} & 6921.49 & \numvvm{1101344.94} & 0.64 & \numvvm{5040643673.00} & -- & -- & \numvvm{1099612.99} \\ 
			\bottomrule
		\end{tabular}
		\caption{Results for semi-log instances with discrete $\Ucal$. \label{table:result_discreteU_semi_log}}
\end{table}

%%%%%%%%%%%%%%%%%%%%%%%%%%%%%%%

\begin{table}[ht]
	\centering
	\small
		\begin{tabular}{llllllllllll}
			\toprule
			$I$&$\Gamma$ & $t_{\RR}$ & $\hat{Z}_{\RR}$ & $\Ebb[R(\hat{\pb}_{\RR},\ub_0)]$ & $t_{\DR}$ & $\hat{Z}_{\DR}$ & $\hat{\RI}$ (\%) & $R(\hat{\pb}_{\DR},\ub_0)$ & $t_{\Nom}$ & $Z_{\Nom}^*$ & $Z_{\Nom,\WC}$ \\ 
			\midrule
			5 &   3 & 13.33 & \numvvm{326068.23} & \numvvm{1857192.29} & 31.23 & \numvvm{279630.96} & 22.03 & \numvvm{1889663.60} & 0.96 & \numvvm{4306037.70} & \numvvm{149373.68} \\ 
    		5 &   6 & 0.20 & \numvvm{63709.13} & \numvvm{4026167.55} & 20.39 & \numvvm{61461.82} & 3.16 & \numvvm{4288399.24} & -- & -- & \numvvm{60999.60} \\ 
    		5 &   9 & 0.27 & \numvvm{49063.48} & \numvvm{4194592.91} & 18.39 & \numvvm{48739.47} & 0.65 & \numvvm{4236343.37} & -- & -- & \numvvm{47720.72} \\ 
    		5 &  12 & 0.35 & \numvvm{47095.58} & \numvvm{4202272.77} & 16.82 & \numvvm{47001.63} & 0.20 & \numvvm{4181090.80} & -- & -- & \numvvm{45705.83} \\ 
    		5 &  18 & 0.17 & \numvvm{46630.10} & \numvvm{4180420.21} & 11.67 & \numvvm{46628.29} & 0.01 & \numvvm{4181090.80} & -- & -- & \numvvm{45201.86} \\ 
    		5 &  24 & 0.14 & \numvvm{46611.38} & \numvvm{4181090.80} & 8.73 & \numvvm{46587.19} & 0.05 & \numvvm{4162529.42} & --& -- & \numvvm{45180.77} \\ 
			\midrule
			10 &   5 & 4.64 & \numvvm{2330359.56} & \numvvm{38902505.73} & 161.57 & \numvvm{1382127.68} & 83.54 & \numvvm{64513223.67} & 0.31 & \numvvm{76582442.47} & \numvvm{1228690.24} \\ 
   			10 &   7 & 5.39 & \numvvm{1306883.08} & \numvvm{52187976.23} & 165.21 & \numvvm{858800.25} & 63.34 & \numvvm{71234671.17} & -- & -- & \numvvm{838521.55} \\ 
   			10 &   9 & 5.30 & \numvvm{873913.12} & \numvvm{51527766.65} & 149.63 & \numvvm{623158.47} & 38.99 & \numvvm{72789417.40} & -- & --& \numvvm{612204.88} \\ 
   			10 &  14 & 5.10 & \numvvm{506916.01} & \numvvm{56673465.14} & 138.67 & \numvvm{382227.14} & 29.31 & \numvvm{73019886.42} & -- & -- & \numvvm{380471.28} \\ 
   			10 &  19 & 3.52 & \numvvm{373545.98} & \numvvm{63520067.90} & 136.70 & \numvvm{330830.96} & 12.90 & \numvvm{75148858.89} & -- & -- & \numvvm{333475.96} \\ 
   			10 &  26 & 3.56 & \numvvm{326822.06} & \numvvm{73754394.91} & 138.77 & \numvvm{318495.04} & 4.28 & \numvvm{75885562.90} & -- & -- & \numvvm{320621.90} \\ 
   			10 &  33 & 4.79 & \numvvm{317753.17} & \numvvm{75977959.88} & 137.07 & \numvvm{314570.87} & 2.34 & \numvvm{75885562.90} & -- & -- & \numvvm{316675.64} \\ 
   			10 &  44 & 3.19 & \numvvm{314631.83} & \numvvm{76242567.78} & 131.47 & \numvvm{319396.86}* & {1.19}* & \numvvm{75349501.72}* & -- & -- & \numvvm{314382.86} \\ 
			\midrule
			15 &   6 & 35.68 & \numvvm{17110848.95} & \numvvm{432684796.90} & 552.64 & \numvvm{7752987.37} & 175.18 & \numvvm{783147066.70} & 1.20 & \numvvm{838470429.10} & \numvvm{7204471.98} \\ 
   			15 &  12 & 62.70 & \numvvm{5545220.55} & \numvvm{492889975.10} & 582.20 & \numvvm{2939628.98} & 102.06 & \numvvm{777587471.30} & -- & -- & \numvvm{2764547.45} \\ 
  			15 &  18 & 76.60 & \numvvm{2924465.88} & \numvvm{455426557.50} & 561.96 & \numvvm{1923908.40} & 65.73 & \numvvm{823588835.90} & -- & -- & \numvvm{1879023.14} \\ 
   			15 &  24 & 54.30 & \numvvm{2027780.73} & \numvvm{679387985.70} & 572.08 & \numvvm{1705034.31} & 32.76 & \numvvm{819187783.30} & -- & -- & \numvvm{1707832.46} \\ 
   			15 &  36 & 48.34 & \numvvm{1692310.31} & \numvvm{810775120.30} & 600.30 & \numvvm{1630492.16} & 8.07 & \numvvm{831592791.90} & -- & -- & \numvvm{1634824.62} \\   
			\bottomrule
		\end{tabular}
		\caption{Results for log-log instances with discrete $\Ucal$.  \label{table:result_discreteU_log_log}}
\end{table}

%%%%
% END OF ORIGINAL DISCRETE UCAL DISCRETE PCAL SYNTHETIC EXPERIMENTS
%%%%

%%%%
% START OF EXACT VS HEURISTIC RESULTS
%%%%

\subsection{Comparison of heuristic and exact solutions of RRPO in the finite $\Pcal$, finite $\Ucal$ case}
\label{subsec:extra_results_exact_vs_heuristic}

%\comebacktothis -- is notation for LBs, UBs consistent with earlier? 

In this section, we perform a detailed comparison of the performance of the heuristic form of RRPO, where the dual separation problem is solved approximately using the random improvement heuristic, against the exact form of RRPO, where both the primal and dual separation problem are solved exactly. We consider the log-log and semi-log demand models, and we test the same instances described in Section~\ref{subsec:results_discreteUcal_loglog_semilog}. 

For the exact approach, we incorporated the following termination logic. When the relative gap between the lower bound from the restricted primal and the upper bound from the restricted dual is within $10^{-6}$. In certain cases, due to numerical precision issues arising from computing the log-log and semi-log objective functions and solving the primal and dual separation problems within Mosek, we would alternatively terminate when the set of price vectors $\hat{\Pcal}$ no longer changes (i.e., when the dual separation problem is unable to produce a price vector that is not contained in $\hat{\Pcal}$). Finally, we also imposed a time limit of 3600 seconds for the double column generation algorithm, i.e., we exit from the loop in Algorithm~\ref{alg:doubleCG} if the time limit exceeds 3600 seconds. We note here that the first two termination criteria were also implemented for the heuristic approach. 

In Tables~\ref{table:exact_vs_heuristic_loglog} and \ref{table:exact_vs_heuristic_semilog}, we display the optimality gap $G_E$ for the exact approach, the computation time $t_{E}$ for the exact approach, the computation time $t_H$ for the heuristic approach, and the relative improvement $RI_{E \to H}$ of the heuristic approach over the exact approach, for the log-log and semi-log demand models respectively. The optimality gap $G_E$ is defined as 
\begin{equation}
G_E = 100\% \times \frac{ Z^*_{\RR,E,UB} - Z^*_{\RR,E,LB}}{Z^*_{\RR,E,UB}},
\end{equation}
where $Z^*_{\RR,E,UB}$ and $Z^*_{\RR,E,LB}$ are the final upper and lower bounds from the double column generation approach. The relative improvement $RI_{E \to H}$ is defined as
\begin{equation}
\RI_{E \to H} = 100\% \times \frac{Z^*_{\RR,H,LB} - Z^*_{\RR,E,LB}}{Z^*_{\RR,E,LB}},
\end{equation}
where $Z^*_{\RR,H,LB}$ is the final lower bound from the heuristic approach. The relative improvement metric measures how much the price distribution from the heuristic double CG method improves on the price distribution from the exact double CG method. Note that this metric can be negative, in cases where the heuristic approach terminates with a suboptimal solution to the RRPO problem. 

Before continuing, we note that in some instances, we found the gap $G_E$ to be either slightly larger than the tolerance of $10^{-6}$, or otherwise to be a very small negative number; both of these phenomena are due to the aforementioned numerical issues we encountered. We also note that for larger instances, the computation time sometimes significantly exceeds 3600 seconds. As noted before, the outer loop of the double CG algorithm is terminated when the run time exceeds 3600 seconds, and in some cases, the final iteration of the outer loop before this criterion is triggered can require a huge amount of time (due to solving both the restricted primal and restricted dual problems, each of which can involve multiple solves of the primal and dual separation problems). 

For the log-log results shown in Table~\ref{table:exact_vs_heuristic_loglog}, we see that in general, the exact double CG algorithm terminates with a provably optimal solution, as the average optimality gap is very close to zero across all values of $I$ and $\Gamma$. With regard to the heuristic double CG method, we can see the average RI metric is in general zero or no smaller than -0.36\%. This implies that the heuristic double CG method returns solutions that are in general either the same or only very slightly suboptimal compared to the exact double CG solution. We additionally see that the average computation time of the heuristic double CG method is much smaller than that of the exact double CG method. For example, for $I = 15$, the average computation time of the exact double CG method can range anywhere from around 15 minutes to about 45 minutes, whereas for the heuristic double CG method it is no more than 1.5 minutes. These results indicate that for the log-log model, the heuristic method can obtain solutions of essentially the same quality as the exact method but in substantially less time.

For the semi-log results shown in Table~\ref{table:exact_vs_heuristic_semilog}, we see similar behavior as we do for the log-log models for $I \in \{5, 10\}$. For $I = 15$, we find that the exact double CG method in general does not terminate with a provably optimal solution (i.e., a gap of 0\%) under the soft 3600 second time limit. The average optimality gap for the exact method ranges from about 14\% to 40\%, and the average computation time can be as high as 90 minutes. In contrast, the average computation time for the heuristic double CG method is no more than two minutes. In addition, for $\Gamma \geq 12$, the heuristic solution has an objective value for which the average improvement over the (suboptimal) exact double CG solution ranges from approximately 6\% to 20\%. This highlights that in settings where the exact double CG method requires a large amount of computation time and must be terminated prematurely, the heuristic method obtains higher quality solutions than the exact method.

\begin{table}[ht]
\centering
\begin{tabular}{rrrrrr} \toprule
$I$ & $\Gamma$ & $G_{E}$ & $t_E$ & $t_H$ & $\RI_{E \to H}$ \\ \midrule
    5 &   3 & 0.00 & 16.16 & 13.33 & 0.00 \\ 
    5 &   6 & 0.00 & 1.58 & 0.20 & 0.00 \\ 
    5 &   9 & -0.34 & 2.21 & 0.27 & -0.35 \\ 
    5 &  12 & -0.06 & 2.45 & 0.35 & -0.07 \\ 
    5 &  18 & -0.03 & 1.88 & 0.17 & -0.05 \\ 
    5 &  24 & 0.00 & 1.49 & 0.14 & 0.00 \\ \midrule
   10 &   5 & 0.00 & 63.35 & 4.64 & -0.12 \\ 
   10 &   7 & 0.01 & 77.43 & 5.39 & 0.00 \\ 
   10 &   9 & 0.01 & 80.18 & 5.30 & 0.00 \\ 
   10 &  14 & 0.01 & 83.99 & 5.10 & -0.01 \\ 
   10 &  19 & 0.01 & 72.05 & 3.52 & -0.03 \\ 
   10 &  26 & 0.01 & 74.75 & 3.56 & 0.00 \\ 
   10 &  33 & 0.01 & 70.30 & 4.79 & 0.00 \\ 
   10 &  44 & 0.01 & 68.23 & 3.19 & 0.00 \\ \midrule
   15 &   6 & 0.01 & 952.93 & 35.68 & -0.25 \\ 
   15 &  12 & 0.02 & 2356.89 & 62.70 & -0.30 \\ 
   15 &  18 & 0.23 & 2739.09 & 76.60 & -0.36 \\ 
   15 &  24 & 0.04 & 2033.25 & 54.30 & -0.24 \\ 
   15 &  36 & 0.04 & 1351.31 & 48.34 & -0.02 \\ \bottomrule
\end{tabular}
\caption{Comparison of exact double CG and heuristic double CG methods for log-log demand model instances from Section~\ref{subsec:results_discreteUcal_loglog_semilog}. \label{table:exact_vs_heuristic_loglog}} 
\end{table}

\begin{table}[ht]
\centering
\begin{tabular}{rrrrrr}
  \toprule
$I$ & $\Gamma$ & $G_{E}$ & $t_E$ & $t_H$ & $\RI_{E \to H}$ \\ \midrule
  5 &   3 & 0.00 & 18.30 & 13.23 & 0.00 \\ 
    5 &   6 & 0.00 & 7.57 & 0.43 & 0.00 \\ 
    5 &   9 & 0.00 & 8.52 & 0.55 & 0.00 \\ 
    5 &  12 & 0.00 & 7.64 & 0.63 & 0.00 \\ 
    5 &  18 & 0.00 & 4.17 & 0.27 & 0.00 \\ 
    5 &  24 & 0.00 & 1.91 & 0.08 & 0.00 \\ \midrule
   10 &   5 & 0.00 & 151.81 & 1.59 & -0.01 \\ 
   10 &   7 & 0.01 & 344.41 & 2.59 & 0.00 \\ 
   10 &   9 & 0.02 & 585.85 & 3.49 & 0.00 \\ 
   10 &  14 & 0.04 & 1279.63 & 6.44 & -0.29 \\ 
   10 &  19 & 0.07 & 1845.95 & 10.23 & 0.02 \\ 
   10 &  26 & 0.04 & 1526.40 & 12.14 & 0.00 \\ 
   10 &  33 & 0.51 & 1084.73 & 11.37 & 0.07 \\ 
   10 &  44 & 0.02 & 686.13 & 7.84 & 0.00 \\ \midrule
   15 &   6 & 0.19 & 1930.95 & 37.63 & 0.14 \\ 
   15 &  12 & 14.02 & 4163.44 & 18.86 & 5.82 \\ 
   15 &  18 & 38.81 & 5119.94 & 33.08 & 19.47 \\ 
   15 &  24 & 39.67 & 4808.48 & 48.52 & 19.86 \\ 
   15 &  36 & 37.17 & 4508.40 & 101.38 & 15.97 \\ 
   \hline
\end{tabular}
\caption{Comparison of exact double CG and heuristic double CG methods for semi-log demand model instances from Section~\ref{subsec:results_discreteUcal_loglog_semilog}. \label{table:exact_vs_heuristic_semilog}} 
\end{table}

%%%%
% END OF EXACT VS HEURISTIC RESULTS
%%%%

%%%%
% START OF WIDTH EXPERIMENT RESULTS
%%%%

\subsection{Effect of parameter bounds on relative improvement in discrete $\Ucal$ experiments}
\label{subsec:extra_results_width_param}

In this section, we explore the effect of the width of the budget uncertainty set in our experiments in Section~\ref{subsec:results_discreteUcal_loglog_semilog}. Recall that in those experiments, we assumed a budget uncertainty set, where $\overline{\ub}$ and $\underline{\ub}$ denote the vectors of upper and lower bounds, respectively, for the uncertain parameter vector $\ub$, and these parameter vectors were defined as $\overline{\ub} = 1.3 \ub_0$, $\underline{\ub} = 0.7 \ub_0$. By setting the bounds in this way, we assume that each uncertain parameter can vary by +30\% or -30\% relative to its nominal value.

We now consider setting these bound vectors as 
\begin{align}
\overline{\ub} & = (1 + \delta) \ub_0, \\
\underline{\ub} & = (1 - \delta) \ub_0,
\end{align}
where $\delta \in (0,1)$ is the \emph{width parameter}. 

We consider the same log-log and semi-log demand model instances from Section~\ref{subsec:results_discreteUcal_loglog_semilog}, and allow for $\delta$ to vary in the set $\{0.1, 0.15, 0.2, 0.25, 0.3, 0.35, 0.4, 0.45, 0.5\}$. For simplicity, we restrict our focus to instances with $I = 5$, and allow the uncertainty set budget $\Gamma$ to vary in the same way as before, i.e., $\Gamma \in \{3,6,9,12,18,24\}$. For each instance, we solve the RRPO problem and DRPO problem using the same methods used in Section~\ref{subsec:results_discreteUcal_loglog_semilog} (i.e., the heuristic double column generation approach and the random improvement heuristic, respectively). We calculate the average value of the $\RI$ metric across all 24 instances for each $\Gamma$ and $\delta$.

\begin{figure}
\centering
\includegraphics[width=0.7\textwidth]{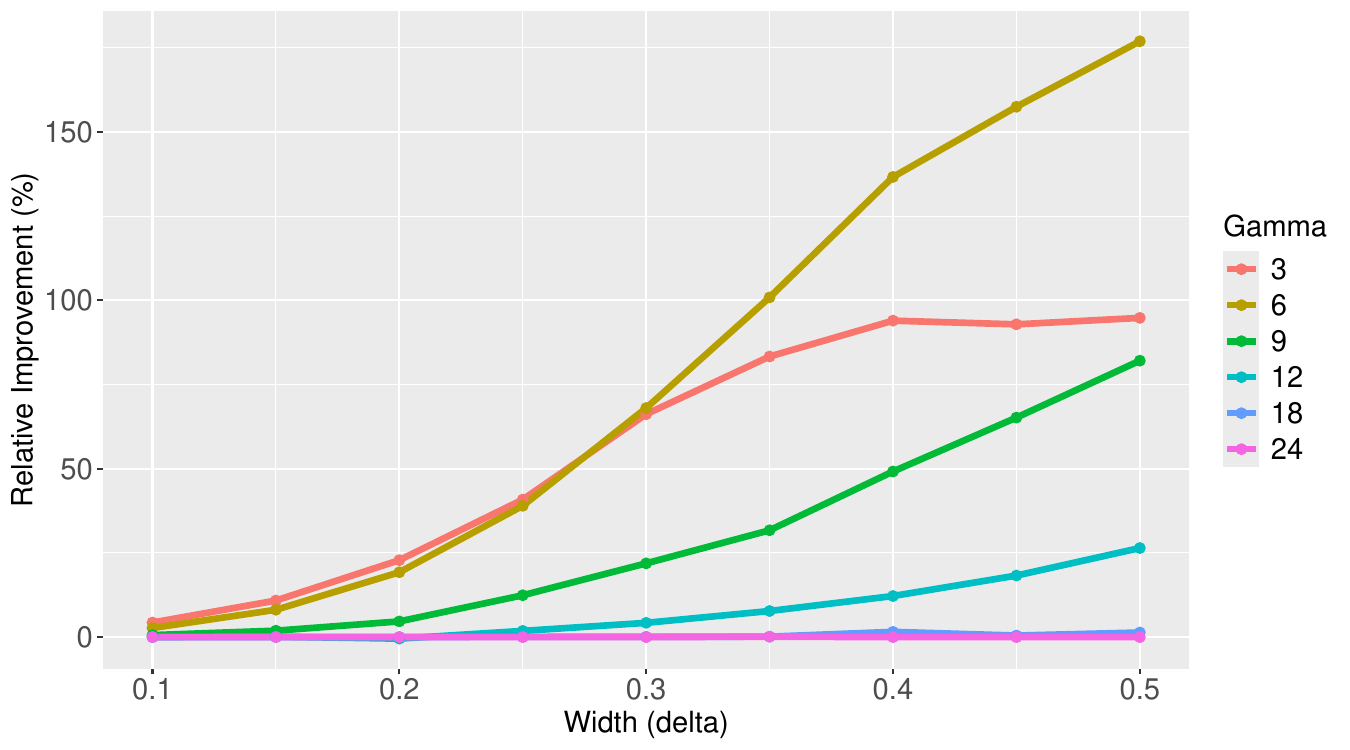}
\caption{Plot of $\RI$ versus $\delta$ for semi-log demand models. \label{figure:R1_semilog_RI_vs_width}}
\end{figure}

\begin{figure}
\centering
\includegraphics[width=0.7\textwidth]{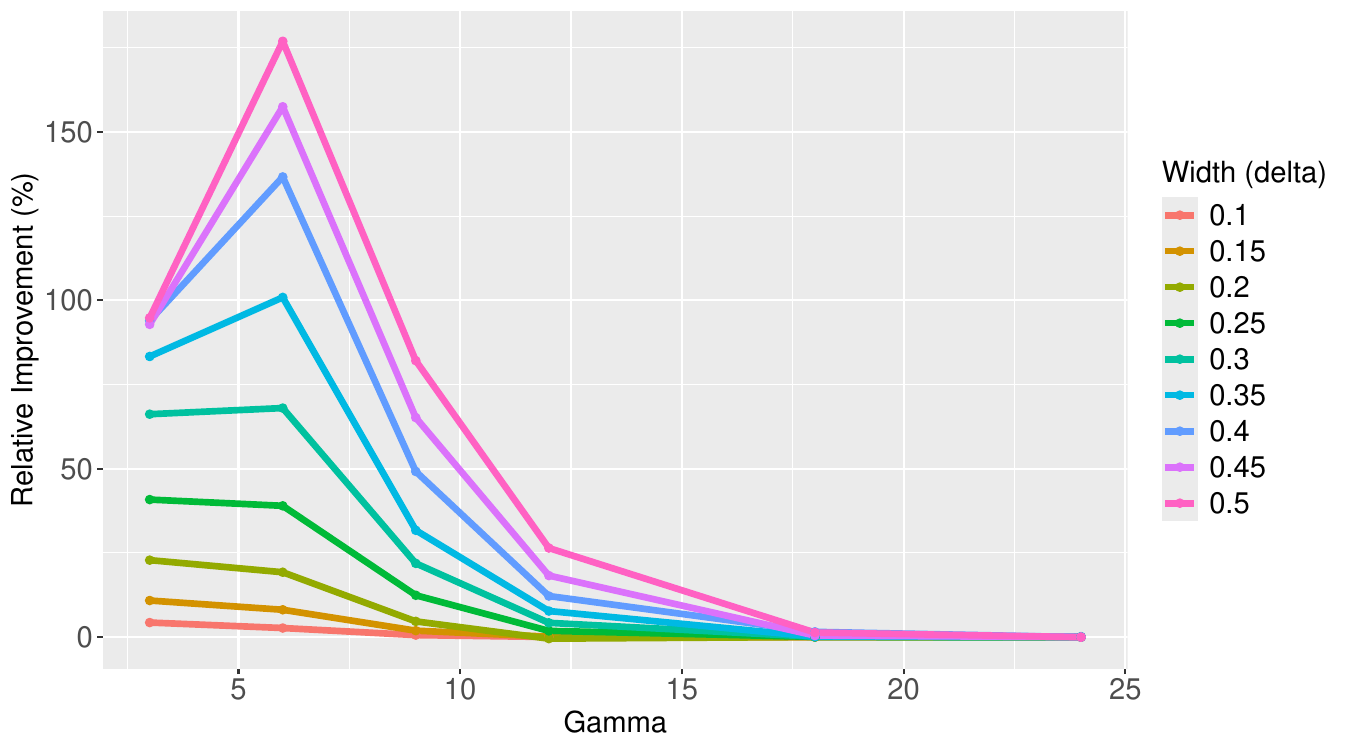}
\caption{Plot of $\RI$ versus $\Gamma$ for semi-log demand models. \label{figure:R1_semilog_RI_vs_Gamma}}
\end{figure}

For the semi-log demand model, Figure~\ref{figure:R1_semilog_RI_vs_width} plots the $\RI$ metric against $\delta$ for different values of $\Gamma$, while Figure~\ref{figure:R1_semilog_RI_vs_Gamma} plots the $\RI$ metric against $\Gamma$ for different values of $\delta$. In general, we can see that the relative improvement of the randomized robust price distribution over the deterministic robust price vector increases as the width parameter $\delta$ increases, suggesting that randomization becomes more beneficial when there is more uncertainty in each individual parameter. Interestingly, as observed in Section~\ref{subsec:results_discreteUcal_loglog_semilog}, the $\RI$ metric is initially increasing with $\Gamma$ but decreases after a certain point, and this peak seems to be the same across all values of $\delta$. Lastly, another interesting insight is that as $\Gamma$ increases, the rate at which the $\RI$ increases with $\delta$ is lower, and for sufficiently high $\Gamma$ values, the $\RI$ seems to stay very close to zero across all $\delta$ values.

\begin{figure}
\centering
\includegraphics[width=0.7\textwidth]{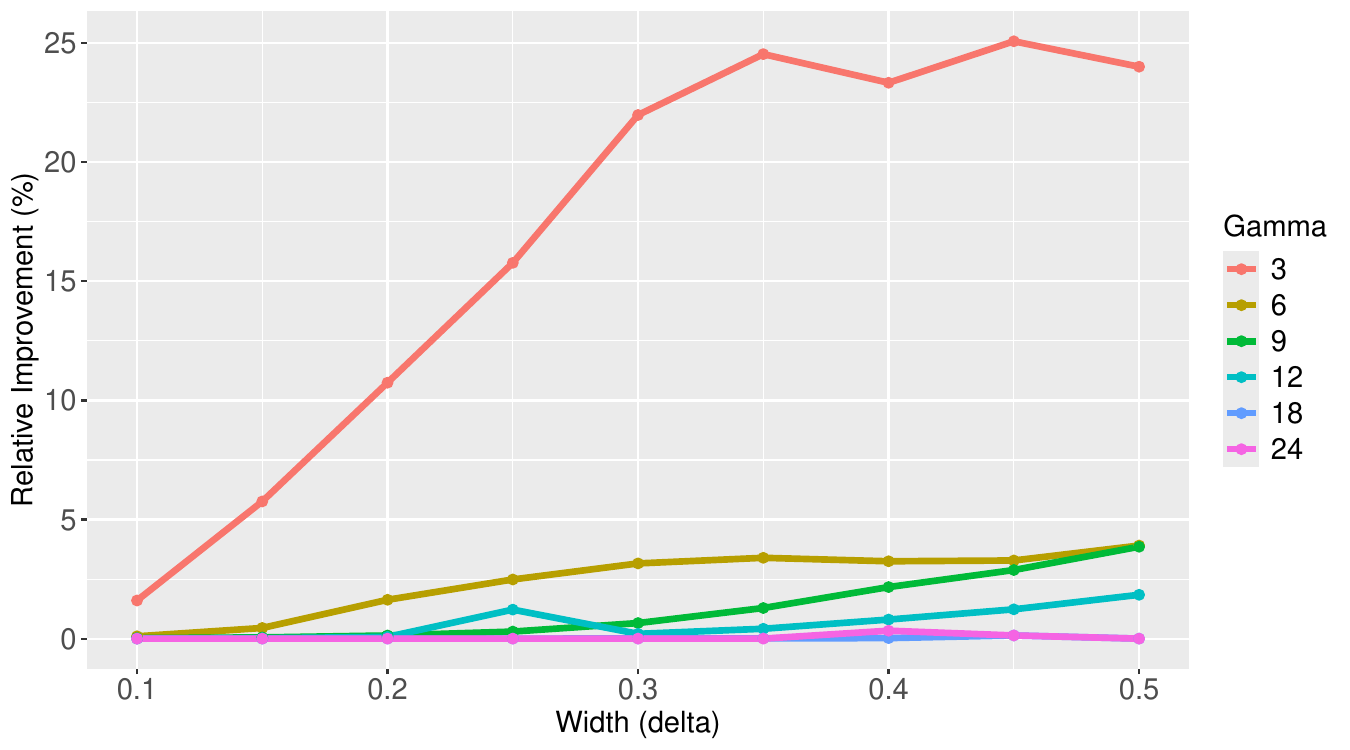}
\caption{Plot of $\RI$ versus $\delta$ for log-log demand models. \label{figure:R1_loglog_RI_vs_width}}
\end{figure}

\begin{figure}
\centering
\includegraphics[width=0.7\textwidth]{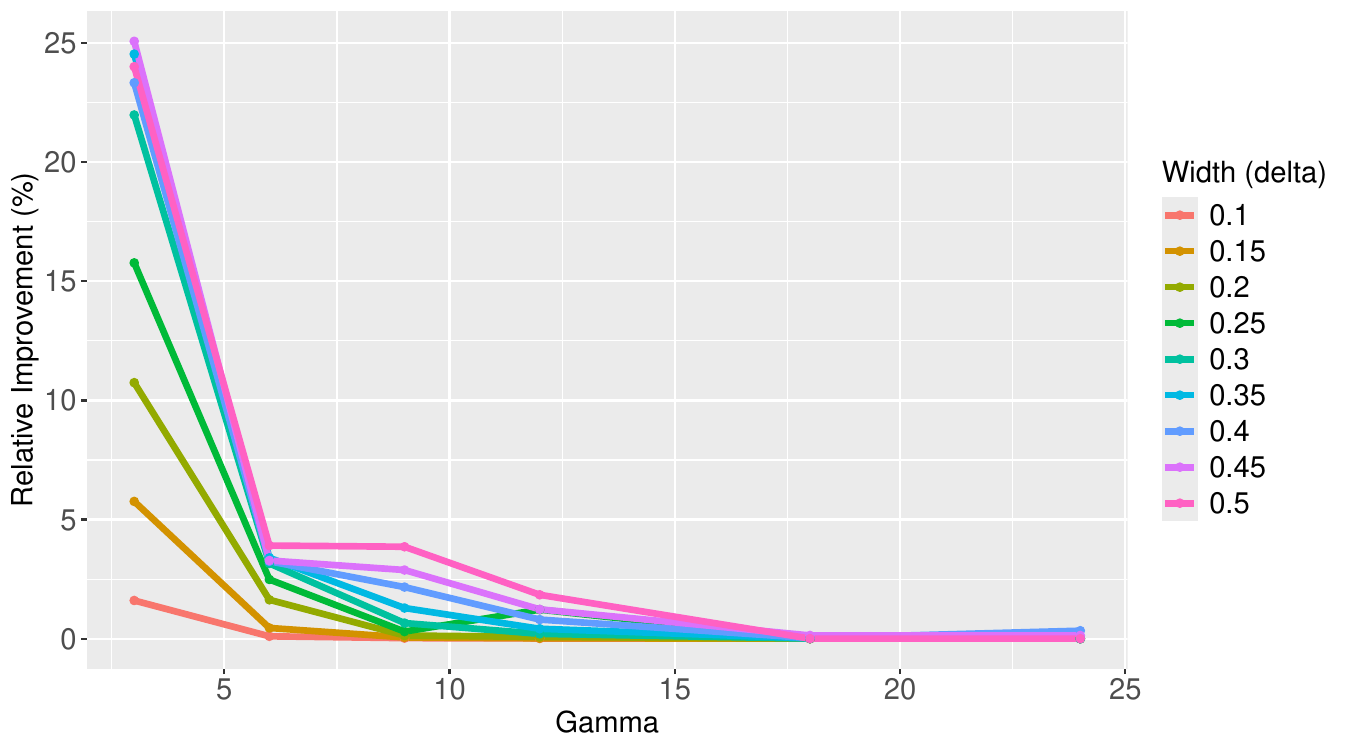}
\caption{Plot of $\RI$ versus $\Gamma$ for log-log demand models. \label{figure:R1_loglog_RI_vs_Gamma}}
\end{figure}

For the log-log demand model, Figure~\ref{figure:R1_loglog_RI_vs_width} plots the $\RI$ metric against $\delta$ for different values of $\Gamma$, and Figure~\ref{figure:R1_loglog_RI_vs_Gamma} plots the $\RI$ metric against $\Gamma$ for different values of $\delta$. Like in the semi-log case, the benefit of randomization over the deterministic robust solution does appear to increase as $\delta$ increases. The rate at which the $\RI$ increases becomes smaller when $\Gamma$ is higher, and for the highest values of $\Gamma$ ($\Gamma = 18$ or $\Gamma = 24$), changing $\delta$ (at least within the 0.1 to 0.5 range) does not seem to change the $\RI$.

%%%%
% END OF WIDTH EXPERIMENT RESULTS
%%%%

\subsection{Estimation results for \orangeJuice data set}
\label{subsec:extra_results_orangejuice_estimates}

Tables~\ref{table:R Estimation Results_1} and \ref{table:R Estimation Results_2} display the point estimates of $\alphab$, $\betab$ and $\gammab$ for the semi-log and log-log demand models for the \orangeJuice data set. 
\begin{table}[ht]
	\centering
	\small
	
		\begin{tabular}{cccccccccccc}
			\toprule
			& \multicolumn{11}{c}{Product} \\
			Parameters & 1 & 2 & 3 & 4 & 5 & 6 & 7 & 8 & 9 & 10 & 11 \\
			\midrule
			$\alphab_{\text{semi-log}}$ & 9.873 & 9.829 & 8.598 & 9.504 & 9.024 & 9.828 & 8.582 & 7.901 & 7.152 & 11.161 & 10.896  \\
			$\betab_{\text{semi-log}}$ & 1.0222 & 0.4581 & 1.2735 & 1.7888 & 1.3354 & 0.6507 & 1.6491 & 1.3945 & 2.0809 & 1.6290 & 0.0383  \\ \midrule
			$\alphab_{\text{log-log}}$ & 10.140 & 10.956 & 8.266 & 8.421 & 9.045 & 10.613 & 7.832 & 7.127 & 6.563 & 11.326 & 11.198  \\
			$\betab_{\text{log-log}}$ &2.7195 & 2.0410 & 3.3037 & 3.8855 & 2.9357 & 2.6101 & 3.6063 & 2.8209 & 3.9717 & 2.7942 & 0.1542 \\	  \bottomrule
	\end{tabular}
	\caption{Estimation results for $\alphab$ and $\betab$. \label{table:R Estimation Results_1} }
\end{table}

\begin{table}[ht]
	\centering
	\small
		\begin{tabular}{ lrrrrrrrrrrr} 
			\toprule
			$\gammab_{\text{semi-log}}$ & $j = 1$ & $j = 2$ & $j = 3$ & $j = 4$ & $j = 5$& $j = 6$ & $j = 7$& $j = 8$& $j = 9$& $j = 10$& $j = 11$  \\  \midrule
			$i = 1$ & --  & 0.0571 & 0.0813 & 0.0966 & 0.0193 & -0.0232 & 0.1305 & 0.1904 & 0.1490 & 0.0582 & 0.0815   \\
			$i = 2$ & 0.1384 & -- & 0.0041 & 0.0009 & 0.0204 & 0.0153 & 0.0090 & 0.1040 & -0.0023 & 0.0491 & 0.0394 \\ 
			$i = 3$ & 0.3386 & 0.0916 & -- & 0.1943 & 0.0702 & -0.0062 & 0.0051 & 0.0950 & -0.0310 & 0.0690 & 0.0950 \\ 
			$i = 4$ & 0.4313 & 0.0976 & -0.1112 & -- & 0.4089 & 0.3518 & 0.2085 & -0.0777 & -0.0352 & 0.0383 & -0.2290 \\   
			$i = 5$ & 0.1916 & 0.0490 & 0.3026 & 0.2966 & -- & -0.1538 & 0.1547 & -0.0314 & 0.1034 & 0.3338 & 0.0370 \\   
			$i = 6$ & 0.0211 & 0.0493 & -0.0194 & -0.0018 & 0.0888 & -- & 0.0340 & 0.0472 & -0.0167 & 0.0297 & 0.1119 \\   
			$i = 7$ & 0.2007 & 0.0388 & 0.0706 & 0.0672 & 0.3233 & 0.0837 & -- & 0.0377 & 0.2216 & -0.0504 & 0.1405 \\   
			$i = 8$ & 0.0117 & 0.0119 & 0.0932 & 0.0757 & 0.1023 & -0.0160 & 0.1345 & -- & 0.1372 & 0.2143 & 0.2699 \\   
			$i = 9$ & 0.0955 & 0.0373 & -0.0211 & 0.3651 & 0.4176 & 0.0358 & 0.2127 & 0.1462 & -- & 0.2337 & 0.1627 \\   
			$i = 10$ & 0.0412 & -0.3941 & 0.0764 & 0.4867 & 0.4810 & 0.0109 & -0.0814 & -0.1047 & 0.0878 & -- & 0.0274 \\  
			$i =  11$ & -0.0893 & -0.1587 & -0.1358 & -0.0252 & -0.0690 & 0.0079 & -0.0574 & -0.1117 & -0.1271 & 0.0809 & -- \\ \midrule
			$\gammab_{\text{log-log}}$ & $j = 1$ & $j = 2$ & $j = 3$ & $j = 4$ & $j = 5$ & $j = 6$ & $j = 7$ & $j = 8$ & $j = 9$& $j = 10$& $j = 11$  \\ \midrule
			$i = 1$ & --   & 0.2196 & 0.1631 & 0.2129 & 0.0646 & -0.0577 & 0.2576 & 0.3338 & 0.2494 & 0.0621 & 0.2939 \\   
			$i = 2$ & 0.3474 & -- & 0.0403 & 0.0004 & 0.0338 & 0.0492 & 0.0193 & 0.1879 & -0.0042 & 0.0739 & 0.1257 \\   
			$i = 3$ & 0.8673 & 0.5123 & -- & 0.4400 & 0.1482 & -0.0338 & 0.0527 & 0.1480 & -0.1001 & 0.1267 & 0.3683 \\   
			$i = 4$ & 1.1581 & 0.3822 & -0.2283 & -- & 0.8367 & 1.2659 & 0.4495 & -0.1569 & -0.0115 & 0.1003 & -0.7321 \\   
			$i = 5$ & 0.4624 & 0.2241 & 0.8344 & 0.6406 & -- & -0.6800 & 0.3223 & 0.0646 & 0.1426 & 0.5815 & 0.1782 \\   
			$i = 6$ & 0.0462 & 0.2424 & -0.0343 & -0.0173 & 0.2086 & -- & 0.0975 & 0.1187 & -0.0364 & 0.0561 & 0.4159 \\   
			$i = 7$ & 0.4644 & 0.2531 & 0.0971 & 0.1204 & 0.6997 & 0.2946 & -- & 0.1728 & 0.4912 & -0.0564 & 0.4497 \\   
			$i = 8$ & 0.0652 & 0.1430 & 0.1980 & 0.1587 & 0.2705 & -0.0988 & 0.3198 & -- & 0.3034 & 0.2992 & 0.9436 \\   
			$i = 9$ & 0.2971 & -0.1190 & -0.0216 & 0.7986 & 0.8825 & 0.3045 & 0.5869 & 0.1706 & 0.0& 0.3171 & 0.4223 \\
			$i = 10$ & 0.1406 & -1.7987 & 0.1061 & 1.0453 & 1.0852 & 0.0811 & -0.1213 & -0.1760 & 0.0454 & -- & -0.0371 \\   
			$i = 11$ & -0.2246 & -0.7519 & -0.3177 & -0.0483 & -0.1635 & -0.0570 & -0.1100 & -0.1780 & -0.2646 & 0.1341 & -- \\  \bottomrule
	\end{tabular}
	\caption{Estimation results for $\gammab$ for \orangeJuice data set. \label{table:R Estimation Results_2} }
\end{table}

\subsection{Performance results for \orangeJuice data set}
\label{subsec:extra_results_orangejuice_discreteUcal}

Tables~\ref{table:orangejuice_discreteU_semi_log} and \ref{table:orangejuice_discreteU_log_log} below compare the performance of the nominal, deterministic robust and randomized robust pricing solutions under a discrete budget uncertainty set for the \orangeJuice data set.

\begin{table}[ht]
	\centering
	\small
		\begin{tabular}{lllllllllll}
			\toprule
			$\Gamma$ & $t_{\RR}$ & $Z_{\RR}^*$ & $\Ebb[R(\pb^*_{\RR},\ub_0)]$ & $t_{\DR}$ & $\hat{Z}_{\DR}$ & RI(\%) & $R(\hat{\pb}_{\DR},\ub_0)$ & $t_{\Nom}$ & $Z_{\Nom}^*$ & $Z_{\Nom,\WC}$ \\ 
			\midrule
			5 & 14.64 & 162753.97 & 290939.28 & 225.36 & 102626.41 & 58.59 & 260321.17 & 0.81 & 590547.01 & 85304.36 \\    
			10 & 6.29 & 70401.48 & 404458.22 & 208.42 & 47969.46 & 46.76 & 350396.25 & -- & -- & 38815.61 \\    
			15 & 3.93 & 39567.50 & 349936.30 & 209.14 & 32757.43 & 20.79 & 334211.84 & -- & -- & 22798.44 \\    
			20 & 16.70 & 31438.76 & 328664.19 & 197.37 & 25348.77 & 24.02 & 299970.20 & -- & -- & 15940.16 \\ 
			\bottomrule
		\end{tabular}
	\caption{Results for \orangeJuice pricing problem with semi-log demand and discrete $\Ucal$. \label{table:orangejuice_discreteU_semi_log}}
\end{table}

\begin{table}[ht]
	\centering
	\small
		\begin{tabular}{lllllllllll}
			\toprule
			$\Gamma$ & $t_{\RR}$ & $Z_{\RR}^*$ & $\Ebb[R(\pb^*_{\RR},\ub_0)]$ & $t_{\DR}$ & $\hat{Z}_{\DR}$ & RI(\%) & $R(\hat{\pb}_{\DR},\ub_0)$ & $t_{\Nom}$ & $Z_{\Nom}^*$ & $Z_{\Nom,\WC}$ \\ 
			\midrule
			5 & 12.85 & 272399.89 & 605265.00 & 306.63 & 174478.12 & 56.12 & 811254.69 & 0.87 & 1110000.00 & 117186.58 \\    
			10 & 10.41 & 135761.31 & 750084.92 & 260.96 & 77297.42 & 75.63 & 896972.12 & -- & -- & 50458.15 \\    
			15 & 15.59 & 72930.45 & 761785.89 & 193.13 & 44914.20 & 62.38 & 896972.12 & -- & --& 27389.15 \\    
			20 & 8.56 & 45153.74 & 770505.32 & 190.76 & 27675.07 & 63.16 & 409330.70 & -- & --& 17502.32 \\ 
			\bottomrule
		\end{tabular}
	\caption{Results for \orangeJuice pricing problem with log-log demand and discrete $\Ucal$. \label{table:orangejuice_discreteU_log_log}}
\end{table}

\end{document}